\documentclass[11pt]{article}
\usepackage{a4}
\usepackage{ntheorem}
\usepackage{amssymb,latexsym}
\usepackage{citesort}
\usepackage{url}
\usepackage{eqname}
\usepackage{mhequ}
\usepackage{deleq}

\newcommand{\ztr}{{\mycal{K}}}

\newcommand{\mpsi}{\left(\begin{array}{c}\psi_0 \\ \psi_- \\ \psi_1 \\ \vdots \\
\psi_{n-1} \end{array}\right)}
\newcommand{\mphi}{\left(\begin{array}{c}\varphi_+ \\ \varphi_1 \\ \vdots \\
\varphi_{n-1} \end{array}\right)}
\newcommand{\dyw}{\underline{\mathrm{div}}\,}

\newcommand{\otimess}{\mathring\otimes}%
\newcommand{\stau}{y}%
\newcommand{\mcN}{{\mycal N}}%
\newcommand{\gbeta}{{\gamma}}%
\newcommand{\mcM}{{\mycal M}}%
\newcommand{\psize}{{\psi_0}}%
\newcommand{\alphahyp}{{\sigma}}

\newcommand{\newV}{\mycal V}
\newcommand{\newN}{\mycal O}
\newcommand{\newH}{\mycal N}

\newcommand{\fmale}{\Upsilon}
\newcommand{\Gamm}{\Gamma}

\newcommand{\tpsi}{{\tilde \psi}}
\newcommand{\tvarphi}{{\tilde \varphi}}

\newcommand{\cprim}{c_{\varphi}'}

\newcommand{\argw}{w}

\newcommand{\yhypini}{\lambda}

\newcommand{\Axp}{{\mycal A}^{+}_{\{x=0\}}}
\newcommand{\Axyp}{{\mycal A}^{+}_{\{\zlxy\}}}
\newcommand{\zlxy}{0\le x \le y}

\ifx\macrosloaded\relax\endinput\else\let\macrosloaded\relax\fi

\newcommand{\localmu}{\mu}%

\newcommand{\localkdelta}{p\delta}%

\newcommand{\ffT}[2]{\mcT^{#1,#2}_{\{\zlxy \},\infty}}
\newcommand{\fffT}[3]{\mcT^{#1,(#2;#3)}_{\{\zlxy \},\infty}}

\newcommand{\chleq}[1]{\cite[Equation~(#1)]{ChLengardnwe}}
\newcommand{\chl}{\cite{ChLengardnwe}}
\newcommand{\chcite}[1]{\cite[#1]{ChLengardnwe}}
\newcommand{\mcH}{{\mycal H}}
 \newcommand{\psAxf}[1]{\fAx{\dot \oplus_{i} x^{i\delta}\fzTi{#1-i\delta}}}%
\newcommand{\isAxf}[1]{\psAxf{#1}}%

\newcommand{\phg}{\mbox{\scriptsize \rm phg}}

\newcommand{\sCysF}{{\displaystyle \dot \oplus_n} F_n}

\newcommand{\mcU}{{\mycal U}}
\newcommand{\mcS}{{\mycal I}}
\newcommand{\mcA}{{\mycal A}}
\newcommand{\mcO}{{\mycal O}}
\newcommand{\thm}[1]{{Theorem~\ref{#1}}}

\newcommand{\zvarphi}{\mathring{\varphi}}

\newcommand{\bel}[1]{\begin{equation}\label{#1}}
\newcommand{\beal}[1]{\begin{eqnarray}\label{#1}}

\newcommand{\eean}{\end{eqnarray*}}
\newcommand{\bean}{\begin{eqnarray}\nonumber}
\newcommand{\eeal}[1]{\label{#1}\end{eqnarray}}
\newcommand{\mcF}{{\mycal F}}

\newcommand{\bcM}{\,\,\,\,\widetilde{\!\!\!\!\cM}}

\newtheorem{defi}{\sc Definition\rm}[section]

\newtheorem{Theorem}[defi]{\sc Theorem\rm}

\newtheorem{Proposition}[defi]{\sc Proposition\rm}

\newtheorem{Lemma}[defi]{\sc Lemma\rm}

\theorembodyfont{\upshape}
\newtheorem{Remark}[defi]{{\sc Remark}\rm}

\newtheorem{remark}[defi]{{\sc Remark}\rm}

{\catcode `\@=11 \global\let\AddToReset=\@addtoreset}
\AddToReset{equation}{section}
\renewcommand{\theequation}{\thesection.\arabic{equation}}

\newcommand{\be}{\begin{equation}}
\newcommand{\bea}{\begin{eqnarray}}
\newcommand{\eea}{\end{eqnarray}}
\newcommand{\beaa}{\begin{eqnarray*}}
\newcommand{\eeaa}{\end{eqnarray*}}
\newcommand{\bseq}{\begin{subeq}}
\newcommand{\eseq}{\end{subeq}}

\def \rectangle#1#2{\hbox{\vrule\vbox to #2
{\hrule\hbox to
#1{\hfil}\vfil\hrule}\vrule}}

\newcommand{\nn}{\nonumber}

\newcommand{\supp}{\mbox{\rm supp}} 
\newcommand{\eps}{\epsilon} 

\def\frac#1#2{{{#1}\over{#2}}}

\newcommand{\mxtau }{(M_{x_1-2\tau})}

\newcommand{\hyp}{{\mycal S}}

     \newcommand{\setZ}{{\mathord{\mathbb Z}}}

\newcommand{\Z}{\setZ}

\newcommand{\qed}{\hfill $\Box$\bigskip}
\newcommand{\proof}{\noindent {\sc Proof:\ }}
\newcommand{\eeq}{\end{equation}}
\newcommand{\ee}{\end{equation}}
\newcommand{\beqa}{\begin{eqnarray}}
\newcommand{\beqar}{\begin{deqarr}}
\newcommand{\eeqa}{\end{eqnarray}}
\newcommand{\eeqar}{\end{deqarr}}
\newcommand{\beqarl}[1]{\begin{deqarr}\label{#1}}
\newcommand{\eeqarl}[1]{{\label{#1}\end{deqarr}}}
\newcommand{\eeqarll}[2]{{\label{#1}\arrlabel{#2}\end{deqarr}}}
\newcommand{\beqan}{\begin{eqnarray*}}
\newcommand{\eeqan}{\end{eqnarray*}}
\newcommand{\ba}{\begin{array}}
\newcommand{\ea}{\end{array}}

\newcommand{\mcC}{{\mycal C}}

\newcommand{\cG}{{\mycal G}}

\newcommand{\Id}{\mbox{\rm Id}} 
\newcommand{\const}{\mbox{\rm const}} 

\newcommand{\hide}[1]{}

\newcommand{\HH}{{\mycal H}}

\newcommand{\pmu}{\partial_\mu}
\newcommand{\px}{\partial_x}

\newcommand{\err}{\zeta} 

\DeclareFontFamily{OT1}{rsfs}{}
\DeclareFontShape{OT1}{rsfs}{m}{n}{ <-7> rsfs5 <7-10> rsfs7 <10->
rsfs10}{} \DeclareMathAlphabet{\mycal}{OT1}{rsfs}{m}{n}
\def\scri{{\mycal I}}%
\def\scrip{\scri^{+}}%
\newcounter{mnotecount}[section]

\newcommand{\R}{\mathbb R}
\newcommand{\N}{\mathbb N}

\newcommand{\eq}[1]{(\ref{#1})}

\newcommand{\Eq}[1]{Equation~\eq{#1}}

\newcommand{\cM}{\mycal M}

\newcommand{\Ayd}{{\mycal A}^\delta_{\{y=0\}}}
\newcommand{\Axd}{{\mycal A}^\delta_{\{x=0\}}}
\newcommand{\Axyd}{{\mycal A}^\delta_{\{\zlxy \}}}
\newcommand{\Ay}{{\mycal A}_{\{y=0\}}}
\newcommand{\Ax}{{\mycal A}_{\{x=0\}}}
\newcommand{\Axy}{{\mycal A}_{\{\zlxy \}}}

\newcommand{\stsg}{{\mathfrak g}}
\newcommand{\tf}{\widetilde f}

\newcommand{\beqd}{\begin{deqarr}}
\newcommand{\eeqd}{\end{deqarr}}

\newcommand{\Cp}{{\mathtt A}}
\newcommand{\soln}[1]{{\mycal S}^{#1}}

\newcommand{\zmcAk}{\,\,\,\,\mathring{\!\!\!\!\mcA_k}{}}
\newcommand{\Cpdk}{{\Cp}^{\delta}_k}
\newcommand{\zCpdk}{{\mathring{\Cp}}^{\delta}_k}
\newcommand{\eCp}{{\Cp}^{\delta}_{k\pm}}
\newcommand{\zeCp}{{\mathring{\Cp}}^{\delta}_{k\pm}}
\newcommand{\mczT}{\,\,\mathring{\!\!\mycal F}}
\newcommand{\mczF}{\,\,\mathring{\!\!\mycal F}}
\newcommand{\mcT}{{\mycal T}}

\newcommand{\mcFa}{{\mcF}^{\alpha}}

\newcommand{\mcTB}{\,\widehat{\!\mycal F}}
\newcommand{\mcTaB}{{\mcTB}^{\alpha}}
\newcommand{\mcTinfB}{{\mcTB}_{\infty}}
\newcommand{\mczTa}{{\mczT}^{\alpha}}
\newcommand{\calk}{\mathbb Z}

\newcommand{\zpsi}{\mathring{\psi}}
\newcommand{\zphi}{\mathring{\varphi}}

\newcommand{\fT}[2]{\mcF_{\{\zlxy \},{#2}}^{#1}}
\newcommand{\fTi}[1]{\mcF_{\{\zlxy \},{\infty}}^{#1}}
\newcommand{\fzTi}[1]{\mczF_{\{\zlxy \},{\infty}}^{#1}}
\newcommand{\fCy}[2]{\mcC_{\{y=0\},{#2}}^{#1}}
\newcommand{\fCx}[2]{\mcC_{\{x=0\},{#2}}^{#1}}
\newcommand{\fCxy}[2]{\mcC_{\{\zlxy \},{#2}}^{#1}}
\newcommand{\fCxi}[1]{\mcC_{\{x=0\},\infty}^{#1}}
\newcommand{\cAxyd}{\mcA_{\{\zlxy \}}^\delta}
\newcommand{\fAxy}[1]{\mcA_{\{\zlxy \},#1}^\delta}
\newcommand{\fAx}[1]{\mcA_{\{x=0\},#1}^\delta}

\newcommand{\fAyd}{\mcA_{\{y=0\}}^\delta}
\newcommand{\cDapskmxy}{{\mycal C}_{\{\zlxy \},k}^{\alpha+1,\sigma}}
\newcommand{\cDapsk}{{\mycal C}_{\{\zlxy \},k}^{\alpha+1,\sigma}}
\newcommand{\cDask}{{\mycal C}_{\{\zlxy \},k}^{\alpha,\sigma}}
\newcommand{\cDazk}{{\mycal C}_{\{\zlxy \},k}^{\alpha,0}}
\newcommand{\cDak}{{\mycal C}_{\{\zlxy \},k}^\alpha}
\newcommand{\cDskp}{{\mycal C}_{\{\zlxy \},k+1}^\sigma}
\newcommand{\cDstwo}{{\mycal C}_{\{\zlxy \},2}^\sigma}
\newcommand{\cDsone}{{\mycal C}_{\{\zlxy \},1}^\sigma}
\newcommand{\cDsk}{{\mycal C}_{\{\zlxy \},k}^\sigma}
\newcommand{\cDsmk}{{\mycal C}_{\{\zlxy \},k}^{\sigma-1}}
\newcommand{\cDapspkmy}{{\mycal C}_{\{y=0\},k}^{\alpha+\sigma+1}}%
\newcommand{\cDapspkmx}{{\mycal C}_{\{x=0\},k}^{\alpha+\sigma+1}}%
\newcommand{\cDaspkmxy}{{\mycal C}_{\{\zlxy \},k}^{\alpha,\sigma+1}}%

\newcommand{\cDai}{{\mycal C}_{\{\zlxy \},\infty}^\alpha}

\newcommand{\cDapmone}{{\mycal C}_{\{\zlxy \},1}^{\min(\alpha+\delta-\epsilon,0)}}

\newcommand{\cDmz}{{\mycal C}_{\{\zlxy \},0}^{\min(\alpha+\delta-\epsilon,0)}}

\newcommand{\cDmonem}{{\mycal C}_{\{\zlxy \},1}^{\min(\alpha+\delta-\epsilon,-\epsilon)}}

\newcommand{\cDmim}{{\mycal C}_{\{\zlxy \},\infty}^{\min(\alpha+\delta-\epsilon,-\epsilon)}}
\newcommand{\cDmimnp}{{\mycal C}_{\{\zlxy \},\infty}^{-\epsilon}}
\newcommand{\cDmimn}%
{{\mycal C}_{\{\zlxy
\},\infty}^{\min(\alpha+\delta-\epsilon,(mq-p)\delta-\epsilon,-\epsilon)}}
\newcommand{\cDmimns}%
{{\mycal C}_{\{\zlxy
\},\infty}^{\min(\sigma+\delta-\epsilon,(mq-p)\delta-\epsilon,-\epsilon)}}

\newcommand{\cDzim}{{\mycal C}_{\{\zlxy \},\infty}^{-\epsilon}}
\newcommand{\cDzimne}{{\mycal C}_{\{\zlxy \},\infty}^{-1+\delta-\epsilon}}

\begin{document}

\title{Polyhomogeneous solutions of nonlinear wave equations
without corner conditions}
\author{Piotr T.\ Chru\'sciel\thanks{E-mail
\protect\url{Piotr.Chrusciel@lmpt.univ-tours.fr}, URL
\protect\url{www.phys.univ-tours.fr/}$\sim$\protect\url{piotr}}\\
LMPT, F\'ed\'eration de recherche Denis Poisson\\ Parc de Grandmont\\
F-37200 Tours, France \\[2ex] Szymon \L\c{e}ski \\
Centrum Fizyki Teoretycznej PAN, Al. Lotnik\'{o}w 32/46,\\ 02-668
Warsaw, Poland}

\maketitle

\begin{abstract} The study of Einstein equations leads naturally to Cauchy
  problems with initial data on hypersurfaces which closely resemble
  hyperboloids in Minkowski space-time, and with initial data with
  \emph{polyhomogeneous} asymptotics, that is, with asymptotic
  expansions in terms of powers of $\ln r$ and inverse powers of $r$.
  Such expansions also arise in the conformal method for analysing
  wave equations in odd space-time dimension.  In recent work~\chl\ it
  has been shown that for non-linear wave equations, or for wave maps,
  polyhomogeneous initial data lead to solutions which are also
  polyhomogeneous \emph{provided} that an infinite hierarchy of corner
  conditions holds. In this paper we show that the result is true
  regardless of corner conditions.
\end{abstract}

\tableofcontents

\section{Introduction}

The \emph{hyperboloidal
  Cauchy problem} provides one of the methods of studying global properties of solutions of
the vacuum Einstein equations. Here one prescribes Cauchy data on
a spacelike hypersurface with asymptotic behavior somewhat similar
to that of a hyperboloid in a Minkowski space-time. Such Cauchy
problems can be used to prove non-linear stability of Minkowski
space-time within a restricted class of initial data sets,
see~\cite{ChDelay,Friedrich:Pune,AndersonChruscielConformal} and
references therein.  Because the initial data for Einstein
equations satisfy the constraint equations, one faces the need to
consider \emph{polyhomogeneous} initial data in generic
situations~\cite{AndChDiss}, that is, initial data with asymptotic
expansions in terms of powers of $\ln r$ and inverse powers of
$r$, where $r$ is a luminosity parameter. Such initial data are
too singular at the conformal boundary at infinity to be handled
by the standard theory of hyperbolic equations.

As a first step towards handling that question, in this paper we
study the hyperboloidal Cauchy problem for simpler nonlinear wave
equations. Indeed, a closely related approach for understanding
the global behavior of solutions of wave equations on Minkowski
space-time is provided by the ``conformal method'', which proceeds
as follows: Let $({\mycal
  M},\stsg)$ be an $(n+1)$--dimensional space-time and let
\begin{equation}
  \label{C.1}
\tilde{\stsg} = \Omega^{2}\stsg \;.
\end{equation}
Let $\Box _h$ denote the wave operator associated with a
Lorentzian metric $h$,
$$\Box_h f= {1\over\sqrt{|\det h_{\rho\sigma}|}} \pmu (\sqrt{|\det
  h_{\alpha\beta}|}h^{\mu\nu} \partial_\nu f).$$
Recall that the scalar curvature $R =R(\stsg)$ of $\stsg $ is
related to the corresponding scalar curvature $\tilde{R} =
\tilde{R}(\tilde{\stsg })$ of $\tilde{\stsg }$ by the formula
\begin{equation}
\tilde{R}\Omega^2 = R -2n \left\{ {1\over
\Omega}\Box_{\stsg}\Omega +{n-3 \over
  2} {|\nabla\Omega|^2_\stsg  \over \Omega^2}\right\}\;.
\label{C.2}
\end{equation}
It then follows from (\ref{C.2}) that we have the identity
\begin{equation}
\Box_{\tilde{\stsg }} (\Omega^{-{n-1\over 2}}f)=
\Omega^{-{n+3\over 2}}\left( \Box_\stsg  f +{n-1 \over
4n}(\tilde{R}\Omega^2 -R)f\right)\;. \label{C.3}
\end{equation}
When $\stsg$ is the Minkowski metric $\eta$, a non-trivial
conformal factor $\Omega\sim 1/r$ can be chosen so that $\tilde
{\stsg}$ is also the Minkowski metric, leading to
\begin{equation}
\Box_{\eta} (\Omega^{-{n-1\over 2}}f)= \Omega^{-{n+3\over 2}}
\Box_\eta  f \;. \label{C.3b}
\end{equation}
Because the conformal boundary $\{\Omega=0\}$ corresponds to points
which are infinitely far away for the original metric $\stsg$,
this technique allows one to reduce global-in-time existence
problems to local ones; this has been exploited by various authors
\cite{ChBPisa,ChristodoulouCPAM,ChBGu,ChBglwa,ChBdeSitter,ChBglwa,ChBNou}
for wave equations on a fixed background space-time.

On a more modest level, the identity \eq{C.3} can be used as a
starting point for the analysis of the asymptotic behavior of
solutions of the scalar wave equation, as it reduces the problem to a
study of the rescaled equation \eq{C.3b} near the set $\Omega=0$.
There is, however, a difficulty that arises for non-linear equations
in \emph{even} space-dimension because of half-integer powers of
$\Omega$ in \eq{C.3b}. This introduces singular terms\footnote{It
  should be mentioned that for special non-linearities the singular
  terms do not arise. Indeed, this is the case for the equation $\Box
  f=H(f)$, where $H$ is a polynomial containing only odd powers of
  $f$. Similarly no singularities arise for the wave map equation with
  special targets~\cite{ChBGu,ChBglwa}. However, \emph{e.g.} wave maps
  will have such singular terms in general.} in the equations, so that
the usual theory of hyperbolic PDEs does not apply. This problem has
been studied in~\cite{ChLengardnwe}, where it was shown that solutions
for a class of semi-linear wave equations, or of the wave-map
equation, with smooth or polyhomogeneous initial data on hyperboloids
in Minkowski space-time, will have a complete asymptotic expansion in
terms of half integer powers of $\Omega\sim 1/r$ and of powers of $\ln
\Omega$, \emph{provided} that an infinite hierarchy of ``corner
conditions" is satisfied.  The object of this paper is to show that no
corner conditions are necessary for polyhomogeneity of solutions of
those equations. We also extend the results in~\cite{ChLengardnwe} to
general asymptotically flat space-times with smooth, or
polyhomogeneous, conformal completions. We expect to be able to prove
corresponding results for Einstein equations in a near future, see
also \cite{Lengard}.

Our main results, in Minkowski space-time, are the following (we
refer the reader to Appendix~\ref{Sphgf} for the definition of the
spaces involved; the hyperboloids $\mcH_\tau$ are defined as
$$\mcH_\tau=\{(t-\tau)^2-|\vec r|^2=1\}\;;$$
the coordinate $x$ in
which the polyhomogeneous expansions are carried out is $1/|\vec r|$):

\begin{Theorem}\label{T2phg}  Let $\delta =1$ in odd space dimensions,
  and let $\delta =1/2$ in even space dimensions. Consider the
  equation \be \Box_\stsg f = H(f)\;,\label{SE.1o} \ee on
  ${\R}^{n,1}$, $n\geq 2$, with smooth or polyhomogeneous initial data
  on a hyperboloid $\mcH_0$: \bel{indaHf} f|_{\mcH_0}\in x^{(n-1)/2}
  \Big(\Axd\cap L^\infty\Big)\;,\quad \partial_\tau f|_{\mcH_0}\in
  x^{(n-1)/2} \Axd\;. \ee
 Suppose that $H(f)$ is smooth in $f$ and has
  a zero of order $\ell$ at $f=0$, with $\ell$ satisfying \be \ell\geq
  \cases{ 4\; , & $n=2$ , \cr 3\; , & $n=3$ , \cr 2\; , & $n\geq 4$ .}
\label{condH} \ee
There exists $\tau^*>0$ and a solution of \eq{SE.1o}-\eq{indaHf}
defined on $\cup_{\tau\in[0,\tau_*]}\mcH_\tau$ such that for every
$\epsilon>0$ we have
$$f|_{\cup_{\tau\in[\epsilon,\tau_*]}\mcH_\tau}\in x^{(n-1)/2}
\Big(\Axd\cap L^\infty\Big)\;.$$
\end{Theorem}

Theorem~\ref{T2phg} is a special case of Theorem~\ref{T2phgg}
below where asymptotically vacuum space-times with polyhomogeneous
$C^1$ conformal completions at null infinity satisfying  mild
restrictions on the scalar curvatures are considered, with general
polyhomogeneous initial data hypersurfaces, and with somewhat more
general non-linearities.

\begin{Theorem}\label{Twavemap}
  Let $\delta $ be as in Theorem~\ref{T2phg}. Consider the wave map
  equation on ${\R}^{n,1}$, $n\geq 2$. For any smooth or polyhomogeneous
  initial data on a hyperboloid $\mcH_0$,
\begin{eqnarray*}
 &
  x^{-(n-1)/2} f^a|_{\mcH_0}\in \Axd\cap L^\infty
 \;,
 &
 \\
 &
 \quad \partial_\tau (x^{-(n-1)/2} f^a)|_{\mcH_0}
  \in \left\{%
\begin{array}{ll}
       \Axd\cap
L^\infty\;, &  n=
  2\;,
  \\
  \Axd\;, &  n\ge
  3\;,\\
\end{array}%
\right.
 &
  \end{eqnarray*}
  there exists $\tau^*>0$, and a solution $f$ defined on
  $\cup_{\tau\in[0,\tau_*]}\mcH_\tau$, such that for every
  $\epsilon>0$ we have
  $$f|_{\cup_{\tau\in[\epsilon,\tau_*]}\mcH_\tau}\in x^{(n-1)/2} \Big(\Axd\cap
  L^\infty\Big)\;.$$
\end{Theorem}

Theorem~\ref{Twavemap} is a special case of
Theorem~\ref{Twavemapg} below.

This paper is organised as follows. In Appendix~\ref{SA} the reader
will find a description of our notation and conventions, as well as
the definitions of the function spaces used, together with some
auxiliary results needed in the body of the paper. In
Section~\ref{Sweci} we show how to apply our main results, derived in
later sections, to the nonlinear wave equation, or to the wave map
equation, on manifolds with conformal completions at null infinity.
The discussion in that section motivates the hypotheses of our
polyhomogeneity theorems in Section~\ref{SsCp}, which is the
cornerstone of the paper. There we prove polyhomogeneity of solutions
of a class of systems of first order PDE's. In Appendix~\ref{SAB} we
give a considerably simpler proof of our main polyhomogeneity theorem
for linear systems, under the supplementary assumption that all
coefficients in the equation are smooth and bounded.

We acknowledge useful discussions with Olivier Lengard at an early
stage of this work: he suggested a proof of the result in a simple
case, which was at the origin of this paper.

\section{Wave equations near conformal infinity}
 \label{Sweci}
 In this section we verify that those results
 of~\cite{ChLengardnwe}
 which have been proved for the Minkowski metric can be extended to
 space-times with completions at null infinity. We allow
 conformal completions with polyhomogeneous metrics, under a condition
 on the scalar curvature near the conformal boundary.  (In fact,
 metrics in weighted Sobolev spaces suffice for existence of solutions
 with weighted Sobolev regularity.) We also show how to write
 quasi-linear wave equations in a form to which our polyhomogeneity
 arguments in the next sections apply.

 For convenience of typesetting we will denote by $g$ the metric $\tilde{\stsg}$
 from the introduction.  We assume, as usual (\emph{cf.,
   e.g.},~\cite{Friedrich:Pune}), that the function $\Omega$ in
 \eq{C.1} vanishes precisely at a null hypersurface which we denote by $\newH$, with
 $d\Omega$ nowhere vanishing on $\newH$. In the construction that
 follows we assume that $\Omega$ is either smooth or $C^1$ and
 polyhomogeneous; equivalently, $\newH$ is a smooth or polyhomogeneous
 $C^1$ hypersurface.

 Near $\newH$ we can always find,
locally, a convenient coordinate system $(\stau ,x,v^A)$ so that
$\newH$ is given by the equation $x=0$, with $\partial_\stau $
tangent to the generators of $\newH$, as follows
(compare~\cite{VinceJimcompactCauchyCMP}; this construction can also
be carried out for $\Axd\cap L^\infty$-polyhomogeneous metrics,
see~\cite[Appendix~B]{ChMS}): Let $\newN \subset
\newH$ be any $(n-1)$--dimensional submanifold of $\newH$,
transverse to the null generators of $\newH$. Let $v^A$ be any
local coordinate system on $\newN $, and let $\ell|_\newN $ be any
field of null vectors, defined on $\newN $, tangent to the
generators of $\newH$. Solving the equation $\nabla_\ell \ell =
0$, with initial values $\ell|_\newN $ on $\newN $, one obtains a
null vector field $\ell$ defined on a $\newH$--neighborhood
$\newV\subset \newH $ of $\newN $, tangent to the generators of
$\newH$. One can extend $v^A$ to $\newV $ by solving the equation
$\ell(v^A)=0$. The function $\stau|_\newH$ is defined by solving
the equation $\ell(\stau)=1$ with initial value $\stau|_\newN =0$.
Passing to a subset of $\newV $ if necessary, this defines a
global coordinate system $(\stau,v^A)$ on $\newV $. By
construction we have $\ell=\partial_\stau$ on $\newV $, in
particular $g_{\stau\stau}=0$ on $\newV $. This justifies the
first equation in \eq{cangnn2} below. Further, $\ell$ is normal to
$\newH$ because $\newH$ is a null surface, which implies $g_{\stau
A}=0$ on $\newV $. This justifies the second equation in
\eq{cangnn2} below.

Let, next, $\bar \ell|_\newV $ be a field of null vectors on
$\newV $ defined uniquely by the conditions \bel{inicondbarl} g(\bar
\ell|_\newV , \ell)=\frac 12\;,\quad g(\bar \ell|_\newV ,
\partial_{A})=0\;.\ee
The first equation implies that $\bar \ell|_\newV $  is everywhere
transverse to $\newV $. Then we define $\bar \ell$  in a
space-time neighborhood $\mcU\subset \mcM$ of $\newV $ by solving
the geodesic equation $\nabla_{\bar \ell} \bar \ell = 0$ with
initial value $\bar \ell|_{\newV }$ at $\newV $. The coordinates
$(\stau,v^A)$ are extended to $\mcU$ by solving the equations
$\bar \ell (\stau)=\bar \ell(v^A)=0$, and the coordinate $x$ is
defined by solving the equation $\bar \ell(x)=1$, with initial
value $x=0$ at $\newV $. Passing to a subset of $\mcU$ if
necessary, this defines a global coordinate system $(\stau,x,v^A)$
on $\mcU$.

 By
construction we have
\bel{barelx} \bar\ell =\partial_x\;,\ee hence $\partial_x$ is a null, geodesic,
vector field on $\mcU$. In particular
$$g_{xx}\equiv g(\partial_x,\partial_x)=0\;.$$
Let $(z^a)=(x,v^A)$, and note that
 \beaa\bar \ell\Big( g(\bar \ell,
\partial_a)\Big)&=& g(\bar \ell, \nabla_{\bar \ell} \partial_a) =g(\bar
\ell, \nabla_{\partial_x} \partial_a)= g(\bar \ell,
\nabla_{\partial_a}
\partial_x)\\ &= &g(\bar \ell,
\nabla_{\partial_a} \bar \ell) =\frac 12
\partial_a\Big(g(\bar \ell, \bar \ell)\Big)=0\;.\eeaa
This shows that the components $g_{xa}$ of the metric are
$x$--independent. On $\newN$ we have $g_{x\stau}=1/2$ and
$g_{xA}=0$ by \eq{inicondbarl}, leading to the following form of
the metric
 \bel{cangnn} g=  dx \otimess  d\stau+
\chi d\stau \otimess d\stau + 2 \gbeta \otimess d\stau  + \mu\;,\ee
with \bel{cangnn2} \chi|_{x=0}=0\;,\quad\gbeta|_{x=0}=0\;.\ee The
symbol $\otimess$ denotes a symmetrized tensor product:
$$\alpha\otimess \gbeta = \frac 12\Big(\alpha \otimes \gbeta + \gbeta \otimes \alpha\Big)\; .$$
Furthermore, $\gbeta=\gbeta_Adv^A$ is a one-form field, and
$\mu=\mu_{AB}\,dv^A\otimess dv^B$ is a symmetric tensor field.  By
inspection of \eq{cangnn}, the tensor field induced from $g$ on the
surfaces $\stau=\const$ has signature $(0,+,\ldots,+)$, thus these
are null hypersurface.

An example of the coordinate system above is obtained by taking
$\newH$ to be the light-cone of the origin in $(n+1)$--dimensional
Minkowski space-time, with $x=r-t$, $y=t+r$, then the Minkowski
metric $\eta$ takes the form
$$\eta = -dt\otimes dt+dr\otimes dr+r^2d\Omega^ 2 = dx\otimess dy +r^2d\Omega^ 2\;,$$
so that $\chi\equiv\gbeta\equiv 0$.

In asymptotically flat solutions of asymptotically vacuum Einstein
equations, a slight variation of the construction above leads to the
\emph{Bondi coordinates} near Scri~\cite{Tamburino:Winicour}, with
the metric taking instead the form \bel{cangnnB} \tilde
g_B=e^{2\beta} dx \otimess d\stau+ \chi d\stau \otimess d\stau + 2
\gbeta \otimess d\stau + \mu\;,\ee for some function $\beta$ (and,
in general, different $\chi$ and $\gbeta$). (Here $y$ corresponds to
the Bondi retarded time $u$, and $x=1/2r$ is half the inverse of the
luminosity distance $r$.) In $3+1$ dimensions, for smoothly
compactifiable metrics, the Einstein equations imply, for matter
fields decaying sufficiently fast, that $\beta=O(x^2)$ as well as
\bel{Bonddec} \chi = O(x^2)\;,\quad \gbeta_A = O(x^2)\;,\ee with
derivatives behaving in the obvious way. One can go from
\eq{cangnnB} to \eq{cangnn} by redefining $x$, without changing the
remaining variables, in a way which will preserve \eq{Bonddec}.
\Eq{Bonddec} remains valid for asymptotically vacuum metrics which,
after conformal rescaling, are polyhomogeneous and $C^1$
(see~\cite[Section~6]{ChMS} or \cite[Appendix~C.1.2]{CJK}), while
for general $\Axd\cap L^\infty$--polyhomogeneous asymptotically
vacuum metrics one has~\cite[Equations~(2.15)-(2.19) with
$H=X^a=0$]{ChMS} the asymptotic behaviors $\beta=O(x^2\ln^Nx)$ and
\bel{Bonddec2} \chi = O(x^2)\;,\quad \gbeta_A = O(x^2\ln^Nx)\;,\ee
for some $N$. Here ``asymptotically vacuum'' requires, for
polyhomogeneous metrics, that the components of the energy-momentum
tensor in asymptotically Minkowskian coordinates satisfy
(see~\cite[end of Section~2]{ChMS}) \bel{asvac}
T_{\mu\nu}=o(r^{-2})\;.\ee

 Returning to \eq{cangnn}, we have $$\det
g = - \frac 14 \det \mu\;,$$ which shows that $\mu$ is a non-degenerate
$(n-1)\times (n-1)$ tensor field. It is simple to check that the
inverse metric $g^\sharp=g^{\alpha\gbeta}\partial_\alpha\otimess
\partial_\gbeta$ is given by the formula
\bean  g^\sharp&=&
{4(-\chi+ |\gbeta|_\mu^2)}
\partial_x \otimess
\partial_x + 4\partial_x  \otimess \partial_\stau
  - 4\gbeta^\sharp \otimess \partial_x + \mu^\sharp
  \\
  &=&
 4\partial_x  \otimess \Big(\partial_\stau + {(-\chi+ |\gbeta|_\mu^2)}
\partial_x -\gbeta^\sharp \Big)+ \mu^\sharp\;,\eeal{cangnni}
  with $\mu^\sharp= \mu^{AB}\partial_A \otimess \partial_B$, where $\mu^{AB}$
   is the matrix inverse to $\mu_{AB}$, $\gbeta^\sharp = \mu^{AB}\gbeta_A\partial_B$,
   and
   $|\gbeta|_\mu^2=\mu^\sharp(\gbeta,\gbeta)=\mu^{AB}\gbeta_A\gbeta_B$. We
note
   $$g(\nabla \stau, \nabla \stau)=g^{\stau \stau}=0\;,$$
which confirms the null character of the level sets of $\stau$, and
can also be used for an alternative justification of the fact that
the integral curves of
$$\nabla \stau=g^{\alpha\gbeta}\partial_\alpha \stau
\partial_\gbeta=g^{\stau
\gbeta}\partial_\gbeta=2\partial_x $$ are null geodesics.

Consider a system of $N$ second order PDE's for an $m$--components
field $u$, with diagonal principal part $(g^{\alpha\gbeta}\otimes
\Id)\partial_\alpha\partial_\gbeta$, where $\Id$ is the $m\times
m$ identity matrix: \bel{geneqsys}
g^{\alpha\gbeta}\partial_\alpha\partial_\gbeta u =
F(x^\mu,u,\partial u)\;.\ee (Recall that the metric $g$ here
should be thought of as the metric $\tilde{\stsg}$ of the
introduction, with the field $u$ being an appropriately rescaled
equivalent of some field defined on the original space-time.)  We
start by checking whether such a system can be written, locally,
in the form considered in Section~\ref{SsCp}
below. For this we let $e_A=e_A{}^C\partial_C$, $A=1,\ldots,n-1$
be a local ON basis for $\mu$, and set \beal{defred1}
&e_-=2\Big(\partial_\stau + {(-\chi+
  |\gbeta|_\mu^2)} \partial_x -\gbeta^\sharp\Big)\;,\quad
e_+=-2\partial_x\;,& \\ &\psize = u\;, \quad
\psi_A=\varphi_A=e_A(u)\;,\quad \varphi_+=e_+( u)\;,\quad
\psi_-=e_-(u)\;,\label{defred2} & \\ & \psi=(\psize
,\psi_-,\psi_A)\;,\quad \varphi=(\varphi_+,\varphi_A)\;.
\eeal{defred3} (Factors of two and signs  in \eq{defred1} have
been introduced for consistency of notation with
\cite{ChLengardnwe}.)  Then \beaa
g^{\mu\nu}\partial_\mu\partial_\nu u & =& 2e_-(\partial_x u) +
\sum_A\Big( e_A (e_A(u))-e_A(e_A{}^C)\partial_C u \Big)
\\
& = & 2\partial_x (e_-(u)) +2[e_-,\partial_x]u+ \sum_A \Big( e_A
(e_A(u))-e_A(e_A{}^C)f_C{}^B e_B u \Big) \;,\eeaa where $f_C{}^B$
denotes the matrix inverse to $e_A{}^B$, and \eq{geneqsys} implies the
following system of equations for $(\psi,\varphi)$: \bel{symhyp1}
\begin{array}{lll}\partial_x \psize & +\frac 12 \varphi_+&= 0\;,
  \\
  \partial_x \psi_-&- [\partial_x,e_-] u &= -\frac 12\sum_A ( e_A
  \varphi_A-e_A(e_A{}^C)f_C{}^B \psi_B ) + \frac
  12F(x^\mu,\psi,\varphi)\;,
  \\
  \partial_x \psi_A &- [\partial_x,e_A] u &= -\frac 12 e_A
  \varphi_+\;,
  \\
  \partial_y \varphi_+ &&= -\Big({(-\chi+ |\gbeta|_\mu^2)}
  \partial_x -\gbeta^\sharp\Big) \varphi_+
  \\
  \nonumber
 && \phantom{=}+ \frac 12\sum_A ( e_A \psi_A-e_A(e_A{}^C)f_C{}^B \psi_B  ) -\frac 12F(x^\mu,\psi,\varphi)\;, 
 \\
 \partial_y \varphi_A&- \frac 12 [e_-,e_A]u &= -\Big({(-\chi+
|\gbeta|_\mu^2)}
\partial_x -\gbeta^\sharp\Big) \varphi_A +\frac 12e_A\psi_-\;.
 \end{array}\ee 
 It is understood that the commutator terms have to be expressed as
 linear functions of $\varphi$ and $\psi$, and that those
 substitutions have also
 been done when writing $F(x^\mu,u,\partial u)$ as
 $F(x^\mu,\psi,\varphi)$. It follows that all the requirements
 concerning the principal part of \eq{geneqsys}, set forth in
 Section~\ref{SsCp} will be satisfied for polyhomogeneous, or smooth, metrics
 if there exists $m$ such that \bel{decrestr}
 \chi=O(x^{1+\delta}\ln^mx)\;, \quad \gbeta_A=O(x^{1+\delta}\ln^m
 x)\;.\ee In situations of main physical interest with $\newH=\scrip$
 the dimension is $3+1$, the space-times are asymptotically vacuum,
 the appropriate choice of $\delta$ is one, and \Eq{Bonddec2} shows
 that \eq{decrestr} holds.

 In situations where \eq{decrestr} does not hold one can proceed as
 follows. Suppose, first, that the metric $g$ is smooth.  One can then
 introduce new coordinates $(\hat x,\hat y, \hat v^A)$ by setting
 $\hat y=y$, and letting $(\hat x,\hat v^A)$ be defined as solutions
 of the equations \bel{hatvar0} e_-(\hat v^A)= e_-(\hat x) = 0\;, \ee
 with initial values $\hat v^A|_{y=0}=v^A$, $\hat x|_{y=0}= x$.  In
 the new coordinate system we have $e_-=\partial_{\hat y}$;
 equivalently \bel{hatvar}\frac {\partial y}{\partial \hat
   y}=1\;,\quad \frac {\partial x}{\partial \hat y}=-\chi+
 |\gbeta|_\mu^2\;,\quad \frac {\partial v^A}{\partial \hat
   y}=-\mu^{AB}\gamma_B\;.  \ee This gets rid of all derivatives of
 $\varphi$ in the right-hand-sides of the last two lines of
 \eq{symhyp1}, when those equations are rewritten in terms of the
 hatted coordinates; \emph{e.g.}, the last equation in \eq{symhyp1}
 becomes
$$ \partial_{\hat y} \varphi_A- \frac 12 [e_-,e_A]u = \frac 12e_A\psi_-\;.
$$
 (Note that we do not change the
definitions \eq{defred1}-\eq{defred3}, which still
 use the coordinates $(x,y,v^A)$ constructed above.)
One can  check that the remaining derivatives in \eq{symhyp1} will
change  in a way compatible with our hypotheses in
Section~\ref{SsCp}. In any case this follows from the analysis of
$C_1\cap \Axd$ metrics, to which we pass now: We need \eq{hatvar}
to be satisfied to first order of $x$:
\bel{hatvar2}
\frac {\partial x}{\partial \hat y}=-\chi+
|\gbeta|_\mu^2+O(x^{1+\delta}\ln^mx)\;,\quad \frac {\partial
v^A}{\partial \hat
y}=-\mu^{AB}\gamma_B+O(x^{1+\delta}\ln^mx)\;.
 \ee
  We have
 $$-\chi+
|\gbeta|_\mu^2\;,\  \mu^{AB}\gamma_B \in  xC_\infty \cap
x^{1+\delta} \Axd\;,
 $$
  therefore there exist $\fmale, \Gamm^A$, smooth functions
of $(v^A, y)$, such that
$$
(-\chi+ |\gbeta|_\mu^2)(x,v^A,y) = \fmale(v^A,y)x +
O(x^{1+\delta}\ln^mx)\;,$$
$$-(\mu^{AB}\gamma_B)(x,v^A,y) = \Gamm^A(v^A,y)x+
O(x^{1+\delta}\ln^mx)\;.$$
We make the ansatz \bel{chhatcrd} \hat x =
x + x\xi(v^A,y)\;, \quad \hat v^A = v^A + x\eta^A(v^B,y)\;, \quad \hat
y = y\;, \ee and \eq{hatvar2} leads to the following equations for
$\xi, \eta^A$: \bel{hateqn} \fmale(1+\xi) +
\frac{\partial\xi}{\partial y} =0\;,\quad \Gamm^A +
\frac{\partial\eta^A}{\partial y} + \eta^A \fmale = 0\;.\ee We impose
$\xi(v^A,0)=0=\eta^A(v^B,0)$. Conversely, smooth solutions to
\eq{hateqn} with those initial values give (at least for $x,y$ close
to $0$)
$$\frac {\partial x}{\partial \hat y} = \fmale(v^A,y)x + O(x^2)\;,
\quad \frac {\partial v^A}{\partial \hat y} = \Gamm^A(v^A,y)x +
O(x^2)\;,$$ and we end up with the following transformation rules
for the derivatives:
\begin{eqnarray*}
\partial_x & = &(1+\xi)\partial_{\hat x} + \eta^A\partial_{\hat
v^A}\;,\\
\partial_{v^A} & = &\partial_{\hat v^A} + x\eta^B{}_{,A}\partial_{\hat
v^B}+x\xi_{,A}\partial_{\hat x}\;,\\
\frac 12e_{-} & = & \partial_{\hat y} + O(x^{1+\delta}\ln^m x)\;.
\end{eqnarray*}
This shows that, in the hatted coordinates, for all
$\Axd$-polyhomogeneous metrics, \eq{symhyp1} takes a form which is
compatible with our hypotheses in Section~\ref{SsCp}.

We consider now a scalar field $u$ satisfying a non-linear
massless wave equation in the physical space-time. What has been
said shows that the second-derivative terms in the equation
satisfy the requirements of Section~\ref{SsCp}. The zero order
terms arising from the supplementary curvature term in \eq{C.3}
will be compatible with the requirements of Theorem~\ref{Tmnl} or
Theorem~\ref{Tlemme1n} if, \emph{e.g.},\begin{equation}
  \Omega^{-{n+3\over 2}}\left(R(\tilde{\stsg})\Omega^2
    -R(\stsg)\right)\in \Axyd \quad \Longleftrightarrow\quad
  R(\tilde{\stsg})\Omega^2 -R(\stsg)\in x^{{n+3\over 2}}\Axyd\;.
  \label{C.3n}
\end{equation}
One can always achieve $R(\tilde{\stsg})=0$ by a smooth redefinition
of the conformal factor $\Omega$, when $\tilde \stsg$ is smooth, by
solving the Yamabe equation. Similarly if the metric is $C_1\cap \Axd$
and $R(\tilde \stsg)$ satisfies \bel{Rtgfc} R(\tilde \stsg) \in \Axd
\ee (in fact $R(\tilde \stsg) \in \Axyd$ suffices), then our results
below show that this redefinition can again be made, leading to a new
metric of
$\Axyd$--differentiability class.\footnote{A sufficient condition for
  $R(\tilde \stsg) \in \Axyd$ is that $\tilde \stsg$ is of
  $(C_\infty+x^{2}\Axd)$--differentiability class, but this
  restriction is certainly not necessary; note that for $\delta=1$ we
  have $C_\infty+x^{2}\Axd=C_1\cap \Axd$. It is likely that the
  condition $g\in C_1\cap \Axd$ suffices for polyhomogeneity of solutions
  for all $\delta$, see Remark~\ref{Rldc} below.}
Then \eq{C.3n} becomes a restriction on the rate of fall-off of
scalar curvature of the space-time metric $\stsg$ near null
conformal infinity $\scrip=\mcN$, in particular \eq{C.3n} will be
satisfied when $\stsg$ is vacuum. Note also that \eq{Bonddec2} is
preserved by any conformal transformation with bounded conformal
factor, so that such a redefinition of $\tilde \stsg$ will not
affect that equation.

Next, still for a scalar field, in order to find the first order
terms which arise from $\Box_{\tilde \stsg}$ we write \beaa \Box_g
u &= & \frac 1 {\sqrt {|\det g|}} \partial_\alpha\Big( \sqrt
{|\det g|} g^{\alpha \gbeta} \partial_\gbeta u\Big) \\ & = & \frac
1 {\sqrt
  {|\det \mu|}} \Big( 2\partial_\stau ( \sqrt {|\det \mu|} \partial_x
u) -2\partial_A\big ( \sqrt {|\det \mu|}\mu^{AB}\gamma_B \partial_x u)
\\ && +2\partial_x \big( \sqrt {|\det \mu|}( \partial_\stau
u+2{(-\chi+ |\gbeta|_\mu^2)} \partial_xu -\gbeta^\sharp u\big)
\Big) + \Delta_\mu u \;. \eeaa For $C_0\cap\Axd$--polyhomogeneous
metrics it only remains to check whether the terms multiplying
$\varphi$, as defined in (\ref{defred3}), are in $C_0\cap \Axd$.
Since all occurrences of a single $e_A$--derivative of $u$ can be
put into $\psi_A$, one only needs to keep track of terms which
arise as coefficients of $\varphi_+=e_+(u)=-2\partial_x u$; the
relevant collected contribution to $\Box_g$ reads
$$
2\Big( \partial_\stau ( \sqrt {|\det \mu|}) -\partial_A\big ( \sqrt
{|\det \mu|}\mu^{AB}\gamma_B) +2\partial_x \big( \sqrt {|\det
  \mu|}{(-\chi+ |\gbeta|_\mu^2)} \big) \Big)\frac{\partial_xu} {\sqrt
  {|\det \mu|}} \;.$$
The coefficient of $\partial_x u $ will be in
$C_\infty+x^\delta\Axyd$ if, \emph{e.g.},\/ $g\in \Axd\cap
C^1( \bcM)$, where $\bcM$ denotes the space-time manifold with the
conformal boundary $\scrip$ added.

Summarising, we have proved:

\begin{Proposition}
\label{Plinwave0} Consider a space-time $(\mcM,\stsg)$ with
$\Axd\cap C^1(\bcM)$-polyhomogeneous conformal completion at null
infinity, and suppose that \eq{C.3n} holds.
 The scalar wave equation on $(\mcM,\stsg)$ implies a system of
equations to which Theorem~\ref{Tmnl} applies.
\end{Proposition}

Similarly if $u$ is a solution of the non-linear wave equation
\eq{SE.1} as in Theorem~\ref{T2phg}, or of the wave-map equation
as in
Theorem~\ref{Twavemap}, then 
the lower order terms discussed so far are compatible with the
requirements of Theorem~\ref{Tlemme1n}.

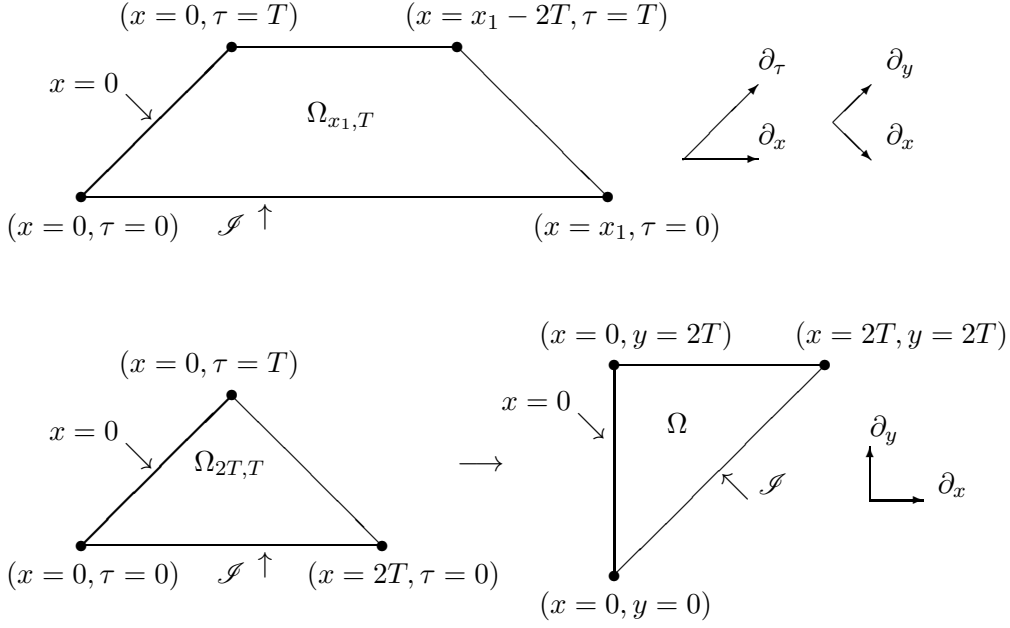
\begin{figure}
\setlength{\unitlength}{1cm} \noindent
\begin{picture}(10,4)(-1,-1)
\thicklines \put(0,0){\line(1,1){2}} \put(0,0){\circle*{.15}}
\thinlines \put(2,2){\line(1,0){3}} \put(2,2){\circle*{.15}}
\put(-0.43,1.4){$x=0$} \put(0.57,1.1){$\searrow$}
\put(7,0){\line(-1,0){7}} \put(7,0){\circle*{.15}}
\put(5,2){\line(1,-1){2}} \put(5,2){\circle*{.15}}
\put(-1,-0.5){$(x=0, \tau=0)$} \put(6,-0.5){$(x=x_1, \tau=0)$}
\put(3,1){$\Omega_{x_1,T}$} \put(0.5,2.3){$(x=0, \tau=T)$}
\put(4.3,2.3){$(x=x_1-2T, \tau=T)$}
\put(1.8,-0.5){$\mcS$} \put(2.35,-0.35){$\uparrow$}
\put(8,0.5){\vector(1,1){1}} \put(9,1.7){$\partial_\tau$}
\put(8,0.5){\vector(1,0){1}} \put(9,0.7){$\partial_x$}
\put(10,1){\vector(1,1){0.5}} \put(10.7,1.7){$\partial_y$}
\put(10,1){\vector(1,-1){0.5}} \put(10.7,0.7){$\partial_x$}
\end{picture}
\noindent
\begin{picture}(7,4)(-1,-1.4)
\thicklines \put(0,0){\line(1,1){2}} \put(0,0){\circle*{.15}}
\thinlines \put(4,0){\line(-1,0){4}} \put(4,0){\circle*{.15}}
\put(-0.43,1.4){$x=0$} \put(0.57,1.1){$\searrow$}
\put(2,2){\line(1,-1){2}} \put(2,2){\circle*{.15}}
\put(-1,-0.5){$(x=0, \tau=0)$} \put(3,-0.5){$(x=2T, \tau=0)$}
\put(1.5,1){$\Omega_{2T,T}$}
\put(5,1){$
{\longrightarrow}$} \put(0.5,2.3){$(x=0, \tau=T)$}
\put(1.8,-0.5){$\mcS$} \put(2.35,-0.35){$\uparrow$}
\end{picture}
\begin{picture}(10,5)(-1,-1)
\thicklines \put(0,0){\line(0,1){2.8}} \put(0,0){\circle*{.15}}
\thinlines \put(0,2.8){\line(1,0){2.8}} \put(0,2.8){\circle*{.15}}
\put(2.8,2.8){\line(-1,-1){2.8}} \put(2.8,2.8){\circle*{.15}}
\put(-1,-0.5){$(x=0, y=0)$} \put(-1,3.1){$(x=0, y=2T)$}
\put(2.4,3.1){$(x=2T, y=2T)$} \put(0.7,1.9){$\Omega$}
\put(-1.5,2.2){$x=0$} \put(-0.5,1.9){$\searrow$}
\put(1.42,1.08){$\nwarrow\mcS$} \put(3.4,1.0){\vector(0,1){0.71}}
\put(3.4,1.85){$\partial_y$} \put(3.4,1.0){\vector(1,0){0.71}}
\put(4.3,1.1){$\partial_x$}
\end{picture}
\caption{\label{F1}The variables $(x,\tau)$ and $(x,y)$.}
\end{figure}

We show, next, that any system of the form \eq{geneqsys} can be
written in a first order form as in Section~3
of~\cite{ChLengardnwe}. Recall that in that last reference a
coordinate system $(x,v^A,\tau)$ was used, with $\tau$ playing the
role of a time variable. The solutions there were defined on the
set
$$ \Omega_{x_1,T} = \{ (x,v^A,\tau) : 0<x<x_1-2\tau,\;
v\in\mcO, \;0<\tau<T\}\;.$$ This will be made compatible with
\eq{symhyp1} by replacing the coordinate $y$ with a coordinate
$\tau$,
\bel{chofvn} (x,y)\longrightarrow (x,\tau=\frac{y-x}2)\;, \ee so
that
\bel{chofvn2}
\partial_x\longrightarrow \partial_x -\frac{1}2\partial_\tau\;,
\quad \partial_y=\frac 12 \partial_\tau\;\ee (see
Figure~\ref{F1}). In this work we are interested in the behavior
of solutions near the set $\{x=0\}$. For this purpose it suffices
to study the properties of the solution $u$ on the following
subset of $\Omega_{x_1,T}$:
\bel{Omdef}\Omega = \{(x,v^A,y) : 0<x<y,\; v\in\mcO,\;
0<y<2T\}\;.\ee
In view of \eq{chofvn2}, \eq{defred1} takes the form
\beal{defred4} &e_-=\partial_\tau + {(-\chi+ |\gbeta|_\mu^2)}
(2\partial_x -\partial_\tau)-2\gbeta^\sharp\;,\quad
e_+=\partial_\tau-2\partial_x\;.&\eea
     We keep the definitions \eq{defred2}-\eq{defred3},
and rewrite \eq{symhyp1} as
\beal{symhyp2}
  &\left.\begin{array}{llll}
(\partial_\tau-2\partial_x) \psi_A &-e_A
 \varphi_+&-
[\partial_\tau-2\partial_x,e_A] u &=  0\;,
 \\
 -\sum_A e_A
 \psi_A  &+e_- \varphi_+& +\sum_A e_A(e_A{}^C)f_C{}^B \psi_B   &=  -F(x^\mu,\psi,\varphi)\;,
 \end{array}\right\}&
 \phantom{xxx}
  \\
  &(\partial_\tau-2\partial_x) \psize
-\varphi_+= 0\;,&
  \\
  &\left.\begin{array}{llll}
  (\partial_\tau-2\partial_x) \psi_-& -\sum_A e_A
 \varphi_A &- [\partial_\tau-2\partial_x,e_-] u +\sum_A e_A(e_A{}^C)f_C{}^B \psi_B  &= -F(x^\mu,\psi,\varphi)\;,
 \\
 -e_A\psi_-&+e_-\varphi_A&- [e_-,e_A]u &=0\;.
 \end{array}\right\}\phantom{}& \eeal{symhyp5}

 If we assume that the
 metric coefficients are $(\Axd\cap C_1)$--polyhomogeneous, then all
 coefficients in \eq{symhyp2}-\eq{symhyp5} are of at least $\Axd\cap C_0$--differentiability class.

 The system \eq{symhyp2}-\eq{symhyp5} is of the form considered in
 \cite[Section~3]{ChLengardnwe}: For example, the matrices $E^\mu_\pm$
 of Equation~3.5 of \cite{ChLengardnwe} take the form $E^\mu_\pm =
 e^\mu_\pm\otimes \Id$, with $e_\pm=e^\mu_\pm\partial_\mu$ as in
 \eq{defred4}.  The covector field $n_\mu$ in \chleq{3.12} can be
 taken to be equal to $d\tau$. The matrices $\ell^A$ can be read-off from \eq{symhyp2}-\eq{symhyp5}:
$$\ell^A\partial_A\psi = \left(
\begin{array}{ccccc}
  0 & 0 & -e_1 & \ldots & -e_{n-1} \\
  0 & -e_1 & 0 &  & 0 \\
   & \vdots &  & \ddots &  \\
  0 & -e_{n-1} & 0 &  & 0 \\
\end{array}\right) \mpsi\;.$$
A natural choice for $L$ is then
$$L\psi = \left(
\begin{array}{ccccc}
  0 & 0 & -e_1-\frac{e_1(\sqrt{\det\mu})}{\sqrt{\det\mu}} & \ldots &
   -e_{n-1}-\frac{e_{n-1}(\sqrt{\det\mu})}{\sqrt{\det\mu}} \\
  0 & -e_1 & 0 &  & 0 \\
   & \vdots &  & \ddots &  \\
  0 & -e_{n-1} & 0 &  & 0 \\
\end{array}\right)\mpsi\;,
$$
so that
$$\ell\psi =\left(
\begin{array}{ccccc}
  0 & 0 & -\frac{e_1(\sqrt{\det\mu})}{\sqrt{\det\mu}} &
  \ldots & -\frac{e_{n-1}(\sqrt{\det\mu})}{\sqrt{\det\mu}} \\
  0 & 0 & 0 &  & 0 \\
   & \vdots &  & \ddots &  \\
  0 & 0 & 0 &  & 0 \\
\end{array}\right)\mpsi\;.$$
This implies
$$-L^\dagger\varphi = \left(\begin{array}{cccc}
0 & 0 &  & 0\\
0 & -\dyw e_1 - e_1 & \ldots & -\dyw e_{n-1} - e_{n-1} \\
\frac{e_1(\sqrt{\det\mu})}{\sqrt{\det\mu}} - \dyw e_1 - e_1 & 0
&&0 \\
\vdots &&\ddots& \\
\frac{e_{n-1}(\sqrt{\det\mu})}{\sqrt{\det\mu}} - \dyw e_{n-1} -
e_{n-1} & 0
&&0 \\
\end{array}\right)\mphi\;,$$
where $\dyw$ denotes the divergence with respect to the measure
$\sqrt{\det \mu_{AB}}\,dv^1\cdots dv^{n-1}$.  It immediately follows
that \chcite{Equations~(3.9) and (3.11)} hold for $\Axd\cap C_1$
metrics.

We have
 $$E^\tau_+ = \Id\;,\qquad E^\tau_-= (1+\chi- |\gbeta|_\mu^2)\Id\;,$$
 which
 shows that the first condition in \chleq{3.12} holds. On the other
 hand, \cite[Equation~(3.13)]{ChLengardnwe} holds only for $E^\tau_+$.
 However, we note the following result, where the notation from
 \cite{ChLengardnwe} is used:
 \begin{Proposition}
 \label{Pspden}
 Under the remaining conditions of
 \cite[Proposition~3.1]{ChLengardnwe}, suppose that instead of
 \cite[Equation~(3.13)]{ChLengardnwe} we have
 \bel{spden1} \| E^\tau_-
(\tau)\|_{{\cG}^0_{k}\mxtau } \le C\;,\qquad   \| E^\tau_+
(\tau)\|_{{\cG}^{-1}_{k}\mxtau }+ \| \partial_x E^\tau_
+
 (\tau)\|_{{\cG}^{-1}_{k-1}(M_{x})}\leq  \err(x)\;.\ee
Then the conclusions of \cite[Proposition~3.1]{ChLengardnwe} hold.
 \end{Proposition}

 \proof If \chleq{3.13} does not hold, then the commutators in the
 unnumbered equation after \chleq{3.28} will have supplementary terms
 such as, \emph{e.g.}, \bel{suptcon}\sum_{i=1}^k ({}^i_k) x^i
 (\partial_ x ^i E_{\pm}^\tau)x^{k-i} \px^{k-i}\partial_\tau
 \chi\;.\ee From the evolution equations \chleq{3.5} one has
 \beal{tauder1} \partial_\tau \varphi & = &
 (E^\tau_-)^{-1}\Big(-E^x_-\partial_x \varphi - E^A_- \partial_A
 \varphi -(B_-+B_{11})\varphi- L\psi - B_{12}\psi+a\Big)\;,
 \\
 \partial_\tau \psi & = & (E^\tau_+)^{-1}\Big(-E^x_+\partial_x \psi -
 E^A_+ \partial_A \psi -(B_++B_{22})\psi+L^\dagger \varphi -
 B_{21}\varphi+b\Big)\;.\phantom{xxxxx} \eeal{tauder2} This is
 inserted in \eq{suptcon}. Some care has to be taken, because the norm
 used requires that each $x$--derivative comes with a factor of $x$.
 This works well for $\varphi$, because of the requirement in
 \chleq{3.14} that $E^x_- \in {\cG}^{1}_{k}\mxtau$. On the other hand,
 the $x$-derivatives of $\psi$ that result in \eq{suptcon} are
 rewritten as \bel{suptcon2}\sum_{i=1}^k ({}^i_k) (x^{i-1}\partial_ x
 ^{i-1} \partial_x E_{+}^\tau)x^{k-i+1}
 \px^{k-i}\Big((E^\tau_+)^{-1}E^x_+\partial_x \psi\Big) \;,\ee
and  the  estimates are done as in \chl\ starting from this
formula. The remaining supplementary terms are handled in  a
similar way. \qed

An identical argument gives:

 \begin{Proposition}
 \label{Pspden2}
 Under the remaining conditions of
 \cite[Proposition~3.2]{ChLengardnwe}, suppose that instead of
 \cite[Equation~(3.13)]{ChLengardnwe} we have
 \bel{spden3} \| E^\tau_-
(\tau)\|_{{\cG}^0_{k}\mxtau } + \| E^\tau_+
(\tau)\|_{{\cG}^0_{k}\mxtau } +\| \partial_x E^\tau_ +
 (\tau)\|_{{\cG}^{0}_{k-1}(M_{x})}\le C\;.\ee
Then the conclusions of \cite[Proposition~3.2]{ChLengardnwe} hold.

 \qed \end{Proposition}

The next question we wish to address is that of the choice of the
initial data surface $\hyp$. The estimates of~\chl\ were made
under the assumption that
$$\hyp=\{\tau=0\}\;.$$
This might be overly restrictive for several
purposes. The simplest generalisation is to assume that $\hyp$
coincides with a member of the family of hypersurfaces $\hyp_s$
defined as \bel{hhyp0}\hyp_s=\{\tau= s+ \alphahyp(x,v^A)\}\;,\ee where
$\alphahyp$ is a function satisfying
 \bel{hhyp} 
 \alphahyp \in \fCx1{k}\;.\ee
The field of conormals to the $\hyp_s$'s is $d\tau - \partial_x
\alphahyp\, dx - \partial_A \alphahyp \,dv^A$, so that these
hypersurfaces will be space-like for the system considered in
\chl\ if
\bel{hhyp2}
\left(
\begin{array}{cc} E^\tau_- - E^x_-\partial_x\alphahyp
- E^A_- \partial_A \alphahyp & -\ell^A\partial_A\alphahyp
\\
-(\ell^A)^t\partial_A\alphahyp & E^\tau_+- E^x_+ \partial_x\alphahyp -
 E^A_+ \partial_A \alphahyp
\\
\end{array}
\right)
\ge \epsilon>0\;,\ee for some constant $\epsilon$. We have
\begin{Proposition}
\label{PgencS} Assume \eq{hhyp0}-\eq{hhyp2}, and  suppose that the
coefficients in the equation remain in the original space after
differentiation by $(x\partial_\tau)^i$, with $1\le i \le k-1$. Then
the estimates of \cite[Propositions~3.1 and 3.2]{ChLengardnwe}, as
generalised in Propositions~\ref{Pspden} and \ref{Pspden2}, remain
valid for solutions with initial data on $\hyp_0$, with the norms
$\|f(t)\|$ in \chl\ understood as norms of functions on $\hyp_t$.
\end{Proposition}

\proof  One integrates the
divergence \chleq{3.22} on the family of sets, parameterised by
$s$,
$$\{\alphahyp(x,v^A)\le \tau \le s+\alphahyp(x,v^A)\;,\ 0\le x\le x_1-2\tau\}
\;.$$ The argument in \chl\ leads then to an estimate on the
integrals over the $\hyp_s$'s of $x$-- and $v^A$--derivatives of
the fields. To relate those integrals to intrinsic Sobolev spaces
on the $\hyp_s$'s it is convenient to pass to a new coordinate
system $(\tilde\tau,\tilde x,\tilde v^A)$ defined by the equations
$$\tilde \tau=\tau-\alphahyp(x,v^A)\;,\quad \tilde x=x\;,\quad \tilde v^A =
v^A\;.$$ We then have
$$\partial_{\tilde x} = \partial_x + \partial_x \alphahyp\,
\partial_\tau\;,\qquad \partial_{\tilde v^A} = \partial_{v^A} +
\partial_{v^A}\alphahyp \,\partial_\tau\;,$$
so that
\begin{eqnarray*} \lefteqn{\tilde{x}^i \partial_{\tilde v}^\gamma \partial_{\tilde x}^i f
= x^i\partial_{v}^\gamma \partial_{x}^i f +}&&\\&&
\sum C(j_0, \ldots, j_k, \gamma_0, \ldots,
\gamma_k) (x^{j_0+k}
\partial_{v}^{\gamma_0} \partial_{x}^{j_0}\partial_\tau^k f)
(x^{j_1}\partial_x^{j_1} \partial_{v}^{\gamma_1} \sigma) \ldots
(x^{j_k}\partial_x^{j_k} \partial_{v}^{\gamma_k}
\sigma)\;,\end{eqnarray*} where the sum is taken over all indices satisfying $j_0 + \ldots + j_k = i$,
$|\gamma_0| + \ldots + |\gamma_k| = |\gamma|$, with $1\le k \le i+|\gamma|$.
One can use \eq{tauder1}-\eq{tauder2} to eliminate
$\tau$--derivatives; the fact that such derivatives of $\psi$ are
replaced by $x$-derivatives introduces the restriction
$\partial_A\sigma=O(x)$ when $k=1$, and to \eq{hhyp} for higher $k$.
\qed

 To proceed
 further, one needs to make a choice for the scalar products and covariant derivative
 operators. We
assume that
 $u=(u^a)$ maps the original space-time into $\R^m$, with canonical scalar
 product.
  For a scalar field we take $m=1$. In case of wave maps, we will
  be interested in maps that asymptote a point $p$ in the target
  manifold when we approach $\{x=0\}$. We choose normal
  coordinates near $p$, and the $u^a$'s are the coordinate components of the
  map.

  We take $E^\mu_\pm \nabla_\mu = E^\mu_\pm \partial_\mu$, so that the
  matrices $B_\pm$ in \chleq{3.16} vanish. We choose
 \beaa |\psi|^2&=& \sum_{a,A} (\psi^a_A)^2 
 + \sum_a
 \Big((\psize ^a)^2+(\psi^a_-)^2\Big)\;,
 \\|\varphi|^2&=& \sum_{a,A}(\varphi^a_A)^2  
 + \sum_a
(\varphi^a_+)^2\;.
 \eeaa
The ``energy-momentum vector'' $X$ in
\cite[Equation~(3.21)]{ChLengardnwe} equals then \bel{X}X=
|\varphi|^2 e_- + |\psi|^2 e_+\;,\ee and we use the usual measure
$\sqrt{|\det g|}d^{n+1}x=\frac 12 \sqrt{|\det \mu|}d^{n+1}x$ to
integrate over domains in the (conformally rescaled) space-time,
so that
 $$\mathrm{div} X = \nabla_\mu X^\mu = \partial_\mu (\sqrt{|\det
 \mu|} X^\mu)/\sqrt{|\det \mu|}\;.$$

 Rather than deriving
 $\nabla_\alpha E^\alpha_\pm$ from first principles, it is simpler to calculate
directly $\mathrm{div} X$ using \eq{X} and estimate the
corresponding terms. This gives
$$|\nabla_\alpha E^\alpha_\pm| \le |\nabla_\alpha e^\alpha_\pm|
\;.$$
It now follows from \eq{defred1} that the relevant estimates in
\chleq{3.14} and \chleq{3.15} hold for any $(\Axd\cap
C_1)$--polyhomogeneous metric.

Summarising, we have proved:

\begin{Proposition}\label{Plinwave} The estimates of \chcite{Propositions~3.1
and 3.2} hold for the scalar wave equation
$$\Box_{\stsg}u=0$$
on any space-time with $(\Axd\cap
C_1)$-polyhomogeneous conformal completion at null infinity for which
\begin{equation} \Omega^{-{n+3\over
2}}\left(R(\tilde{\stsg})\Omega^2 -R(\stsg)\right)\in
\Axd\;. \label{C.3nn}
\end{equation}
\end{Proposition}

\subsection{Semi-linear wave equations}\label{sSslwe}

We continue with existence of solutions of the semi-linear wave
equation: \be \Box_\stsg f = H(x^\mu,f)\;.\label{SE.1} \ee The spaces
$ \mcH^\beta_{\ell}(\hyp_0)$ below are defined in~\chcite{Appendix~A}:

\begin{Theorem}\label{Tnlweexpt} Consider a space-time $(\mcM,\stsg)$ with
  $\Axd\cap C_1(\bcM)$-polyhomogeneous conformal completion at null
  infinity satisfying (\ref{C.3nn}), where $\delta$ is the inverse of
  a strictly positive integer. Let $\hyp_0$ be a spacelike hypersurface in
  $\mcM$ of the form \eq{hhyp0}, with $\alphahyp$
  satisfying \eq{hhyp} together with
 \bel{hhyp2b}\liminf_{x\to 0}
  \partial_x\sigma >-1/2\;.\ee Set $$\tilde f = \Omega^{-(n-1)/2}f\;,$$
  and suppose
  that the initial data for \eq{SE.1} on $\hyp_0$ satisfy, for some
  $k>\lfloor \frac n2 \rfloor +1$ and $-1<\alpha<-1/2$,
$$\tf|_{\hyp_0} \in \mcH^\alpha_{k+1}(\hyp_0)\;,\qquad
\partial_x\tf|_{\hyp_0} \in
 \Big(\mcH^{\alpha-1/2}_{k}\cap \fCx{\alpha}{0}\Big)(\hyp_0)\;,$$
$$\partial_\tau\tf|_{\hyp_0} \in \mcH^\alpha_{k}(\hyp_0)\;.$$
Suppose  that
  $H(x^\mu,f)$ is smooth in $f$ at fixed $x^\mu$, bounded and
  $\delta$-polyhomogeneous in $x^\mu$ at constant $f$, and has a zero of order
  $\ell$ at $f=0$, with $\ell$ satisfying \eq{condH}.  There exists $T>0$ and a
solution of \eq{SE.1} defined on
 $\cup_{s\in[0,T]}\hyp_s$ such that
 \minilab{timderreg}
\begin{equs}
  \label{timderreg1}
   \tf &\in
   \fCxy{\alpha}{k-1-\lfloor \frac n2 \rfloor}(\cup_{s\in[0,T]}\hyp_s)\cap L^\infty\;,
 \\
  \label{timderreg2}
    \partial_x\tf &\in
     \fCxy{\alpha-1/2}{k-1-\lfloor \frac n2 \rfloor}(\cup_{s\in[0,T]}\hyp_s)\cap
     \fCx{\alpha}{0}
    \;,
 \\
  \label{timderreg3}
    \partial_\tau\tf &\in
     \fCxy{\alpha}{k-1-\lfloor \frac n2\rfloor}(\cup_{s\in[0,T]}\hyp_s)
    \;.
\end{equs}
 \end{Theorem}

 \proof
 The proof of existence is essentially a repetition of that of
 \chcite{Theorem 4.1}. Condition~\eq{hhyp2b} guarantees that
 \eq{hhyp2} holds for $x$ small enough. Using the results established
 earlier in this section, one proceeds as in \chl\ until
 Equation~\chleq{4.32}.  When $\tilde \stsg$ is smooth, that last
 equation is handled by introducing the variables $(\hat x, \hat y,
 \hat v^A)$ as in \eq{hatvar0}.  This leads again to the estimate
 \chleq{4.33}.  In order to obtain \chleq{4.34} one can proceed as
 follows: in the coordinate system $(x,y,v^A)$ we have $\partial_x \tf
 = -\frac 12 \varphi_+$, and integration in $x$ until the initial data
 hypersurface $\{\tau=\alphahyp(x,v^A)\}$ gives
 \bel{newinteq}\tf(x,y,v^A)= \tf(x_0(y,v^A),y,v^A)+\frac 12
 \int_x^{x_0(y,v^A)} \varphi_+(s,y,v^A)ds\;, \ee where $x_0$ is
 implicitly defined by the equation \bel{xoqe}
 x_0+2\alphahyp(x_0,v^A)=y\;\ee (note that \eq{hhyp2b} implies that
 the left-hand-side is a strictly monotonous function of $x_0$,
 guaranteeing solvability of \eq{xoqe}).  Reverting to the variables
 $(\tau,x,v^A)$, this leads to the estimate
\bel{revestt}|\tf(\tau,x,v^A)|\le  \|\tf|_{\hyp_0}\|_{L^\infty}
+\frac 12 \sup_{s\in[0,\tau]}\|\varphi_+(s)\|_{\fCx{\alpha}{0}}
 \int_x^{x_1} s^{-\alpha}ds\;.\ee
 Clearly, $\sup_{s\in[0,\tau]}\|\varphi_+(s)\|_{\fCx{\alpha}{0}}$
 satisfies the same estimate \chleq{4.33} as
$ \|\varphi_+(\tau)\|_{\fCx{\alpha}{0}}$, which allows us to
recover \chleq{4.34}. The remaining arguments of
\chcite{Theorem~4.1} apply without changes. This provides existence,
and space-derivative estimates, for smooth metrics $\tilde \stsg$.

When $\tilde \stsg$ is polyhomogeneous, a little more work is needed:
the variables $(\hat x, \hat y, \hat v^A)$ are introduced as in
\eq{chhatcrd}.  This allows us to rewrite the fourth equation in
\eq{symhyp1} in the form
$$\partial_{\hat y}\varphi_+ + \underbrace{\Gamma}_{\in C_\infty}
\varphi_+ = \rho\;,$$
where, in the notation of \chl, the source term
can be estimated as, for any $\lambda<0$, \beaa \|\rho(\tau)\|_{\fCx
  \lambda 0} & \le & C\Big(\underbrace{\|\psi(\tau)\|_{\fCx
  {\lambda} 1} 
+\|G(\tau)\|_{\fCx {\lambda} 0}}_{\le CE_\lambda(\tau)}
\\
&& +\underbrace{\|\varphi(\tau)\|_{\fCx
    {\lambda-\delta+\epsilon}1}}_{\le C \|\varphi(\tau)\|_{\fCx
    {\lambda-\frac 12 +\epsilon}1}\le C
  E_{\lambda+\epsilon}(\tau)}+\|\tilde f(\tau)\|_{L^\infty}\Big)\;.
\eeaa Here $G=-\Omega^{-(n+3)/2}H$. The last term above arises
from the scalar curvature terms in \eq{C.3}. The number
$\epsilon>0$, due to possible powers of $\ln x$ in the
coefficients in the equation, can be chosen to be arbitrarily
small. We choose $\epsilon$ so that $\alpha-\epsilon>-1$, set
$\lambda=\alpha-\epsilon$ and obtain an obvious modification of
\chleq{4.33} with $\alpha$ there replaced by $\lambda$.
\Eq{revestt} with $\alpha$ replaced by $\alpha-\epsilon$ leads
then to the following variation of \chleq{4.34}:
$$|\tf(\tau,x,v^A)|\le  \|\tf|_{\hyp_0}\|_{L^\infty}
+Ce^{C\tau}\|\varphi_+(0)\|_{\fCx {\alpha-\epsilon} 0}
+ \int_0^{\tau} e^{C(\tau-s)}C(E_{\alpha}(s),E_{\alpha-\epsilon}(s),
\|\tf(s)\|_{L^\infty})ds\;,$$
for a certain continuous function $C$, bounded on bounded sets.
We combine the estimate \chleq{4.29} with its equivalent where $\alpha$ is replaced by $\alpha-\epsilon$ to obtain
\begin{eqnarray}
  \|\widetilde f(\tau)\|_{L^\infty}+
E_\alpha(\tau)+ E_{\alpha-\epsilon}(\tau) &\leq& \|\tf|_{\hyp_0}\|_{L^\infty} +
  Ce^{C\tau}\left(E_\alpha({0})+ E_{\alpha-\epsilon}(0)
+\|\varphi_+(0)\|_{\fCx {\alpha-\epsilon} 0} \right) \nonumber \\ &&
 + \int_{0}^{\tau}
\Phi\left(\tau,s,\|\widetilde f(s)\|_{L^\infty},
E_\alpha(s), E_{\alpha-\epsilon}(s)\right) ds\;,
  \label{eq:glbd}\end{eqnarray}
and one concludes as in \chl.


The estimates \eq{timderreg} on the time-derivatives are obtained by a
repetition of the argument of the case $m=0$ of \chcite{Theorem~4.4},
together with a weighted Sobolev embedding.  \qed

We are ready now to prove the first main result of this paper:
\begin{Theorem}\label{T2phgg}  Let
$\delta =1$ in odd space dimensions,
  and let $\delta =1/2$ in even space dimensions. Under the hypotheses of
  Theorem~\ref{Tnlweexpt},
  suppose moreover that  $H$ is $\Axyd$-polyhomogeneous in $x^\mu$ with a
  uniform zero of order $l$, and that
  both the function $\alphahyp$ defining $\hyp_0$, and the initial data are
  polyhomogenous:
  $$\sigma \in \Axd\cap C_{1}\;,$$
$$ f|_{\hyp_0}\in x^{(n-1)/2} \Big(\Axd\cap L^\infty\Big)\;,\quad
\partial_\tau f|_{\hyp_0}\in x^{(n-1)/2} \Axd\;. $$
  Then the solution given by Theorem~\ref{Tnlweexpt} satisfies
   $$f\in x^{(n-1)/2} \Big(\Axyd\cap L^\infty\Big)(\cup_{s\in[0,T]}\hyp_s)\;.$$
   In particular for every $\epsilon>0$ we have $f|_{\cup_{s\in[\epsilon,T]}\hyp_s}\in \Axd\cap
  L^\infty$.
\end{Theorem}

\proof Setting $\mcO=\hyp_0\cap \mcN$ in the construction at
the beginning of this section leads to a coordinate system in which
$\hyp_0\cap \mcN\subset \{y=0\}$.  As explained at the beginning of
Section~\ref{Spwcc}, we can then without loss of generality assume
that $\hyp_0=\{y=x\}$.  Theorem~\ref{Tnlweexpt} with $k=\infty$ and
Proposition~\ref{Plinwave0} show that the hypotheses of
Theorem~\ref{Tlemme1n} below are satisfied, and the result follows.
\qed

More precise information on the behavior of $\tf$ can be obtained from
Theorem~\ref{Tlemme1n} and Remarks~\ref{Rldcb}-\ref{Rldcc} below.

\subsection{Wave maps}\label{sSwme}
Let $(N,h)$ be a smooth Riemannian manifold, and let
$f:(\cM,\stsg)\to(N,h)$ solve the wave map equation. As already
pointed out, we will be interested in maps $f$ which have the
property that $f$ approaches a constant map $f_0$ as $r$ tends to
infinity along lightlike directions, $f_0(x)=p_0\in N$ for all
$x\in\cM$. For the purposes of the analysis in \chl\ it was useful
to use normal coordinates around $p_0$, and if we write $f=(f^a)$,
$a=1,\ldots,\dim N$, then the functions $f^a$ satisfy the set
of equations
\begin{equation}
  \label{W.1}
\Box_\stsg f^a + \stsg^{\mu\nu}\Gamma^{a}_{bc}(f) \frac{\partial
  f^b}{\partial x^\mu} \frac{\partial f^c}{\partial x^\nu} = 0\;,
\end{equation}
where the $\Gamma^{a}_{bc}$'s are the Christoffel symbols of the
metric $h$. Setting as before $\tf^a=\Omega^{-{n-1 \over 2}}f^a$,
$\tilde\stsg=\Omega^2\stsg$, we then have from \eq{C.3},
\begin{equation}
  \label{W.2}
\Box_{\tilde\stsg} \tf^a =- \Omega^{-{n-1 \over 2}}
\tilde\stsg^{\mu\nu}\Gamma^{a}_{bc}(\Omega^{n-1 \over 2}\tf)
\frac{\partial   (\Omega^{n-1 \over 2}\tf^b)}{\partial x^\mu}
\frac{\partial (\Omega^{n-1 \over 2}\tf^c)}{\partial x^\nu} + {n-1
\over 4n}\Omega^{-{n+3 \over 2}}(\tilde{R} -R\Omega^{-2})\tf^a\;.
\end{equation}
Proceeding as in \eq{defred1}-\eq{defred3} with each component
$f^a$ of $f$, we obtain the system of Equations \eq{symhyp1}, with
the obvious replacements associated with $\tf=\Omega^{-(n-1) \over
2}f \to \tf^a$, and with $F=(F^a)$ taking the form
\begin{eqnarray}
\lefteqn{F^a:= - \Gamma^{a}_{bc} (\Omega^{n-1 \over 2}\tf)
\Omega^{{n-1 \over 2}} \bigg\{-\varphi_+^b\psi_-^c + \varphi_A^b
\varphi_A^c } & & \nonumber\\&& {}-
\frac{n-1}{2\Omega}\left(e_+(\Omega) \psi_0^b \psi_-^c +
e_-(\Omega) \psi_0^c \varphi_+^b\right) + \frac{n-1}{
\Omega}e_A(\Omega)\varphi_A^b\psi_0^c \nonumber \\
&& {}+ \frac{(n-1)^2}{4\Omega^2} |\nabla_{\tilde\stsg}\Omega|^2
\psi_0^b \psi_0^c\bigg\}
  + \frac{n-1}{4n}(\tilde{R} -R\Omega^{-2})\psi_0^a\;.
  \label{W.3}\end{eqnarray}
We have
\begin{Theorem}\label{Twmexpt} Consider a space-time $(\mcM,\stsg)$ with
  $\Axd\cap C_1(\bcM)$-polyhomogeneous conformal completion at null
  infinity satisfying \eq{C.3nn}, where $\delta$ is the inverse of a
  strictly positive integer. Suppose that the conformal factor
  $\Omega$ is of $(C_2\cap \Axd)$--differentiability class. Let $\hyp_0$
  be a spacelike hypersurface in $\mcM$ of the form \eq{hhyp0}, with
  $\alphahyp\in\fCx 1 k $ satisfying \eq{hhyp2b}.  Suppose that the
  initial data for \eq{W.3} on $\hyp_0$ satisfy \beqa\label{cpt1.0w}
  \tf ^a|_{\hyp_{0}} \equiv\Omega^{-(n-1)\over 2}f^a|_{\hyp_{0}} &
  \in& \cases{ \left(\HH_{k+1}^{\alpha}\cap L^\infty\right)(\hyp_0)
    \;, & $n\geq 3\;,$\cr \left(\HH_{k+1}^{\alpha}\cap
      {{\mcC}}^0_1\right)(\hyp_0) \;, & $n=2\;,$} \\ \partial_x \tf
  ^a|_{\hyp_{0}} &\in &
  \HH_k^{\alpha}(\hyp_0)\;, \label{cpt1.1w}\\
  \partial_\tau\tf ^a|_{\hyp_{0}} &\in&
  \cases{ \HH_{k}^{\alpha}(\hyp_0) \;, & $n\geq 3\;,$\cr
  \left(\HH_{k}^{\alpha}\cap L^\infty\right)(\hyp_0) \;, &
  $n=2\;.$}
\label{cpt1w}
 \eeqa  for some $k > {n\over2} +1$, $ -1< \alpha
\le -1/2$.  If $n=2$ we will moreover assume a single ``corner
condition''
\be
\partial_\tau^{2}\widetilde
  f^a|_{\hyp_{0}}\in \HH^{-1}_{k-1}({\hyp_{0}})\;.\label{2dgcpt1.3wgt}
 \ee
 There exists $T>0$ and a solution of \eq{SE.1}
defined on
 $\cup_{s\in[0,T]}\hyp_s$ such that
 \minilab{timderregwm}
\begin{equs}
  \label{timderregwm1}
   \tf &\in
   \fCxy{\alpha}{k-1-\lfloor \frac n2\rfloor}(\cup_{s\in[0,T]}\hyp_s)\cap L^\infty\;,
 \\
  \label{timderregwm2}
    \partial_x\tf\;,\  
    \partial_\tau\tf &\in
     \fCxy{\alpha}{k-1-\lfloor \frac n2\rfloor}(\cup_{s\in[0,T]}\hyp_s)
    \;,
 \end{equs}
 for $n\ge 3$. If $n=2$ the differentiability indices in
 \eq{timderregwm} have to be replaced by $k-3$.
 \end{Theorem}

 \proof The proof of existence is, again, a repetition of that of
 \chcite{Theorem~5.1}. The reader will note that, by definition of
 null infinity, the function $|\nabla \Omega|_{\tilde \stsg}^2$
 vanishes at $\mcN$ so that $|\nabla \Omega|_{\tilde \stsg}^2/\Omega
 \in L^\infty$. Similarly $e_-(\Omega)/\Omega\in L^\infty$, which
 implies that all those terms in \eq{W.3} which are already present in
 \chleq{5.3} have the right structure. The last term in the third line
 of \eq{W.3}, which is new, can be absorbed in the linear part of the
 operator. The last term in the second line of \eq{W.3}, which again
 is not present in \chleq{5.3}, can be controlled using $ab\le
 a^2+b^2$, which allows one to absord this term in the remaining ones
 (note that $e_A(\Omega)/\Omega\in L^\infty$) for the purpose of
 estimates.  Those estimates can be done as in \chl\ until
 Equation~\chleq{5.22}, which is replaced by \eq{newinteq}. This leads
 to the inequality \chleq{5.23} (the $\varphi_-$ term there is actually
 not necessary in the intermediate inequality). In space-dimension two
 one further obtains \chleq{5.24} by differentiating \eq{newinteq}
 with respect to $v^A$.

The estimates \eq{timderregwm} on the time-derivatives are
obtained by a repetition of the argument of the case $m=0$ of
\chcite{Theorems~5.4 and 5.5}. \qed

This leads us to the second main result of this paper, which is a
consequence of Theorems~\ref{Twmexpt} and \ref{Tlemme1n}, as in
the proof of  Theorem~\ref{T2phgg}:
\begin{Theorem}\label{Twavemapg}  Let
$\delta =1$ in odd space dimensions,
  and let $\delta =1/2$ in even space dimensions. Under the hypotheses of
  Theorem~\ref{Twmexpt},
  suppose moreover that
  both the function $\alphahyp$ defining $\hyp_0$, and the initial data are
  polyhomogenous:
  $$\sigma \in \Axd\cap W^{1,\infty}\;,$$ 
\beqa\label{cpt1.0wphg} \tf ^a|_{\hyp_{0}}
  \in
 \Axd\cap L^\infty \;, \qquad
  \partial_\tau\tf ^a|_{\hyp_{0}} \in
  \cases{\Axd \;, & $n\geq 3\;,$\cr
 \Axd\cap L^\infty \;, &
  $n=2\;.$}
\label{cpt1wphg}
 \eeqa
If $n=2$ we will moreover assume the corner condition
 \be
\partial_\tau^{2}\widetilde
  f^a|_{\hyp_{0}}=o(x^{-1})\;.\label{2dgcpt1.3wgtx}
 \ee
  Then the solution given by Theorem~\ref{Twmexpt} satisfies
   $$f\in x^{(n-1)/2} \Big(\Axyd\cap L^\infty\Big)(\cup_{s\in[0,T]}\hyp_s)\;.$$
   In particular for every $\epsilon>0$ we have $f|_{\cup_{s\in[\epsilon,T]}\hyp_s}\in \Axd\cap
  L^\infty$.

\qed
\end{Theorem}

As before, the reader will find more precise information on the
behavior of $\tf$ in Theorem~\ref{Tlemme1n} and
Remarks~\ref{Rldcb}-\ref{Rldcc} below.

\section{Polyhomogeneity without corner conditions}\label{SsCp}

We are ready now to prove the polyhomogeneity of solutions for a class
of linear systems of equations of the form {\rm
\begin{deqarr}
  \partial_y\varphi + B_{\varphi\varphi}\varphi +
  B_{\varphi\psi}\psi&=& L_{\varphi\varphi}\varphi +
  L_{\varphi\psi}\psi + a \;,\label{s1a}\\ \partial_x \psi+
  B_{\psi\varphi}\varphi + B_{\psi\psi}\psi &=& L_{\psi\varphi}\varphi
  + L_{\psi\psi}\psi + b \label{s1b}\;,
  \arrlabel{s1}\end{deqarr}}where $\varphi$ and $\psi$ are vector
valued, the $B$'s are linear maps and the $L$'s are first order linear
differential operators, while $a$ and $b$ are --- possibly non-linear
--- source terms.
Such equations arise from semi-- and quasi--linear wave equations in
Minkowski space-time, say \bel{nwepnc1} \Box u = F(x,u,\partial
u)\;.\ee As discussed in detail in Section~\ref{Sweci}, one obtains a
system of the form \eq{s1} by introducing $$\psi=(u, \partial_y u)\;,
\quad \varphi= \partial_x u\;.$$
Then both source terms $a$ and $b$ in
\eq{s1} are proportional to $F$. For several classes of
non-linearities, solutions of such equations with polyhomogeneous (or
smooth) initial data will satisfy $u\in \cDai\cap L^\infty$
\cite{ChLengardnwe}, see also Section~\ref{Sweci}. It also follows
from the results there that the natural polyhomogeneity space for
solutions of \eq{nwepnc1} is $u\in \Axd$, $\delta =1/2$ in even
space-dimensions, and $\delta =1$ in odd dimensions, at least when an
infinite number of corner conditions is satisfied. For such solutions
$u$ and $\partial_y u$ will typically be bounded, so that $\psi$ will
be bounded, while $\varphi$ will typically blow up, at $x=0$, as
$x^{\delta -1}\ln ^N x$ for some $N$. The source terms $a$ and $b$
behave then in a similar way near $x=0$. Since $b$ appears in an
equation where the polyhomogeneous function $\psi$ is differentiated,
one expects $b$ to behave again as $x^{\delta -1}\ln ^{N'} x$ for some
$N'$.  This discussion justifies the hypotheses of our theorems below.

\subsection{Linear equations}\label{Spwcc}

We consider the Cauchy problem with polyhomogeneous coefficients,
polyhomogeneous sources, as well as polyhomogeneous initial data given
on a polyhomogeneous hypersurface \bel{hypini}\hyp=\{y=\yhypini
(x,v)\}\;,\quad \mbox{with} \ \yhypini(0,v)=0\;,\
\partial_x\yhypini(0,v)>0\;,\ee for some function $\lambda\in\Axd\cap
C^1(\overline \Omega)$.  This can be reduced to the case $\yhypini=x$
by a change of variables
$$\Big(x,v,y\Big) \longrightarrow \Big(\yhypini(x,v),v,y\Big)\;,$$
 which preserves  the structure of the equations.
We further note that the hypotheses of Theorem~\ref{Tmnl} are
invariant under such coordinate transformations. From now on we
assume that this change of variables has been made.

\begin{Theorem}\label{Tmnl} Consider a set $\Omega$ as defined in
\eq{Omdef}, let $\alpha, \beta \in \R$, $k \in \N^*$,
$\delta=1/k$, and let $f=(\varphi, \psi) \in \cDai(\Omega)$ be a
solution of \eq{s1}. Suppose that
 \bel{HLdelta0} L_{ij}=
L_{ij}^{A}\partial_A + xL_{ij}^y
\partial_y + xL_{ij}^x\partial_x\,,
\ee 
with  \be
  L_{\varphi\varphi}^{\mu}\in x^\delta\Axyd\;, \quad L_{\psi\varphi}^{\mu}\;,
  L_{\varphi\psi}^{\mu}\;, L_{\psi\psi}^{\mu}\in \Axyd
\label{HLdelta} \ee (no symmetry hypotheses  are made on the
matrices $L^\mu_{ij}$), and {\rm
  \begin{deqarr}&
    B_{\varphi\varphi}\in C_\infty(\overline \Omega) +x^\delta\Axyd
\;,\qquad
B_{\varphi\psi}, B_{\psi\psi}, B_{\psi\varphi}\in \Axyd \;,&
\label{H3first}\arrlabel{H3}
\\ &a,b \in x^\beta\Axyd\;,&
\label{H3SEC}
\\ &\varphi|_\hyp =\zvarphi\in x^\beta\Axd\;,\
\psi|_\hyp =\zpsi\in x^\beta\Axd\;.&
\heqno\label{H4first}\arrlabel{H4}
  \end{deqarr}}Then
\bel{weakst}f\in x^\beta \Axyd+ y^{\beta} \Axyd+\Axd\;.\ee
In particular for any $\tau>0$ we have $f|_{\{y\ge \tau\}} \in
x^\beta \Axd+ \Axd$ , so that the solution is polyhomogeneous with
respect to $\{x=0\}$ on $\{y\ge\tau\}$.
\end{Theorem}

\begin{Remark}
We have  the more precise
statement \minilab{nwstst}
\begin{equs}\label{nwstst1}\psi &\in x^{\beta+1}
\Axyd+y^{\beta}\Ayd+ xy^\beta\Axyd + C_\infty(\overline\Omega) + x\Axd\;,\phantom{xxxxx}\\
 \varphi&\in x^\beta \Axd + x^\beta y \Axyd+y^{\beta+1}\Axyd+x \Axd  +y\Axd\;.
 \label{nwstst2} \end{equs}
Indeed, inserting \eq{weakst} into the equations and integrating one
obtains 
\beaa 
\psi &\in x^{\beta+1}
\Axyd+y^{\beta}\Ayd+ xy^\beta\Axyd + x\Axyd+y\Ayd\;, 
\\
 \varphi&\in x^\beta \Axd + x^\beta y \Axyd+x \Axd +y^{\beta+1}\Axyd +y\Axyd\;.
\eeaa
Comparing this with \eq{weakst} one obtains \eq{nwstst}.
\end{Remark}
\begin{Remark}\label{RLinfty}
  It is of interest to enquire when the solutions are in $L^\infty$. A
  direct analysis of the equations and \eq{nwstst} lead to the
  following: First, $\psi$ will be in $L^\infty$ whenever $\beta>-1$
  and $\mathring \psi\in L^\infty$. Next, consider the case $\beta =
  0$. If
$$L_{\varphi\psi}^\mu,B_{\varphi\psi},a,\zvarphi,\zpsi\in L^\infty\;,$$
then it also holds that \bel{linfconc} (\varphi,\psi)\in
L^\infty\;.\ee Finally, for $\beta>0$, it suffices to assume that
$L_{\varphi\psi}^\mu$ and $B_{\varphi\psi}$ are in $L^\infty$ to
achieve  \eq{linfconc}.
\end{Remark}

\begin{Remark}
In all physically relevant cases known to us we
have $\beta\in -\delta \N$, in which case the two last terms in
\eq{weakst} can  be absorbed in the first one, leading to the simpler formula
\bel{weakst2}f\in x^\beta \Axyd 
\;.\ee
\end{Remark}

\proof By hypothesis we can decompose $B_{\varphi\varphi}$ as
$$B_{\varphi\varphi}=\mathring B_{\varphi\varphi}+B_{\varphi\varphi}^\delta\;,
\quad \mathring B_{\varphi\varphi}\in C_\infty(\overline
\Omega)\;, \ B^\delta_{\varphi\varphi}\in x^\delta \Axyd\;.$$
 We rewrite the equations at hand as
\begin{eqnarray*}
\partial_y \varphi  + \mathring B_{\varphi\varphi}\varphi& = & c_\varphi\;,
\\
\partial_x \psi  & = & c_\psi\;,
\end{eqnarray*}
where
\begin{eqnarray*}
c_\varphi&:=& L_{\varphi\varphi}\varphi + L_{\varphi\psi}\psi + a
 -B_{\varphi\varphi}^\delta \varphi-B_{\varphi\psi}\psi\;,
\\c_\psi &:=&
L_{\psi\varphi}\varphi + L_{\psi\psi}\psi +
b-B_{\psi\varphi}\varphi -B_{\psi\psi}\psi\;.
\end{eqnarray*}
Integration of  the second equation yields \bel{psiint}
\psi(x,v^A,y) = \underbrace{\psi(y,v^A,y)}_{\mathring \psi(y,v^A)}
+ \int_y^{x} c_\psi(s,v^A, y)\,ds\;.\ee
By hypothesis we have $\zpsi(y,v^A)\in y^\beta\Ayd\subset
y^\beta\Axyd$. Similarly $\varphi$ can be calculated as
\bel{varphi} \varphi(x,v^A,y) =
R(x,v^A;y,x)\underbrace{\varphi(x,v^A, x)}_{\zvarphi(x,v^A)} +
\int_x^y R(x,v^A;s,x)c_\varphi(x,v^A,s) \,ds\;,\ee where
$R(x,v^A;y,y_1)$
 is the family of  resolvents  (smooth up to boundary in all variables) for the
family of ODEs
$$\partial_y\varphi (x,v^A,y)= -\mathring
B_{\varphi\varphi}(x,v^A,y)\varphi(x,v^A,y)\;,$$ with parameters
$(x,v^A)$ and with initial values at $y_1$. By hypothesis we have
$\zvarphi\in x^\beta\Axd$ which implies that the first term
$R\zvarphi$ is in $x^\beta \Axd\subset x^\beta\Axyd$.

In what follows we let $\epsilon>0$ be a positive constant, which
can be made as small as desired, and which may change from line to
line.

We wish, now, to show that
\bel{fsup} f\in \cDzim + x^\beta \Axyd+ y^{\beta} \Axyd\;.\ee
If $\alpha \ge 0$ there is nothing to prove, so assume that
$\alpha<0$. Decreasing $\alpha$ if necessary we can assume
that $\alpha+k\delta\ne -1$ for all $k\in \N$.
For all multi-indices $\gamma$ we have
\bel{stpoint2}
\partial_v^\gamma c_f\in \mcC^{\alpha-\epsilon}_{\{\zlxy \},\infty}+x^\beta \Axyd\ee
(here $\epsilon$ arises because of multiplicative powers of $\ln $
that might occur in some terms), in particular
$\partial_v^\gamma\partial_x \psi\in
\mcC^{\alpha-\epsilon}_{\{\zlxy \},\infty}+x^\beta \Axyd$, and by
\eq{psiint},
\bean \lefteqn{\psi(x,v^A,y) = \underbrace{\mathring \psi(y,v^A) +
\int_y^{x} b(s,v^A, y)\,ds}_{=:\psi_{0,\phg}\in x^{\beta+1}
\Axyd+y^{\beta}\Axyd}}&&\\ &&\qquad{}+ \int_y^{x} \underbrace{(
L_{\psi\varphi}\varphi + L_{\psi\psi}\psi -B_{\psi\varphi}\varphi
-B_{\psi\psi}\psi)(s,v^A,
y)}_{O(s^{\alpha-\epsilon})}\,ds\;.\eeal{psiintag} It follows that
\bel{endpoint0}\psi-\psi_{0,\phg}= O(x^{\alpha+1-\epsilon})+
O(y^{\alpha+1-\epsilon})\;,\ \mbox{similarly for} \
\partial_v^\gamma \psi\;.\ee
The derivatives of $\psi$ with respect to $y$ take the form:
\bel{psider}
\partial^\ell_y \psi(x,v^A,y) =
(\partial^\ell_x\zpsi)(y,v^A) + \int_y^x
\partial^\ell_y c_\psi(s,v^A,y)\,ds - \sum_{i=0}^{\ell-1}
C_{i\ell} (\partial^i_y \partial^{\ell-1-i}_x
c_\psi)(y,v^A,y)\;,\ee for some constants $C_{i\ell}\in \N$.
\Eq{psider} with $\ell=1$ gives
$$y\partial_y(\psi-\psi_{0,\phg})=
O(x^{\alpha+1-\epsilon})+ O(y^{\alpha+1-\epsilon})\;,\
\mbox{similarly for} \ y\partial_v^\gamma \partial_y\psi\;.$$ As
$\delta \le 1$, we have shown that
 \bel{psibet} \psi\in x^{\beta+1}
\Axyd+y^{\beta}\Axyd + \cDapmone\;. \ee

It now follows that
 $$c_\varphi = \underbrace{c_{\varphi,\phg}}_{\in x^{\beta}\Axyd + y^{\beta}\Axyd}
 +  \underbrace{c_{\varphi,\alpha+\delta-\epsilon}}_{O(x^{\alpha+\delta-\epsilon})+O(y^{\alpha+\delta-\epsilon})}\;,$$
 and \eq{varphi} gives
\beaa \varphi(x,v^A,y) &=&\underbrace{
R(x,v^A;y,x)\zvarphi(x,v^A)}_{\in x^\beta \Axd } +
\underbrace{\int_x^y R(x,v^A;s,x)c_{\varphi,\phg}(x,v^A,s)\,ds}_{
\in x^{\beta}\Axyd +y^{\beta}\Axyd}
\\ &&
+ \int_x^y
\underbrace{R(x,v^A;s,x){c_{\varphi,\alpha+\delta-\epsilon}}(x,v^A,s)}_{
O(x^{\alpha+\delta-\epsilon})+
 O(s^{\alpha+\delta-\epsilon})} \,ds\;,\eeaa
 leading to
  \bean\partial_v^\gamma (\varphi-\varphi_{0,\phg})&=&
O(yx^{\alpha+\delta-\epsilon})+O(x^{\alpha+\delta+1-\epsilon})+O(y^{\alpha+\delta+1-\epsilon})
\\ &=& O(x^{\min(\alpha+\delta-\epsilon,0)}) \ \mbox{ for }\ 0\le x \le y \;,\eeal{endpoint}
with
 \bel{varphiphg}\varphi_{0,\phg}(x,v^A,y) := R(x,v^A;y,x)\zvarphi(x,v^A)
+ \int_x^y R(x,v^A;s,x){c_{\varphi,\phg}}(x,v^A,s)\;.\ee
Equations~\eq{psibet} and \eq{endpoint} in particular imply
\bel{finaleq} f\in x^\beta \Axyd+y^{\beta}\Axyd + \cDmz\;.\ee

We note the following:

\begin{Lemma}\label{L4.4}
Under  the hypotheses of Theorem \ref{Tmnl}, the fields
$$\tilde\varphi := \left(\begin{array}{c} \varphi \\\partial_{A}\varphi \\
x\partial_x\varphi \\y\partial_y \varphi
 \end{array}\right)\;,\qquad \tilde\psi := \left(\begin{array}{c}
\psi \\\partial_A \psi \\ x\partial_x\psi \\y\partial_y \psi
 \end{array}\right)$$ satisfy a system of equations of
 the form \eq{s1}, with coefficients $\tilde L_{ij}$, $\tilde
 B_{ij}$ and sources $\tilde a$, $\tilde b$ satisfying the hypotheses of Theorem \ref{Tmnl}, with
 $\tilde\varphi,\tpsi \in \cDai$, and with initial data
 $\tilde\varphi|_\hyp,\tpsi |_\hyp\in x^\beta\Axd$.
\end{Lemma}
\proof  We rewrite the equations \eq{s1} as
\beaa \partial_y\varphi & = & \mathring c_\varphi \;,\\
\partial_x\psi & = & c_\psi\;.
\eeaa Differentiating  $(\tilde\varphi,\tilde\psi)$ one
gets
\bel{seq} \partial_y\tilde\varphi =
\left(\begin{array}{c}
\mathring c_\varphi \\\partial_A\mathring c_\varphi \\ x\partial_x \mathring c_\varphi \\
  \mathring c_\varphi + y\partial_y \mathring c_\varphi
 \end{array}\right)\;,\qquad
\partial_x\tilde\psi =
\left(\begin{array}{c} c_\psi \\\partial_A c_\psi \\ c_\psi + x\partial_x c_\psi \\
y\partial_y c_\psi
\end{array}\right)\;.\ee
With a little work one shows that the right-hand sides of these
equations are first order differential operators of appropriate
form acting on $(\tilde\varphi,\tilde\psi)$. For example, consider
the term $y\partial_y \mathring c_\varphi$:
\beaa y\partial_y \mathring c_\varphi &=&  y\partial_y
(L_{\varphi\varphi}\varphi + L_{\varphi\psi}\psi + a
 -B_{\varphi\varphi} \varphi-B_{\varphi\psi}\psi) \\
 &=& (y\partial_y L^A_{\varphi\varphi})\partial_A\varphi
 + (x\partial_y L^y_{\varphi\varphi})y\partial_y\varphi
 + (y\partial_y L^x_{\varphi\varphi})x\partial_x\varphi \\
 && {} +  L^A_{\varphi\varphi}\partial_A(y\partial_y \varphi)
 + L^y_{\varphi\varphi}x\partial_y(y\partial_y\varphi
 -\varphi)
 + L^x_{\varphi\varphi}x\partial_x(y\partial_y\varphi)\\
 &&{}+(\mbox{analogous } \partial \psi \mbox{ terms})\\
 &&{}+y\partial_y a- B_{\varphi\varphi} (y\partial_y\varphi)
 -B_{\varphi\psi}(y\partial_y\psi) \\
 && {}- (y\partial_yB_{\varphi\varphi})
 \varphi-(y\partial_yB_{\varphi\psi})\psi\;.\eeaa
The terms from the first two lines of the second equality above
can be included in $\tilde L_{\tilde\varphi\tilde\varphi}$, the
third line can be included in $\tilde L_{\tvarphi  \tpsi }$,
$\tilde a$ will include $y\partial_y a$, whatever remains going
into the $\tilde B_{ij}$ terms. We note that $\Axyd$ and
$x^\delta\Axyd$ are invariant under $y\partial_y$ and
$x\partial_x$, which makes it simple to check that these choices
lead to the functional spaces as required in
\eq{HLdelta}--\eq{H3}.

The other terms in \eq{seq} are treated analogously.

It remains to show that the initial data are in the right space.
This is obvious for $\varphi$, $\psi$, as well as their derivatives
with respect to $v^A$. For the remaining derivatives, we pass to the
limit $x\to y$ in \eq{psider} with $\ell=1$ to obtain
\begin{deqn}y\partial_y \psi|_{\hyp}= y(\partial_x \zpsi)
-yc_\psi|_\hyp\;.\label{eqders1}\arrlabel{eqders}\end{deqn}
Similarly from
$$ \varphi(x,v^A,y) = \underbrace{\varphi(x,v^A,x)}_{\mathring
\varphi(x,v^A)} + \int_x^{y} \mathring c_\varphi(x,v^A, s)\,ds\;.$$
we find (again for $x=y$)
\begin{ddeqn}y\partial_x \varphi|_{\hyp}= y(\partial_x
\zphi)-y\mathring
c_\varphi|_\hyp\;.\label{eqders2}\end{ddeqn}
Equations~\eq{s1} further give
\begin{ddeqn}y(\partial_x\psi)|_\hyp = yc_\psi|_\hyp\;,\end{ddeqn}
 \begin{ddeqn}y(\partial_y\varphi)|_\hyp = y\mathring c_\varphi|_\hyp\;.\end{ddeqn}
The terms $y(\partial_x \zpsi)(y,v^A)$ and $y(\partial_x
\zphi)(y,v^A)$ in \eq{eqders1}-\eq{eqders2} are in $y^\beta \Ayd$.
Now,
$$y \mathring c_\varphi|_\hyp=y(L_{\varphi\varphi}\varphi +
L_{\varphi\psi}\psi + a
 -B_{\varphi\varphi} \varphi-B_{\varphi\psi}\psi)|_\hyp\;.$$
 The restrictions to $\hyp$ of the terms $a$, $B_{\varphi\varphi} \varphi$,
 $B_{\varphi\psi}\psi$ and the derivatives of $\varphi$ and $\psi$ with respect to
 $v^A$, give a contribution which is in  $y^\beta \Ayd$. The remaining
 terms are of the form $y(\partial_y\psi)|_\hyp$,
 $y(\partial_y\varphi)|_\hyp$, $y(\partial_x\psi)|_\hyp$,
 $y(\partial_x\varphi)|_\hyp$ multiplied by coefficients from
 $\Ayd$. The same analysis applies to $y
 c_\psi|_\hyp$, so that we can write the system of equations \eq{eqders}
 as
$$(\mathrm{Id}-yK)\left(\begin{array}{c} y(\partial_y\psi)|_\hyp\\
y(\partial_y\varphi)|_\hyp\\ y(\partial_x\psi)|_\hyp\\
y(\partial_x\varphi)|_\hyp\end{array}\right) \in y^\beta \Ayd\;.$$
Here $K$ is a  matrix with components in $\Ayd$. There exists
$\epsilon>0$ so that for $0\le y<\epsilon$ the matrix
$\mathrm{Id}-yK$ has an inverse in $\mathrm{Id}+y\Ayd$, and
polyhomogeneity (with appropriate power structure) of the initial
data for $(\tvarphi ,\tpsi )$ follows.
 \qed

Applying \eq{endpoint0} and \eq{endpoint} to $(\tvarphi ,\tpsi )$
one finds that
$$(\tilde\varphi,\tilde\psi)\in \cDmz+ x^\beta \Axyd+y^{\beta}\Axyd\;,$$
and Proposition~\ref{Pequiv4} implies
$$(\varphi,\psi)\in \cDmonem+ x^\beta \Axyd+y^{\beta}\Axyd\;.$$ Continuing in this way we
are led to \bel{realend}(\varphi,\psi)\in \cDmim+x^\beta
\Axyd+y^{\beta}\Axyd\;.\ee
Repeating the whole argument $\ell$ times, where $\ell$ is the
smallest number such that $\alpha+\ell \delta >0$, one arrives at
\bel{fsup2} f\in \cDzim + x^\beta \Axyd+y^{\beta}\Axyd\;.\ee

Suppose, for the sake of induction, that \bel{begind}
f=\underbrace{f_{1,k}}_{\in x^\beta
\Axyd+y^{\beta}\Axyd}+\underbrace{f_{2,k}}_{\in
\isAxf{k\delta-\epsilon}}+\underbrace{f_{3,k}}_{\in
\fffT{k\delta-\epsilon}00}
 \;; \ee
\Eq{fsup2} and the embedding $$\cDzim \subset \fffT{-\epsilon}00$$
show that this equation holds with $k=0$ and
$f_{2,0}=0 $. 
We note that all of the above spaces are invariant with respect to
$x\partial_x$, $x\partial_y$, $\partial_v$ and with respect to
multiplication by a function from $\Axyd$.

Integrating the equation for $\psi$ and using
Propositions~\ref{calkA}, \ref{PIoneT} and \ref{PIoneF} one finds
\bean
\psi(x,v^A,y)& =& \underbrace{\psi_{0,\phg}(x,v^A,y)+ I_1
\Big(L_{\psi\varphi}\varphi_{1,k} + L_{\psi\psi}\psi_{1,k}
-B_{\psi\varphi}\varphi_{1,k}
-B_{\psi\psi}\psi_{1,k}\Big)}_{=:\psi_{1,k+1}\in x^{\beta+1}
\Axyd+y^{\beta}\Axyd}\\\nonumber & &
+\underbrace{I_1\Big(L_{\psi\varphi}\varphi_{2,k} +
L_{\psi\psi}\psi_{2,k} -B_{\psi\varphi}\varphi_{2,k}
-B_{\psi\psi}\psi_{2,k}\Big)}_{\in \isAxf{k\delta+1-\epsilon}\subset \isAxf{(k+1)\delta-\epsilon}}\\
& & +\underbrace{I_1 \Big(L_{\psi\varphi}\varphi_{3,k} +
L_{\psi\psi}\psi_{3,k} -B_{\psi\varphi}\varphi_{3,k}
-B_{\psi\psi}\psi_{3,k}\Big)}_{\in \fzTi {k\delta+1-\epsilon} +
\fffT {k\delta+1-\epsilon} 00\subset \fzTi {(k+1)\delta-\epsilon}
+ \fffT {(k+1)\delta-\epsilon} 00} \;,
 \eeal{begind48}
showing that the result is true for $\psi$ with $k$ replaced by
$k+1$. Inserting this information into \eq{varphi} one similarly
finds, using Propositions~\ref{calkA} and \ref{PItwo}, that
\eq{begind} with $k$ replaced by $k+1$ holds for $\varphi$ (here
$R$ stands for the
--- smooth up to the boundary --- resolvent $R(x,v^A;y,x)$):
\bean
\varphi(x,v^A,y)& =& \underbrace{\varphi_{0,\phg}(x,v^A,y)+ I_2
\Big[ R\cdot \big(L_{\varphi\varphi}\varphi_{1,k} +
L_{\varphi\psi}\psi_{1,k} -B^\delta_{\varphi\varphi}\varphi_{1,k}
-B_{\varphi\psi}\psi_{1,k}\big)\Big]}_{=:\varphi_{1,k+1}\in
x^\beta
\Axyd+y^{\beta}\Axyd}\\
& & +\underbrace{I_2 \Big[R\cdot\big(
L_{\varphi\varphi}\varphi_{2,k} + L_{\varphi\psi}\psi_{2,k}
-B^\delta_{\varphi\varphi}\varphi_{2,k}
-B_{\varphi\psi}\psi_{2,k}\big)\Big]}_{\in \isAxf{(k+1)\delta}}
 \nonumber
 \\
& & +\underbrace{I_2
\Big[R\cdot\big(\underbrace{L_{\varphi\varphi}\varphi_{3,k}
-B^\delta_{\varphi\varphi}\varphi_{3,k} +
L_{\varphi\psi}\psi_{3,k} -B_{\varphi\psi}\psi_{3,k}}_{\in \fffT
{(k+1)\delta-\epsilon} 00}\big)\Big]}_{\in \fffT
{(k+1)\delta-\epsilon} 11\subset \fffT {(k+1)\delta-\epsilon} 00}
\;.\eeal{begind50} We have thus shown that \eq{begind} holds for
all $k$. For any $m\in \N$ we can choose $k$ large enough so that
the last term there is in $C_m(\overline \Omega)$, and that all
the coefficients of an expansion of $f_{2,k}$ in terms of powers
of $x$ and $\log x$ also are in $C_m(\overline \Omega)$. The
result follows now by arguments essentially identical to those of
Proposition~\ref{Pequiv}. \qed

\begin{Remark}
  An identical proof shows that ``any'' polyhomogeneous initial data,
  for equations with ``any'' polyhomogeneous and bounded coefficients,
  and ``any'' polyhomogeneous source terms lead to polyhomogeneous
  solutions, without imposing corner conditions. More precisely, let
  us assume, instead of \eq{HLdelta}, that there exists $\lambda>0$
  such that \be L_{\varphi\varphi}^{\mu}\in x^\lambda\Axyp\;, \quad
  L_{\psi\varphi}^{\mu}\;, L_{\varphi\psi}^{\mu}\;,
  L_{\psi\psi}^{\mu}\in \Axyp\;.  \ee By $\Axyp$, $\Axp$, etc.~we mean
  here spaces of polyhomogeneous functions for which $n_i\ge0$ and
  $n_i+\hat{n}_i\ge 0$, where $n_i$ and $\hat{n}_i$ denote the powers
  of, respectively, $x$ and $y$ in the polyhomogeneous expansion. Let
  us further assume that, for some $\lambda'>0$, \be
  B_{\varphi\varphi}\in C_\infty(\overline \Omega) +
  x^{\lambda'}\Axyp\;, \quad B_{\varphi\psi}, B_{\psi\psi},
  B_{\psi\varphi}\in \Axyp \;,\ee
  \be a \in x^{\beta_1} \Axyp\;,\quad b \in x^{\beta_2} \Axyp \;,\ee
  \be \varphi|_\hyp =\zvarphi\in x^{\beta_3}\Axp\;,\ \psi|_\hyp
  =\zpsi\in x^{\beta_4}\Axp\;. \ee We have then the following version
  of the result: \be f\in y^\beta \Axyp+x^\beta \Axyp+\Axp\;,\ee where
  $\beta = \min(\beta_1+1, \beta_2, \beta_3, \beta_4)$. The conclusion
  of Remark~\ref{RLinfty} remains valid if $\beta$ there is replaced
  with $\min(\beta_1+1, \beta_2, \beta_3, \beta_4)$.
\end{Remark}

\subsection{Nonlinear equations} \label{Snepcc}

 In this section
we prove the non-linear equivalent of Theorem~\ref{Tmnl}. The method
of proof is essentially identical, when the nonlinearities are
treated appropriately. We consider semi-linear systems of the form
\eq{s1} with a supplementary non-linearity at the right-hand-side:
{\rm
\begin{deqarr}
  \partial_y\varphi + B_{\varphi\varphi}\varphi
  + B_{\varphi\psi}\psi&=& L_{\varphi\varphi}\varphi +
  L_{\varphi\psi}\psi + a+G_\varphi
  \;,\label{s1an}\\ \partial_x \psi+
  B_{\psi\varphi}\varphi + B_{\psi\psi}\psi &=&   L_{\psi\varphi}\varphi +
  L_{\psi\psi}\psi + b+G_\psi \label{s1bn}\;,
\arrlabel{s1n}\end{deqarr}}It should be clear from the discussion at
the beginning of Section~\ref{Spwcc} that the right hypothesis for a
consistent set-up as in Theorem~\ref{Tmnl} is to assume
$\beta=-1+\delta$, and to have a non-linearity compatible with this
hypothesis. This leads us to assume, following~\cite{ChLengardnwe},
that the nonlinearity
$$G= (G_\varphi,G_{\psi_1},G_{\psi_2})$$
takes the form
\begin{equation}
  \label{m0}
  G= x^{-p\delta} H(x^\mu,x^{q\delta} \psi_1, x^{q\delta+1} \psi_2, x^{q\delta+1} \varphi)
\;.
\end{equation}
Here we have decomposed $\psi$ as
\be\label{psidecomp}\psi=\left(\begin{array}{c} \psi_1\\ \psi_2
\end{array}\right)\;; \ee this is motivated by different {\em a
priori\/} estimates we have at our disposal for the appropriately
defined components $\psi_1$ and $\psi_2$ of $\psi$ in the
applications we have in mind. Polyhomogeneity of solutions of
nonlinear wave equations, or of the wave map equation, follows as
a special case, see the proofs of \thm{T2phg} and of
\thm{Twavemap}.

We will need to impose various restrictions on the function $H$,
in order to do that some terminology will be needed. We shall say
that a function $H(x^\mu,\argw)$ is {\em $\Axyd$-polyhomogeneous
in $x^\mu$ with a uniform zero of order $l$ in $\argw$} if the
following hold: First, $H$ is smooth in $\argw\in\R^N$ at fixed
$x^\mu$. Next, it is required that for all $B\in \R$ and $k\in \N$
there exists a constant $\hat{C}(B)$ such that for all $|\argw|\le
B$ and $0\leq i\leq \min(k,l)$,
\be\label{S2.71} \left\|\frac{\partial^iH(\cdot,\argw)}{\partial
    \argw^i}\right\|_{{\mcC}_{k-i}^0(\Omega)} \leq \hat{C}(B)
|\argw|^{l-i}\;.\ee Further,
\begin{equation}
  \label{eq:m1}
\forall i\in\N \quad \partial^i_\argw H(\cdot, \argw)\in\Axyd
\end{equation}
at fixed constant $\argw$. Finally we demand the uniform estimate
for constant $\argw$'s
\begin{equation}
  \label{m1}
  \forall \epsilon >0,  M \geq 0 , i,k\in \N \ \; \exists\;
C(\epsilon,M,i,k) \ \forall |\argw|\leq M \quad \|
\partial^i_\argw H(\cdot, \argw)\|_{\mcC^{-\epsilon}_k(\Omega)} \leq
C(\epsilon,M,i,k)\;.
\end{equation}
The qualification ``in $\argw$'' in ``uniform zero of order $l$ in
$\argw$'' will often be omitted. The small parameter $\epsilon$
has been introduced above to take into account the possible
logarithmic blow-up of functions in $\Axyd$ at $x=0$; for the
applications to the nonlinear scalar wave equation or to the wave
map equation on Minkowski space-time, the alternative simpler
requirement would actually suffice:
\begin{equation}
  \label{m2}
  \forall  \; M \geq 0 , i,k\in \N \ \; \exists\;
C(M,i,k) \ \forall |\argw|\leq M \quad \| \partial^i_\argw
H(\cdot, \argw)\|_{\mcC^{0}_k(\Omega)} \leq C(M,i,k)\;,
\end{equation}
again for constant $\argw$'s. Clearly functions which are smooth
in $(\argw, x^\mu)$, and have a zero  of appropriate order in
$\argw$ at $\argw=0$, satisfy the above conditions.

\begin{Theorem}
  \label{Tlemme1n}
  Let $p\in \Z,\; q, 1/\delta\in \N^*$, $k\in \N\cup\{\infty\}$, let
  $f=(\psi=(\psi_1,\psi_2),\varphi)$ be a solution of \eq{s1n} with
  $G$ of the form \eq{m0}, where $H$ is $\Axyd$-polyhomogeneous in
  $x^\mu$ with a uniform zero of order
  \begin{equation}
    \label{eq:ord}
    m > \frac{p- \frac 1 \delta}{q}\;.
  \end{equation}%
Assume that for every $\epsilon>0$ we have
$$
\psi_1 \in \Big(\fCxy {-1-\epsilon}\infty\cap
\fCxy{-\epsilon}0\Big)(\Omega)\;,$$
$$\varphi,\psi_2 \in \fCxy {-1-\epsilon}\infty (\Omega)\;.$$
Suppose
that Equations~\eq{HLdelta0}-\eq{H4} hold with $\beta=-1+\delta$, and
that $\psi|_{x=y}\in \Axd$.
 Then $$(\psi,\varphi)\in \Axyd\times x^{\delta-1}\Axyd\;;$$ more precisely
 \minilab{stsn}\begin{equs}\label{stst1n}\psi &\in x^{\delta} \Axyd+ \Ayd\;,\\
   \varphi&\in x^{\delta-1} \Axd + x^{\delta-1} y
   \Axyd\;.\label{stst2n} \end{equs} In particular for any $\tau>0$ we
 have $$(\psi,\varphi)|_{\{y\ge \tau\}} \in
 \Axd\times x^{\delta-1} \Axd\;,$$
 which shows that the solution is
 polyhomogeneous with respect to $\{x=0\}$ on
 $\{y\ge\tau\}$.
\end{Theorem}

\begin{Remark}\label{Rldcb}
  If $\psi|_{x=y}\in \Axd\cap L^\infty$, then $\psi \in x^{\delta}
  \Axyd+\Ayd\cap L^\infty $.
\end{Remark}

\begin{Remark}\label{Rldcc}
  Let us, in addition to $\psi\in L^\infty(\overline{\Omega})$ (see
  Remark \ref{Rldcb}), assume that
$$L_{\varphi \psi}^\mu,B_{\varphi \psi},a,\varphi|_{x=y}
\in L^\infty(\overline{\Omega})\;.$$ Suppose further that
$G_\varphi$ satisfies
$$G_\varphi (x^\mu, x^{q\delta}\psi_1, x^{q\delta+1}\psi_2,
x^{q\delta+1}\varphi) \in L^\infty(\overline{\Omega})$$for every
(fixed) $\psi, \varphi$ as in \eq{stsn} (which is true \emph{e.g.} if
$(mq-p)\delta>0$). Then it holds that
$$\varphi\in
\Axyd\cap L^\infty(\overline{\Omega})\;.$$
\end{Remark}

\begin{Remark}\label{Rldc} There is little doubt that the proof below can be generalised to
allow $\psi_2\in x^{-1+\delta}\Axd$, as well as $B_{\psi\psi}\in
x^{-1+\delta}\Axd$. We have not attempted to check all the details
of this, because of technical complications in the proof.
\end{Remark}


 \remark The theorem remains true if we replace $G$ by
a finite sum of nonlinearities satisfying the above hypotheses,
with different $p$ and $q$ satisfying \eq{eq:ord} for each term in
the sum.

\medskip

\proof We follow the proof of Theorem~\ref{Tmnl}, checking that
the supplementary non-linear terms do not affect the argument. We
decompose $c_\psi$ as $(c_{\psi_1},c_{\psi_2})$.

We have
$$\tf :=(\psi_1,x\psi_2,x\varphi)\in \fCxy {-\epsilon}0\;,$$
with $\epsilon$ as small as desired (we hope that a clash of
notation $\tf$ with Section~\ref{Sweci} will not confuse the
reader). This shows that the non-linearity $G_\psi$ gives a
contribution to $c_\psi$ which is $O(x^{(mq-p)\delta-\epsilon})$
(the easiest way to see that is to view the non-linearity as a
function of $x^{q\delta-\epsilon}\tf_\epsilon$, where $\tf_\epsilon
= x^\epsilon \tf\in L^\infty$).

First,  \eq{endpoint0} becomes
 \bel{endpoint0n}\psi-\psi_{0,\phg}=
     O(x^{-\epsilon})\;. \ee
     Next, to estimate the derivatives of $\psi$, note that the integrand  for $\partial_y \psi$ in
\eq{psider} with $\ell=1$ contains a supplementary term
 $$\partial_y \Big( G_{\psi}(x^\mu,x^{q\delta}\tf )\Big)
 =\underbrace{(\partial_y G_{\psi})(x^\mu,x^{q\delta}\tf )}
 _{O(x^{(mq-p)\delta-\epsilon})y^{-1}}
  +(\partial_\argw G_{\psi})(x^\mu,x^{q\delta}\tf )
  x^{q\delta}\partial_y \tf \;. $$
Write
$\argw=(\argw_1,\argw_2,\argw_\varphi)=(x^{q\delta}
\psi_1,x^{q\delta+1} \psi_2,x^{q\delta+1} \varphi)$,
$\argw_\psi=(\argw_1,\argw_2)$. We have the equation
\bean(\partial_\argw
G_{\psi})(x^\mu,x^{q\delta}\tf )x^{q\delta}\partial_y \tf
  &=& x^{-p\delta}
\Big(\partial_{\argw_2} H_{\psi}\Big)(x^\mu,\argw)
x^{q\delta+1}\partial_y{\psi_2}
 \\
 \nonumber
 && +x^{-p\delta}
\Big(\partial_{\argw_\varphi} H_{\psi}\Big)(x^\mu,\argw)
x^{q\delta+1}\partial_y{\varphi}
 \\
 && +x^{-p\delta}
\Big(\partial_{\argw_1}
H_{\psi}\Big)(x^\mu,\argw)x^{q\delta}\partial_y
{\psi_1}\;. \eeal{Gder}
We know by hypothesis that $$(\psi_1,
x\psi_2,x\varphi,xy\partial_y\psi_2,xy\partial_y\varphi)\in
\fCxy{-\epsilon}0\;,$$ which allows us to estimate all terms above
by $C x^{(mq-p)\delta-\epsilon}y^{-1}$, \emph{except the last one}.
However, we can write a system of equations for $y\partial_y\psi$ of
the form:
 \bea {\partial_x (y\partial_y\psi) -
\underbrace{x^{-p\delta} \Big(\partial_{\argw_1}
H_{\psi}\Big)(x^\mu,\argw)x^{q\delta}}_{=O(x^{(mq-p)\delta-\epsilon})}y\partial_y
{\psi_1}
}=  O(x^{-1-\epsilon})\;. \eqname{3.45a}
 \eeal{odecas}
The following is essentially Proposition~B.3 of
\cite{ChLengardnwe} with $k=0$ there, we revisit the proof because
of the uniformity in $y$ required here:

\begin{Lemma}
\label{LoddEST} For $0\le x\le y\le \mathring y$ consider the
system
$$\partial_x \psi + b \psi = c\;,$$
and suppose that there exists $\epsilon<1$ such that the linear
map $b$ has coefficients in  ${\fCx{-\epsilon}0}$. For $\alpha\in
\R\setminus\{-1\}$ there exists a constant
$C=C(\alpha,\epsilon,\|b\|_{\fCx{-\epsilon}0},\mathring y )$ such
that
\begin{enumerate} \item For $\alpha>-1$ we have
\bel{firstinLodd}\|\psi\|_{L^\infty} \leq C\left (\|\psi|_{x=y}\|_{L^\infty}+
\|c\|_{\fCx{\alpha}0}\right)\;,
 \ee
\item while for $\alpha<-1$ it holds that
 \bel{firstinLodd2}\|\psi\|_{\fCx {\alpha+1}0} \leq C\left (\|\psi|_{x=y}\|_{\fCx {\alpha+1}0}+
\|c\|_{\fCx{\alpha}0}\right)\;.
 \ee
\end{enumerate}
\end{Lemma}

\proof 1. Obvious modifications in the calculations in the proof
of~\chcite{Proposition~B.3} lead to the following replacement of
\chleq{B.19}
\beqa -\partial_xe_a &\leq &
x^{\alpha}\|c\|_{\fCx{\alpha}k} +x^{\- \epsilon}
\|b\psi\|_{\fCx{-\epsilon}k} \;.\label{ineqbase}
 \eeqa
 By integration, after  passing to the limit $a\to 0$, one obtains,
 for $0\le x \le y \le y_1 \le \mathring y$
$$\|\psi\|_{L^\infty([0,y_1])} \leq \|\psi|_{x=y}\|_{L^\infty}+
\frac{\mathring y^{\alpha+1}}{|\alpha+1|}\|c\|_{\fCx{\alpha}k}
+\frac{y_1^{1- \epsilon} }{1-\epsilon}\|b\|_{\fCx{-\epsilon}k}
\|\psi\|_{L^\infty([0,y_1])} \;.
$$
Choosing $y_1$ small enough the last term can be carried over to
the left-hand-side, which yields the inequality on $[0,y_1]$. The
final inequality \eq{firstinLodd} is then standard.

2. Decreasing $\epsilon$ slightly if necessary, we can without
loss of generality assume that
$\alpha+1-\epsilon\ne -1$. We recall the first two lines of
\chleq{B.20}, adapted to our notation, with $x_3$ there replaced
by $y$:
\beqa
 |\psi(x,v^A,y)|&\leq& |\psi(y,v^A,y)| + \left(-{x^{\alpha+1} \over
(1+\alpha)} + {y^{\alpha+1}\over (1+\alpha)}\right)
\|c\|_{\fCx \alpha 0} \nonumber\\
&&+\left({y^{\alpha+2-\epsilon} \over (2+\alpha-\epsilon)} -
{x^{\alpha+2-\epsilon} \over (2+\alpha-\epsilon)}\right)
\|b\|_{\fCx{-\epsilon}0}\|\psi\|_{\fCx {\alpha+1} 0} 
 \;.\phantom{xxxxx}\label{lasteq} \eeqa
 For $\alpha<-1$ this
implies
\beqa
 |\psi(x,v^A,y)|&\leq& \|\psi(y,v^A,y)\|_{\fCx {\alpha+1}0} \underbrace{y^{\alpha+1}}_{\le x^{\alpha+1}}  + {x^{\alpha+1} \over
|1+\alpha|}
\|c\|_{\fCx \alpha 0} \nonumber\\
&&+{1 \over
|2+\alpha-\epsilon|}\Big(\underbrace{y^{\alpha+2-\epsilon}}_{=y^{1-\epsilon}
y^{\alpha+1}\le y_1^{1-\epsilon} x^{\alpha+1}} +
\underbrace{x^{\alpha+2-\epsilon}}_{=x^{1-\epsilon}
x^{\alpha+1}\le y_1^{1-\epsilon} x^{\alpha+1}} \Big)
\|b\|_{\fCx{-\epsilon}0}\|\psi\|_{\fCx {\alpha+1} 0}
 \nonumber
\\&\leq& \Big(\|\psi(y,v^A,y)\|_{\fCx {\alpha+1}0}  + {1
\over |1+\alpha|} \|c\|_{\fCx \alpha 0} \nonumber \\ && +{2
y_1^{1-\epsilon}\over |2+\alpha-\epsilon|}
\|b\|_{\fCx{-\epsilon}0}\|\psi\|_{\fCx {\alpha+1} 0}\Big)x^{\alpha+1} 
 \;.\phantom{xxxxx}\label{lasteq2} \eeqa
 Choosing $y_1$ sufficiently small, one concludes as before.
 \qed

{}From \eq{odecas} and Lemma~\ref{LoddEST} with $\alpha=-1-\epsilon$
we obtain
$$y\partial_y\psi=
     O(x^{-\epsilon})\;. $$
A similar treatment of $\partial_A G_\psi$ allows one to conclude
that
\bel{psibetn} \psi\in
      x^{\delta} \Axyd+\Ayd +\fCxy {-\epsilon}1\;. \ee
      (The first two spaces are actually included in the
      last one, but we exhibit them to  keep
      track of the form of the polyhomogeneous contributions in
      $\psi$, which will become important later.)
It now follows that
 $$c_\varphi +G_\varphi = \underbrace{c_{\varphi,\phg}}_{\in x^{\beta}\Axyd + y^{\beta}\Axyd}
 +  {\cprim}\;,$$
 with
$$ {\cprim}=
     O(x^{-1+\delta-\epsilon})\;. $$
\Eq{varphi} gives
\beaa \varphi(x,v^A,y) &=&\underbrace{
R(x,v^A;y,x)\zvarphi(x,v^A)}_{\in x^{\delta-1} \Axd }
+\underbrace{ \int_x^y R(x,v^A;s,x)c_{\varphi,\phg}(x,v^A,s)
\,ds}_{ \in x^{\delta}\Axd +x^{\delta-1}y\Axyd}
\\ &&
+ \int_x^y {R(x,v^A;s,x){\cprim}(x,v^A,s)} \,ds\;,\eeaa
 leading to
 \bean\varphi-\varphi_{0,\phg}&=&
     O(x^{-1+\delta-\epsilon})\;.
 \eeal{endpointn}
Here, as before, $\varphi_{0,\phg}$ is given by \eq{varphiphg}.

 Next, consider Lemma~\ref{L4.4}.
We rewrite \eq{s1} as
\beaa \partial_y\varphi & = & \mathring c_\varphi +G_\varphi \;,\\
\partial_x\psi & = & c_\psi+G_\psi\;.
\eeaa The field $\tpsi $ is decomposed as $(\tpsi _1, \tpsi _2)$,
with $\tpsi _1 =(\psi_1,\partial_A \psi_1, x\partial_x \psi_1,
y\partial_y \psi_1)$. Differentiating
$(\tilde\varphi,\tilde\psi)$, instead of \eq{seq} one gets
\bel{seqn} \partial_y\tilde\varphi =
\left(\begin{array}{c} \mathring c_\varphi +G_\varphi
\\\partial_A\mathring c_\varphi
          +\partial_A G_\varphi
          \\ x\partial_x \mathring c_\varphi +x\partial_x G_\varphi
          \\
  \mathring c_\varphi+G_\varphi + y\partial_y \mathring c_\varphi+y\partial_y G_\varphi
 \end{array}\right)\;,
 \ee
similarly for $\partial_x\tpsi $, \emph{e.g.},
\bel{seqn2}
 \partial_x\tilde\psi_1 =
 \left(\begin{array}{c} c_{\psi_1}+G_{\psi_1}
 \\\partial_A c_{\psi_1}+\partial_A G_{\psi_1}
 \\ c_{\psi_1}+G_{\psi_1}  + x\partial_x c_{\psi_1}+x\partial_x G_{\psi_1}\\
y\partial_y c_{\psi_1}+y\partial_yG_{\psi_1}
\end{array}\right)\;.
 \ee
One can check that the non-linear terms above
 have the structure claimed. For example, we have
\beaa  x\partial_x G_{\psi_1} & = &  x\partial_x \Big(x^{-p\delta}
H_{\psi_1}(x^\mu,\argw)\Big)
\\ & = &  -px^{-p\delta}
H_{\psi_1}(x^\mu,\argw)  + x^{1-p\delta} \Big(\partial_x
H_{\psi_1}\Big)(x^\mu,\argw)
 \\ && +x^{-p\delta}
\Big(\partial_{\argw_1} H_{\psi_1}\Big)(x^\mu,\argw)\Big(q\delta
\,x^{q\delta}{\psi_1}
 + x^{q\delta}x\partial_x {\psi_1}\Big)\\
  && +x^{-p\delta}
\Big(\partial_{\argw_2}
H_{\psi_1}\Big)(x^\mu,\argw)\Big((q\delta+1)x^{q\delta+1}{\psi_2}
 + x^{q\delta+1}x\partial_x{\psi_2}\Big)\\ && +x^{-p\delta}
\Big(\partial_{\argw_\varphi}
H_{\psi_1}\Big)(x^\mu,\argw)\Big((q\delta+1)x^{q\delta+1}{\varphi}
 + x^{q\delta+1}x\partial_x{\varphi}\Big)\;. \eeaa
It should be clear that each term in the sum above has a zero of
order $m$ in the new variables.

 This analysis allows
us to repeat the induction argument which led to \eq{realend},
obtaining instead
\bean \psi&\in&x^{\delta}
\Axyd+\Ayd+ \cDmimnp\;,
 \\ \varphi&\in& x^{\delta-1}
\Axd+x^{\delta-1}y\Axyd +
     \fCxy{-1+\delta-\epsilon}\infty\;.
 \eeal{realend2}
%
The embeddings
 $$\cDzim \subset \fffT{-\epsilon}00\;,\qquad \cDzimne \subset \fffT{-1+\delta-\epsilon}00\;,$$
 justify the case $k=0$ of  the following
induction hypothesis, where $f_{2,0}=(\varphi_{2,0},\psi_{2,0})=0$
(we hope that the reader will not get confused by a slight clash
of notation, as the decomposition  in the equations below is
unrelated to the decomposition $\psi=(\psi_1,\psi_2)$):
\beal{begindn1} \varphi&=&\underbrace{\varphi_{1,k}}_{\in x^{\delta-1}
\Axd+x^{\delta-1}y\Axyd}+\underbrace{\varphi_{2,k}}_{\in
\isAxf{k\delta-1+\delta-\epsilon}}+\underbrace{\varphi_{3,k}}_{\in
\fffT{k\delta-1+\delta-\epsilon}00}
 \;,\phantom{xxxxxxxxxx}
 \\\label{begindn2}
 \psi&=&\underbrace{\psi_{1,k}}_{\in x^{\delta}
\Axyd+\Ayd}+\underbrace{\psi_{2,k}}_{\in
\isAxf{k\delta-\epsilon}}+\underbrace{\psi_{3,k}}_{\in
\fffT{k\delta-\epsilon}00}
 \;. \eeal{begindn3}
 We note that
\bel{phikre}x\varphi\in x^{\delta} \Axd+ x^{\delta }y\Axyd
+\isAxf{(k+1)\delta-\epsilon}
+\fffT{(k+1)\delta-\epsilon}00\;.  \ee It follows that
$\tf=(\psi_1,x\psi_2,x\varphi)$ belongs to the space \bel{gphgk}\tf
\in \Big( \Axyd +\isAxf{k\delta-\epsilon}
+\fffT{k\delta-\epsilon}00\Big)\;.  \ee Note that our result is purely
local, so without loss of generality we can take the domain $\mcO$ of
the $v$--cordinates to be a closed coordinate ball.
Lemma~\ref{dlemmaGpol} below, with $\lambda=k\delta-\epsilon$, applied
to $g=\tf $ gives \bel{Hconcl4} x^{-p\delta}H(\cdot,x^{q\delta}
\tf)\in x^{(mq-p)\delta}\Axyd+{ \isAxf{(k+mq-p)\delta -\epsilon}}+
\fffT{(k+mq-p)\delta -\epsilon}00\;.\ee We can repeat now the
integration in
\eq{begind48}, recovering 
\eq{begindn3} with
$k=k+1$. This, together with \eq{phikre} shows that \eq{gphgk}
holds with $k$ replaced by $k+1$, hence \eq{Hconcl4} holds with
$k$ replaced by $k+1$. An integration as in \eq{begind50} gives
now \eq{begindn1} with $k$ replaced by $k+1$, and the induction
step is completed.

By an argument analogous to the one in the linear case, we conclude
that
$$\varphi \in x^{\delta-1}\Axd + x^{\delta-1}y\Axyd\;,\quad \psi \in
x^\delta \Axyd + \Ayd + \Axd\;.$$
An analysis as in
Remark~\ref{nwstst} finishes the proof.
\qed

We finish this section with the lemma referred to above:

\begin{Lemma}\label{dlemmaGpol} Assume that $\mcO$ is convex, compact, with interior
points. Let $q\in\N^*$, $\R\ni \lambda\ge 0$, and let $H(x^\mu,\argw)$ be
  $\Axyd$-polyhomogeneous with respect to $x^\mu$ with a zero of order
  $m$ in $\argw$.  If for all $\epsilon>0$ we have \bel{begindnpr}
  g\in \Big( \Axyd+{ \isAxf{\lambda}}+ \fffT{\lambda-\epsilon}00\Big)
  \;,\ee then it also holds, for all $\epsilon>0$, \bel{Hconcl}
  H(\cdot,x^{q\delta} g)\in x^{mq\delta}\Big( \Axyd+{
    \isAxf{\lambda}}+ \fffT{\lambda-\epsilon}00\Big)\;.\ee If
  $\lambda>0$ and \eq{begindnpr} holds with $\epsilon=0$, then
  \eq{Hconcl} also holds with $\epsilon=0$.
\end{Lemma}

\proof  This is a repetition of the proof
of~\cite[Lemma~4.8]{ChLengardnwe}, using an analogue of
\chcite{Lemma~A.5} in the $\mcT$--spaces. That last result is
proved using the usual Moser inequality together with  scaling, as
follows: For each $(\bar x, \bar v^A, \bar y) \in \Omega$ we
define a map $S_{(\bar x, \bar v^A, \bar y)}$ which maps a
standard set of parameters $\ztr\times \mcO$ to a compact subset
of $\Omega$, containing $(\bar x, \bar v^A, \bar y)$. By $\ztr$ we
denote a ``standard triangle'':
$$\ztr := \{ (s,t) : 0\le s\le t, \; 0\le t \le 1 \}\;.$$
The maps $S_{(\bar x, \bar v^A, \bar y)}$ are defined as
follows:
$$(s,t,v) \mapsto (x,y,v) = \cases{ (\bar x + \varsigma s, \bar y +
  \varsigma t,v), \; \varsigma = \bar x & for $\bar y < T$ \cr (\bar x +
  \varsigma (s-1), \bar y + \varsigma (t-1)), \; \varsigma = \bar x/2 & for
  $\bar y \ge T$}.$$
 This definition guarantees that the
$S$-image of $\ztr\times \mcO$ is contained in $\Omega$. We have
$$h\in\mcT^\alpha_k\quad \iff\quad  \forall i+j+|\gamma|\le k
\quad |\partial_v^\gamma \partial_x^i \partial_y^j h| \le C
x^{\alpha-i-j}\;.$$ There exists a smallest such $C$ which we
denote by $\|h\|_{\mcT^\alpha_k}$. For $g\in\mcT^\alpha_k$ we
estimate the norm of $H(\cdot,x^{q\delta} g)$ as follows:
$$\|H(\cdot,x^{q\delta} g)\|_{\mcT^\alpha_k(\overline\Omega)}\le
\sum_{i+j+|\gamma|\le k} \sup |x^{-\alpha+i+j}\partial_v^\gamma
\partial_x^i \partial_y^jH|\;,$$
$$
\partial_v^\gamma \partial_x^i \partial_y^jH =
\left(\frac{1}{\varsigma}\right)^{i+j}(\partial_v^\gamma \partial_s^i
\partial_t^j\tilde{H})\circ S^{-1}_{(\bar x, \bar v^A, \bar y)}\;,$$
where $\tilde H = H \circ S_{(\bar x, \bar v^A, \bar y)}$. We use
the interpolation inequality
$$\|u\|_{C_i(\ztr\times \mcO)} \le \|u\|^{1-i/k}_{L^\infty(\ztr\times \mcO)}
\|u\|^{i/k}_{C_k(\ztr\times \mcO)}\;,\qquad i\le k\;,$$ and
proceed as in \chcite{Proposition~A.2} to estimate the derivatives
of $\tilde H$. Let $z$ stand for $(s,t,v)\in\ztr\times\mcO$, then
for $|\sigma|\le k$ we have:
 \beqan |\partial^\sigma \tilde H|
&\le &C\Big|\sum_{|\gamma|+|\sigma_1|+\cdots +|\sigma_i|=|\sigma|}
  \varsigma^{q\delta(|\sigma_1|+\cdots +|\sigma_i|)}
\\
&& \qquad\qquad\times{\partial^{|\gamma|+i} \tilde H \over
    \partial z^\gamma \partial\argw^i}
  \partial^{\sigma_1}((1+s)^{q\delta} \tilde g)
  \cdots\partial^{\sigma_i}((1+s)^{q\delta} \tilde g)\Big|\\
&\leq& C\varsigma^{mq\delta} \sum_{|\sigma_1|+\cdots+|\sigma_i|\leq
  |\sigma|} | \partial^{\sigma_1} ((1+s)^{q\delta} \tilde g)
|\cdots|\partial^{\sigma_i}((1+s)^{q\delta} \tilde g)|\\
&\leq& C\varsigma^{mq\delta}\|(1+s)^{q\delta}\tilde
g\|_{L^\infty(\ztr\times\mcO)} \|(1+s)^{q\delta}\tilde
g\|_{C_k(\ztr\times\mcO)} \; , \eeqan where $\tilde g = g\circ
S_{(\bar x, \bar v^A, \bar y)}$. Let us note that $\bar x$ and
$\varsigma$ are of the same order, in fact $\frac{\bar
  x}{\varsigma}$ equals $1$ or $2$. Using this and the
  relationship
between $g$ and $\tilde g$ we conclude that
$$\|(1+s)^{q\delta}\tilde g\|_{C_k(\ztr\times\mcO)} \le C
\|g\|_{\mcT^\alpha_k(\overline\Omega)} \varsigma^\alpha\;.$$
Putting all the estimates together we arrive at
$$\|x^{-mq\delta}H(\cdot,x^{q\delta}
g)\|_{\mcT^\alpha_k(\overline\Omega)}\le
C\|g\|_{\mcT^\alpha_k(\overline\Omega)}\;,$$ which is what we need
to proceed with the proof of the Lemma. The point in all the
estimates is that $C$ depends only on $q\delta$, $k$ and
$\|g\|_{L^\infty}$ and \emph{does not} depend on $(\bar x, \bar v^A,
\bar y)$.

If $g\notin L^\infty(\overline\Omega)$ we can use $\hat g_\epsilon
:= x^\epsilon g \in L^\infty(\overline\Omega)$, hence the epsilons
in Lemma \ref{dlemmaGpol}. \qed

\section*{Acknowledgements}
We would like to acknowledge support by a Polish Research
Committee grant 2 P03B 073 24, by the Schr\"odginer Institute,
Vienna, and by the PAN-CNRS exchange programme. Sz.~\L.~was also
supported in part by a scholarship from the Foundation for Polish
Science.

\appendix

\section{Function spaces, auxilliary results}\label{SA}
\subsection{$\mcC$--spaces}

\label{SpwmCc}

For $k\in \N$ we denote by $C_k(\Omega)$ the set of all functions
which are $k$ times continuously differentiable on $\Omega$. We
denote by $C_k(\overline \Omega)$ the set of
$C_k(\Omega)$--functions which can be extended by continuity to
$C_k$ functions defined in an open neighborhood of $\Omega$.

Let $F$ be a space of functions, we shall say that $f\in x^\alpha
y^\beta F$ if $x^{-\alpha}y^{-\beta} f \in F$.

 We need to
introduce various families of function spaces with controlled
singular behavior at $\{x=0\}$, or $\{y=0\}$,  or $\{x=y=0\}$. The
domains we will consider will always be subsets of the set $0\le
x\le y\le y_0$ for some $y_0<\infty$. For $k\in\N$ we define \beaa
 \fCx{\alpha}{k}(\Omega) &=&
\{f: \forall \ i,j\in \N\;,\  \beta\in\N^{r}\;,\ i+j+|\beta|\le
k\;, \ \sup_\Omega|
x^{-\alpha}\partial_v^\beta[\partial_y]^i[x\partial_x]^j
 f | < \infty\}
 \;,\\\fCy{\sigma}{k}(\Omega) &=&
\{f: \forall \ i,j\in \N\;,\  \beta\in\N^{r}\;,\ i+j+|\beta|\le
k\;, \ \sup_\Omega|
y^{-\sigma}\partial_v^\beta[y\partial_y]^i[\partial_x]^j
 f | < \infty\}
 \;,\\
\cDak(\Omega) &=& \{f: \forall \ i,j\in \N\;,\
\beta\in\N^{r}\;,\ i+j+|\beta|\le k\;, \ \sup_\Omega|
x^{-\alpha}\partial_v^\beta[y\partial_y]^i[x\partial_x]^j
 f | < \infty\}
 \;,\\ \cDask(\Omega) &=&
\{f: \forall \ i,j\in \N\;,\  \beta\in\N^{r}\;,\ i+j+|\beta|\le
k\;, \ \sup_\Omega|
x^{-\alpha}y^{-\sigma}\partial_v^\beta[y\partial_y]^i[x\partial_x]^j
 f | < \infty\}
\;. \eean The spaces $ \fCx{\alpha}{k}$ here correspond to the
spaces $ \mcC^{\alpha}_{k}$ of \chl.

 We shall
write
$$  \cDai=\cap_{k\in\N}\cDak\;,$$
similarly for $C_\infty(\Omega)$, $\fCxi {\alpha}$, etc.

We  note the following:
\begin{Proposition}
\label{Ptrivemb} 
$\forall \delta \ge 0$ we have $x^\alpha \fCy
\beta k \subset \fCxy {\alpha-\delta,\beta+\delta} k\;.$
%
\end{Proposition}

\proof
There exist constants $C_{\ell,\alpha}$ and $C$ such that, for
 $0\le x\le y$ and for $f\in \fCy \beta k$ we have
\beaa(x\partial_x)^i (y\partial_y)^j
\partial_v^\beta (x^\alpha f) &= & x^\alpha \sum_{\ell=0}^i
C_{\ell,\alpha} (x\partial_x)^\ell (y\partial_y)^j
\partial_v^\beta f
\\ & \le & C x^\alpha y^\beta = C x^{\alpha-\delta} x^\delta
y^\beta \le C x^{\alpha-\delta} y^{\beta+ \delta}\;.\eeaa
 \qed

\subsection{Polyhomogeneous functions ($\mcA$--spaces)}\label{Sphgf}
 A
function $f\in C_\infty
\left(\Omega\right)
$ will be said to be \emph{polyhomogeneous at} $\{x=y=0\}$ if
there exist  integers $N_{i}$,   real numbers $n_{i}, \hat n_{i}
$,  and functions $f_{ij\ell}\in
C_\infty\left(\overline{\Omega}\right)$ with the property that
 \beal{cphgn} &&\forall m\in \N\quad \exists \
N(m)\
\ \textrm{such that}\ \nn \\
&&\phantom{cocoxxxxxxxx}
f-\sum_{i=0}^{N(m)}\sum_{j,\ell=0}^{N_{i}} f_{ij\ell}\ y^{\hat
n_{i}} x^{n_{i}} \ln^j y \ln^\ell x\in C_m(\overline\Omega)\;.\eea
We then write $f\in\Axy$. To avoid repetitions of terms with
identical powers in \eq{cphgn} it is convenient to impose
$(n_i,\hat n_i)\ne (n_j,\hat n_j)$ for $i\ne j$, and we will
always assume that this condition is satisfied.

We write \bel{phgdrestr}\mbox{$f\in\Axyd$ when $\delta\in 1/\N^*$
and } \ \{n_{i}\}\subset\delta \N\;,\ \{\hat n_{i}\}\subset\delta
\Z\;,\quad \hat n_i\ge -n_i\;.\ee The last inequality (recall that
$0\le x \le y$) guarantees that functions in $\Axyd$ are estimated
by $C(1+|\ln x|^N)$ for some $N\in \N$, and that they are bounded up
to a finite number of logarithmic terms. The need for  negative
powers of $y$ arises from the requirement of invariance of $\Axyd$
under $x\partial _y$. As an example, consider the function $\ln y\in
\Axyd$ for any $\delta$, then $x\partial_y \ln y =x/y$. This
exhibits the necessity of negative powers of $y$, appearing however
in a way consistent with the  inequalities in \eq{phgdrestr}.

Similarly we shall write $f\in \Ax$ if \eq{cphgn} holds with $\hat
n_i=0$ for all $i$, and with no non-trivial
powers of  $\ln y$, 
with the obvious definitions for $\Axd$, $\Ay$, etc.

The following observation will be used repeatedly:
\begin{Proposition}
\label{Pphgbas} 1. We have the inclusion
$$\Axyd\cap L^\infty \subset\fCxy{0}{\infty}\;.$$ It follows that
for any $\epsilon>0$ we have
$\Axyd\subset\fCxy{0-\epsilon}{\infty}\;.$

2. Similarly
$$\Axd\cap L^\infty
\subset\fCx{0}{\infty}\;,$$ and for  any $\epsilon>0$ we have
$\Axd\subset\fCx{0-\epsilon}{\infty}\;.$
\end{Proposition}
\proof 1. The last statement is obvious given the previous ones.
It remains to show that for all $k\in \N$ we have $\Axyd\cap
L^\infty \subset\fCxy{0}{k}\;.$ But clearly each term in the sum
in \eq{cphgn} with $m=k$, as well as $f$ minus the whole sum, is
in $\fCxy{0}{k}$, whence the result. The proof of point 2. is
identical. \qed

We shall need the following characterisation of the space of
polyhomogeneous functions:

\begin{Proposition} \label{Pequiv}
$f\in\Ax$ if and only if for every $m\in \N$ there exist $N(m),
N_i(m), n_i(m)$ and functions $f_{ij}\in C_m({\overline \Omega})$
such that  \beal{cphgold} f-\sum_{i=0}^{N(m)}\sum_{j=0}^{N_{i}(m)}
f_{ij}\ x^{n_{i}(m)} \ln^j x\in C_m(\overline\Omega)\;,\eea with a
similar property for $\Axd$, $\Axy$, etc.
\end{Proposition}

\proof The direct implication is obvious, therefore we only need
to prove the reverse one. We start by noting that in \eq{cphgold}
we can choose the $f_{ij}$'s to be independent of $x$. Indeed, if
some $f_{ij}$ depends on $x$ one can Taylor-expand it with respect
to $x=0$: $$ f_{ij}(x,v^A, y) = \sum_{\ell=0}^n
\frac{1}{\ell!}\frac{\partial^\ell f_{ij}(0,v^A,y)}{\partial
x^\ell} x^\ell +r\;.$$ The coefficients of the expansion are in
$C_{m-n}(\overline\Omega)$. Therefore
$$f_{ij} x^{n_i(m)} \ln^j x = \sum_{\ell=0}^n
\tf_{ij\ell}(v^A,y) x^{\ell+{n_i(m)}} \ln^j x + x^{n_i(m)} \ln^j x
\,r\;.$$ The last term is in $C_k(\overline\Omega)$ with $k$
arbitrary large, provided one started with $m$ large enough.
Similarly, one can arrange things so that the
$C_m(\overline\Omega)$ term in \eq{cphgold} is $o(x^m)$. Hence for
every $m$ there exist $N(m)$, $N_i(m)$, $n_i(m)$ and
$f_{ij}(v^A,y)$ such that \beal{cphgold2}
f=\sum_{i=0}^{N(m)}\sum_{j=0}^{N_{i}(m)} f_{ij}\ x^{n_{i}(m)}
\ln^j x+ r_m\;,\eea with $r_m \in C_m(\overline\Omega)$ and
$r_m=o(x^m)$. If necessary we rearrange the $n_i$'s increasingly,
and we also assume that $f_{ij}$'s do not vanish identically.

The next step in the proof is to show that the $f_{ij}$'s, defined
as in \eq{cphgold2}, do not depend upon $m$, i.e., if one writes
\eq{cphgold2} with some $m^\prime$ then the $f_{ij}$'s
corresponding to $n_i < \min \{m, m^\prime\}$ will be equal. This
is proved by comparing the expansions order by order, starting
with the smallest $n_i$ and largest $j$. For example:
$$ f= \sum_{j=0}^{N_0(m)} f_{0j}x^{n_0(m)}\ln^j x
+ o(x^{n_0(m)+\epsilon}) = \sum_{j=0}^{N_0(m^\prime)}
f^\prime_{0j}x^{n_0(m^\prime)}\ln^j x +
o(x^{n_0(m^\prime)+\epsilon})\;.$$ Dividing this equation by
$x^{n_0(m)} \ln^{N_0(m)} x$ and taking the limit $x\to0$ one gets
$n_0(m) \le n_0(m^\prime)$, similarly dividing by
$x^{n_0(m^\prime)} \ln^{N_0(m^\prime)} x$ one gets $n_0(m) \ge
n_0(m^\prime)$, hence $n_0(m) = n_0(m^\prime)$, $N_0(m) =
N_0(m^\prime)$ and $f_{0N_0} = f^\prime_{0N_0}$. Once the
uniqueness of $N, N_i, n_i$ and $f_{ij}$ is established, the
smoothness of $f_{ij}$ follows. This proves that $f\in\Ax$. The
result in $\Axd$ follows as the property $n_i\in\delta\N$  is
clearly preserved by the above procedure.  The result for $\Axy$
is proved similarly, using Taylor expansions both in $x$ and $y$.
\qed

 The space denoted by $\Axd$
here coincides with the space $\mcA^\delta_\infty$
of~\cite{ChLengardnwe}:

\begin{Proposition}$f\in \Ax$ if and only if for all $ \alpha\in \R$ there exists
  $ \hat N(\alpha)$ and functions $f_{i\ell}\in
  C_\infty(\overline{\Omega})$ such that
\bel{cphgold3}
 f-\sum_{i=0}^{\hat N(\alpha)}\sum_{\ell=0}^{N_{i}}
f_{i\ell}\
 x^{n_{i}}  \ln^\ell x\in \fCx{\alpha}{\infty}\;.\ee
This implies $\Axd=\mcA^\delta_\infty$.
\end{Proposition}

\proof Let $f$ satisfy the condition in the right member of the
equivalence above, then the property $f\in \Ax$ follows from the
inclusion $\fCx{\alpha}{\infty}\subset C_k(\overline \Omega)$ for
$k<\alpha$. Reciprocally, suppose that there exist  integers
$N_{i}$,   real numbers $n_{i}$, and  functions $f_{i\ell}\in
C_\infty\left(\overline{\Omega}\right)$ such that
 for all $m\in \N$ there exists $ N(m)$ such that
\bel{cphgn2} f-\sum_{i=0}^{N(m)}\sum_{\ell=0}^{N_{i}} f_{i\ell}\
x^{n_{i}} \ln^\ell x=:r_m\in C_m(\overline\Omega)\;.\ee Replacing
$f_{il}$ by $x^k f'_{il}$ if necessary, without loss of generality
we may assume that $f_{il}|_{x=0}\not \equiv 0$ for all
$i,\ell\in\N$. From what has been said it should be clear that the
set $\{n_i\le 0\}$ is finite. It suffices to prove the result for
$$f-\sum_{n_i\le 0}\sum_{\ell=0}^{N_{i}} f_{i\ell}\ x^{n_{i}} \ln^\ell
x\;,$$ and thus without loss of generality we may assume that
$f\in L^\infty$.  By Proposition~\ref{Pphgbas} the left-hand side
of the defining equality \eq{cphgn2} is then in $\fCx{0}\infty$,
thus $r_m\in \fCx{0}\infty$.

It follows from the definition that $f|_{x=0}\in C_\infty$, which
clearly implies that $r_m|_{x=0}\in C_\infty$. One similarly shows
that $(\partial_x^ir_m)|_{x=0}\in C_\infty$  for $0\le i\le m$.
Let $r'_m$ be defined as $r_m$ minus its Taylor series in $x$ of
order $m-1$, 
then we still have $r'_m\in \fCx{0}\infty$, and  $r'_m=O(x^m)$,
thus
$$r'_m\in \fCx{0}\infty \cap \fCx m 0\;.$$
Redefining the $f_{i\ell}$'s, \eq{cphgn2} still holds with $r_m$
replaced by $r'_m$.

 Let $m'> m$, we then have
\bel{rpmp}r'_{m}-\sum_{i=N(m)+1}^{N(m')}\sum_{\ell=0}^{N_{i}}
f_{i\ell}\ x^{n_{i}} \ln^\ell x=r'_{m'}\in \fCx{0}\infty \cap \fCx
{m'} 0\;,\ee with each term in the sum being $O(x^m)$ (otherwise
$r'_{m'}$ wouldn't be $O(x^{m'})$). Recall  the usual interpolation
inequality~\cite{Hormander}, for $0 < k<\ell$,
$$\|f\|_{C_k} \le C(k,\ell)\|f\|_{C_0}^{1-{k \over
\ell}}\|f\|_{C_\ell}^{k \over \ell}\;;$$ its weighted equivalent
reads (compare the proof of~\cite[Lemma~A.4]{ChLengardnwe})
$$\|r'_{m'}\|_{\fCx {(1-{k\over \ell}) m'} k} \le C'(k,\ell)\|r'_{m'}\|_{\fCx {m'} 0}^{1-{k \over
\ell}}\|r'_{m'}\|_{\fCx 0 \ell}^{k \over \ell}\;.$$ Given $k\in
\N$ 
we choose $\ell=2k$, $m'=2m$, leading to
$$r'_{m}-\sum_{i=N(m)+1}^{N(m')}\sum_{\ell=0}^{N_{i}} f_{i\ell}\ x^{n_{i}}
\ln^\ell x\in \fCx{m }k\quad \Longrightarrow \quad r'_{m}\in \fCx{m
}k\;.$$ Since $k$ is arbitrary, we find that
$$r'_m\in \fCx{m }\infty\;,$$
and our claim follows. \qed

Let $F$ be a space of functions on $\Omega$ such that $F\subset
C_\infty(\Omega)$. We shall say that $f\in \fAx{F}$ if for any
$k\in\N$ there exists $N(k)$ and functions $\varphi_{ij}\in F$
such that
\bel{fAF} f-\sum_{i,j=0}^{N(k)} \varphi_{ij} x^{\delta
i}\ln^j x \in C_k(\overline \Omega)\;. \ee 
The spaces $\fAxy{F}$ are defined in a similar way; for example we
have the identity $$\cAxyd=\fAxy{C_\infty(\overline \Omega)}\;.$$
In this notation it holds that
$$\cAxyd=\fAx{\fAyd}\;.$$

We will need the following characterisation of functions which are
polyhomogeneous up to lower order terms. To avoid annoying special
cases involving logarithms we assume $\sigma\not \in \N$, though the
proof gives also a corresponding statement in this case:
\begin{Proposition}
\label{Pequiv4}  Suppose that $\sigma\not\in\N$, let
$$f|_\hyp \in x^\beta\Axy\;,\quad f\in x^\beta
\Axyd+y^{\beta}\Axyd + \cDsk\;,$$ and assume that for all $i,j$
satisfying $i+j\le k+1$ there exists
$$g_{i,j} \in x^\beta \Axyd+y^{\beta}\Axyd $$
such that for every multi-index $\gamma$ for which $i+j+|\gamma|=
k+1$ we have
$$|(x\partial_x)^i (y\partial_y)^j\partial^\gamma_v (f-g_{i,j})|\le Cx^\sigma\;.$$
Then \bel{kgh} f\in x^\beta \Axyd+y^{\beta}\Axyd +\cDskp\;.\ee
\end{Proposition}

\proof For $\sigma<\beta $ and,   simultaneously, $\sigma<0$ there
is nothing to prove, as then the first two spaces are included in
the third one, and the claimed decomposition of $f$ is
uninteresting, and can be done in many different ways. We therefore
assume that at least one of those inequalities is violated, and we
proceed by induction on $k$. Suppose we know that the property
\bel{kpro}
\partial_Af\;,\ y\partial_yf\;,\ x\partial_x f\in x^\beta
\Axyd+y^{\beta}\Axyd + \cDsk\ee   implies \eq{kgh}. Then
Proposition~\ref{Pequiv4} is also established for $k=0$. Next,
applying the already established case $k=0$ of the proposition to
$(x\partial_x)^{\ell-1}f$ we find that
$$(x\partial_x)^{\ell-1}f\in x^\beta \Axyd+y^{\beta}\Axyd +\cDsone\;.
$$
Similarly, applying the case $k=0$ of the current proposition to
$y\partial_y(x\partial_x)^{\ell-2}f$ and
$\partial_A(x\partial_x)^{\ell-2}f$ we find that \beaa &
y\partial_y(x\partial_x)^{\ell-2}f\in x^\beta \Axyd+y^{\beta}\Axyd
+\cDsone\; ,
 & \\ &
\partial_A (x\partial_x)^{\ell-2}f\in x^\beta
\Axyd+y^{\beta}\Axyd +\cDsone\; .
 &
 \eeaa
The implication \eq{kpro}$\Rightarrow$\eq{kgh} gives then
$$(x\partial_x)^{\ell-2}f\in x^\beta \Axyd+y^{\beta}\Axyd +\cDstwo\;.
$$
Taking $\ell=2$ establishes  Proposition~\ref{Pequiv4}  with $k=1$.
Continuing in this way, Proposition~\ref{Pequiv4} follows for all
$k\in\N$.

It remains to establish the implication
\eq{kpro}$\Rightarrow$\eq{kgh}. Writing
$$\partial_x f= \underbrace{f_{x,k,\phg}}_{\in x^{\beta-1}
\Axyd+x^{-1}y^{\beta}\Axyd} + \underbrace{f_{x,k,\sigma}}_{\in
\cDsmk}\;,$$ we have\bel{psiint2} f(x,v^A,y) = \underbrace{
f(y,v^A,y) + \int_y^{x} f_{x,k,\phg}(s,v^A,
y)\,ds}_{=:f_{k+1,\phg}\in x^\beta \Axyd+y^{\beta}\Axyd}+
\underbrace{\int_y^{x} f_{x,k,\sigma}(s,v^A,
 y)\,ds}_{=:f_{k+1,\sigma}}\;.\ee
We want to show that $f_{k+1,\sigma}\in\cDskp$; equivalently:
$$\forall \gamma\;, \ i+j+|\gamma|\le k+1 \qquad |\partial^\gamma_v (\partial_x)^i
 (\partial_y)^j f_{k+1,\sigma}| \le C x^{\sigma-i}y^{-j}\;.$$
The inequality is clear  if $i\ge 1$ or if $|\gamma|+j<k+1$, by
differentiating the integral defining $f_{k+1,\sigma}$. Suppose thus
that $i=0$ and $|\gamma|+j= k+1$. Recall that, by hypothesis, we
have
$$\partial_v^\gamma (y\partial_y)^{j}f= \underbrace{\partial_v^\gamma f_{j,y,k+1,\phg}}_{\in x^\beta
\Axyd+y^{\beta}\Axyd}+ O(x^\sigma)\;,$$ and comparing with
\eq{psiint2}, \bel{noico}\partial_v^\gamma
(y\partial_y)^{j}f_{k+1,\sigma}= \underbrace{\partial_v^\gamma
 f_{j,y,k+1,\phg}-\partial_v^\gamma (y\partial_y)^{j}f_{k+1,\phg}}_{=:\underbrace{g_{k+1,1}}_{\in x^\beta
 \Axyd}+\underbrace{g_{k+1,2}}_{\in y^{\beta}\Axyd}}+ O(x^\sigma)\;.\ee
Recall, now, that we only have to consider the cases $\sigma\ge
\beta $, or $\sigma\ge 0$, or both. For $y>0$ consider the collected
polyhomogeneous terms in \eq{noico}, they can be written as a finite
sum
 \beaa g_{k+1,1}(x,y,v^A)&=&\sum_{\beta+n\delta \le \sigma}\sum_{m}^{N_n}f_{mn}(y,v^A)x^{\beta+n\delta}\ln^mx
+ r_{k+1,1}(x,y,v^A)\;,
 \\
 g_{k+1,2}(x,y,v^A)&=&\sum_{n\delta \le \sigma}\sum_{m}^{N'_n}f'_{mn}(y,v^A)x^{n\delta}\ln^mx
+ r_{k+1,2}(x,y,v^A)\;, \eeaa
 with functions $f_{nm}$ and $f'_{nm}$ which are smooth in both variables (as long as $y>0$), and
 with the remainders being
$O(x^\sigma)$. If $\beta\in -\N$ we absorb the second sum in the
first one; similarly if $\beta\in\N$ we absorb the first sum in the
second one. This shows that without loss of generality we can always
assume that all powers of $x$ appearing in the sums above are now
pairwise distinct. Suppose that $f_{mn}\not \equiv 0$ or
$f'_{mn}\not \equiv 0$ for some couple $mn$, then $f_{k+1,\sigma}$
wouldn't be $O(x^\sigma)$ for all $(y,v^A)$'s by integration of
\eq{noico} in $y$ or in $v$. It follows that
$$g_{k+1,1}+g_{k+1,2}=r_{k+1,1}+r_{k+1,2}=O(x^\sigma)\;,$$ which establishes our
claim. \qed

\subsection{$\mcF$-- and  $\mcT$--spaces} \label{STsps}

For $\alpha\in \R$ and $k\in\N$ we set
\begin{eqnarray}\nonumber \lefteqn{\mcFa_{\{\zlxy \},k} =
\Bigg\{f\;|\;\forall\; 0\leq i+j+|\gamma| \leq k \quad \exists N:}
&&\\ && |\partial^i_x\partial^j_y\partial^\gamma_v f| \leq
\cases{Cy^{\alpha-i-j} (1+|\ln y|)^N& \mbox{if } $\alpha-i-j\geq0$
\cr Cx^{\alpha-i-j} (1+|\ln x|)^N& \mbox{if }
$\alpha-i-j<0$}\,\Bigg\}\;. \label{Tadef1}
\end{eqnarray}
We will also need a version of the $\mcF$-spaces where the
functions involved are ``almost independent of $x$ when $\alpha$
is large", in the following sense:
\begin{eqnarray}\nonumber
\lefteqn{\mczTa_{\{\zlxy \},k} = \Bigg\{f\;|\;\forall\; 0\leq
i+j+|\gamma| \leq k \quad \exists N:} &&\\ &&
|\partial^i_x\partial^j_y\partial^\gamma_v f| \leq
\cases{Cy^{\alpha-j} (1+|\ln y|)^N& \mbox{if } $\alpha-j\geq0$,
$i=0$ \cr Cx^{\alpha-i-j} (1+|\ln x|)^N& \mbox{otherwise}
}\,\Bigg\}\;. \label{Tadef}
\end{eqnarray}

Let $\alpha,\beta\in \R$, $k\in \N$. To be able to estimate in terms
of powers of $|\ln x|$ rather than $1+|\ln x|$ it is convenient to
assume $0<y_0<1$. We say that $f\in\fffT{\alpha}\beta k$ if for all
$ i,j,\gamma$ there exist constants $C>0$ and $N\in\N$ such that,
for $0<x\le y \le y_0$ we have \bea
 |\partial_x^i\partial_y^j
\partial_v^\gamma f |&\le& C \Big(x^{\alpha+\beta-i-j} +
x^{\alpha-i} y^{\beta-j} + x^{\alpha
+\beta-i-k}y^{k-j}\Big)|\ln^Nx|\;.
\eeal{ffTdef} We write $f\in\ffT{\alpha}\beta $ for
$f\in\fffT{\alpha}\beta 0$, and we note that for $k=0$, or for
$\beta=k$, the last term in \eq{ffTdef}  is not needed, e.g.:
\bean
f\in \ffT{\alpha}\beta  \quad \Longleftrightarrow \quad
|\partial_x^i\partial_y^j
\partial_v^\gamma f |\le C \Big(x^{\alpha+\beta-i-j} +
x^{\alpha-i} y^{\beta-j}\Big)|\ln^Nx|\;.\nn
\\ && \eeal{ffTdefb}
Finally, for $\beta\le 0$ the last term in \eq{ffTdefb} can be
dropped altogether.

Strictly speaking, the only space out of the $\fffT{\alpha}\beta
k$'s which is absolutely necessary in our proofs is the one with
$k=\beta=0$. However, we have decided to include a short
discussion of the other ones as well, as those spaces appear
naturally in the problem at hand.

Let $\{F_i\}_{i\in \N}$ be any countable family of function
spaces, we shall write
$$\sCysF = \{f: \exists N\in\N\;, f_n \in F_n\;,0
\le n\le N\;, f= \sum_{n=0}^N f_n\}\;.$$ The dot over the symbol
$\oplus$ is meant to emphasise the fact that only finite linear
combinations are considered.

For further use we note the following elementary properties:

\begin{Proposition}
\label{PffT}
\begin{enumerate}\item\label{PffT-1} If  $f\in\fffT{\alpha}\beta k$
then $\partial_xf \in\fffT{\alpha-1}\beta k$ and $\partial_yf
\in\fffT{\alpha}{\beta-1} {\max(k-1,0)}$
\item\label{PffT-4}   For $\alpha'\ge \alpha$  and $\beta'\ge \beta$ we have $\fffT{\alpha'}{\beta'} k\subset \fffT{\alpha}{\beta} k$.
\item\label{PffT-3} For $\sigma\ge 0$ we have $\fffT{\alpha+\sigma}\beta k\subset \fffT{\alpha}{\beta+\sigma} k$.
\item\label{PffT-2} If $\N\ni\ell< \alpha$ and $\ell\le k \le \beta$
we have $\fffT{\alpha}\beta k\subset \fCy{\beta} \ell$.
\item\label{PffT3n} $ \fCxi \alpha \subset  \fCxy{\alpha}\infty \subset \fffT{\alpha}{0} k $, and
$\fCy \beta \infty\subset \fffT{0}\beta k$ for all $k$.
\item\label{PffT1} If  $f\in\fffT{\alpha}\beta k$ and $g\in C^\infty$, then $fg\in\fffT{\alpha}\beta
k$.
\item\label{PffT2} If $g\in \fCx \alpha \infty$ and $h\in \fCy
\beta \infty$ then $gh\in \fffT{\alpha}\beta k$ for all $k$.
\item \label{PffT4} We have  $x^\sigma\fffT{\alpha}\beta k= \fffT{\alpha+\sigma}{\beta} {k}$
 for all $k\in \N $ and  $\sigma \in \R$.
\item \label{PffT5} For $f\in\fffT{\alpha}\beta k$ and $\ell\in \N$ we have
$x^\ell f \in \fffT{\alpha}{\beta+\ell} {k+\ell}$ for all
$k$.
\end{enumerate}
\end{Proposition}

\proof 
Points \ref{PffT1}--\ref{PffT4} follow immediately from the
formula
 \bel{deri} \partial^i_x \partial^j_y (gh) = \sum_{\scriptsize \begin{array}{c}
   r+s=i \\
   m+n=j \\
 \end{array}}C(r,s,m,n)(\partial^
 r_x \partial^m_y g)(\partial^s_x \partial^n_y h)\;.
\ee The remaining claims are direct consequences of the
definition. \qed

\subsection{Extensions of a class of functions}\label{Sext}

Let $0\le \varphi\in C^\infty(\R)$, $\supp \varphi \subset
[-1/2,1/2]$, $\int_\R \varphi(x)dx=1$. For $0<x\le y\le y_0$ we
set
 \minilab{Emdef}
  \begin{equs}\label{Emdef1} E[f](x,y,v) &:=
\int_0^\infty \frac{\varphi(\frac{w-y}x)}x f(w,v)dw \\
&= \int_{y/2}^{3y/2}
\frac{\varphi(\frac{w-y}x)}x f(w,v)dw \label{Emdef2.0}\\
&= \int_{-\infty}^\infty
\frac{\varphi(\frac{w-y}x)}x f(w,v)dw \label{Emdef2}\\
&=\int_{-\infty}^{\infty} \varphi(z)f(y+xz,v)dz \label{Emdef2.1}
\\
&=\int_{-1/2}^{1/2} \varphi(z)f(y+xz,v)dz \label{Emdef3}\end{equs}
(
there is no need to know the values of $f$ for negative $w$ when
using \eq{Emdef2} as $\varphi= 0$ there; a similar comment applies
to \eq{Emdef2.1}).

The results here are an adaptation to the problem at hand of
\cite[Section~3.3]{AndChDiss}. In the lemma that follows one can
think of $\localmu $ as belonging to $[0,1)$, but this restriction
is not necessary for the result:

\begin{Lemma}
\label{LExt1} For $k\in\N$ and $\localmu \in \R$ suppose that
\begin{equa}[LExt1.1]
     |\partial^\gamma_v
\partial^\ell_y f|&\le Cy^{k+\localmu -\ell}(1+|\ln y|)^N\ \mbox{ for } \ 0\le \ell\le
k\;,
\end{equa}
then \bel{ext0} E[f]\in y^\localmu  \fT k \infty\;.\ee
 If moreover
there exists $\lambda>0$ such that \bel{LExt1.2}
   |\partial^\gamma_v
\partial^k_y f(y,v)-\partial^\gamma_v
\partial^k_y f(y',v)|\le Cy^{\localmu -\lambda}(1+|\ln y|)^N |y-y'|^\lambda \ \mbox{ for } \ |y'-y|\le
y/2\;,\ \ee
  then we also have
 \bel{ext0.3} E[f](x,y,v)
 \in y^{\localmu -\lambda} \fT {k+\lambda} \infty \;.\ee
\end{Lemma}

\proof
 The estimate \eq{LExt1.1} together with \eq{Emdef3} gives
\bean
0\le i+j\le k \qquad |\partial_x^i \partial_y^j \partial _v^\gamma
E[f](x,y,v) | &=&\Big| \int_{-1/2}^{1/2}
\varphi(z)\partial_y^{i+j}\partial^\gamma_vf(y+xz,v)z^idz\Big|
\\ &\le & Cy^{k+\localmu -i-j}(1+|\ln y|)^N\;.
\eeal{ext1} On the other hand, if $i+j>k$ we write $i=i_1+i_2$ and
$j=j_1+j_2$ with $i_2+j_2=k$, obtaining
\bean
 \partial_x^i \partial_y^j \partial _v^\gamma
E[f](x,y,v)  &=& \partial_x^{i_1}
\partial_y^{j_1}\int_{-1/2}^{1/2}
\underbrace{\varphi(z)z^{i_2}}_{=:\varphi_{i_2}(z)}\partial_y^{k}\partial^\gamma_vf(y+xz,v)dz
\\\nonumber&=& \partial_x^{i_1}
\partial_y^{j_1}\int_{0}^{\infty}
\frac{\varphi_{i_2}(\frac{w-y}x)}x\partial_y^{k}\partial^\gamma_vf(w,v)dw
\\&=&x^{-i_1-j_1} \int_{0}^{\infty}
\frac{\varphi_{i_1,i_2,j_1}(\frac{w-y}x)}x\partial_y^{k}\partial^\gamma_vf(w,v)
dw
 \;,
\eeal{ext1.4} where we have set
$$\varphi_{i_1,i_2,j_1}(\frac{w-y}x):= x^{1+i_1+j_1}\partial_x^{i_1}
\partial_y^{j_1}\big(\frac{\varphi_{i_2}(\frac{w-y}x)}x\Big)\;.$$
It follows that \bean
 |\partial_x^i \partial_y^j \partial _v^\gamma
E[f](x,y,v) | &\le&Cx^{-i_1-j_1}y^{\localmu } (1+|\ln
y|)^N\int_{0}^{\infty}
\Big|\frac{\varphi_{i_1,i_2,j_1}(\frac{w-y}x)}x\Big|dw
 \\\nonumber &=&Cx^{-i_1-j_1}y^{\localmu } (1+|\ln y|)^N\int_{-1/2}^{1/2}
|{\varphi_{i_1,i_2,j_1}(z)}|dz
\\
 &\le &
CC'x^{-i_1-j_1}y^{\localmu }(1+|\ln y|)^N \nonumber \\ &=&
CC'x^{k-i-j}y^{\localmu }(1+|\ln y|)^N\;, \eeal{ext2} which
establishes \eq{ext0}.

To prove \eq{ext0.3}, we start by noting that for $i_1+j_1>0$ we
have
\beaa 0 &=&
 \partial_x^{i_1}
\partial_y^{j_1}\int_{-\infty}^{\infty}
{\varphi_{i_2}(z)}dz
\\&=& \partial_x^{i_1}
\partial_y^{j_1}\int_{0}^{\infty}
\frac{\varphi_{i_2}(\frac{w-y}x)}x \;dw=\int_{0}^{\infty}
\partial_x^{i_1}
\partial_y^{j_1}\Big(\frac{\varphi_{i_2}(\frac{w-y}x)}x \Big)\; dw
 \;
\eean
This allows us to write (compare \eq{ext1.4})
\beaa
 \partial_x^i \partial_y^j \partial _v^\gamma
E[f](x,y,v)  &=&  \int_{0}^{\infty}
\partial_x^{i_1}
\partial_y^{j_1}\Big(\frac{\varphi_{i_2}(\frac{w-y}x)}x \Big)
\Big(\partial_y^{k}\partial^\gamma_vf(w,v)
-\partial_y^{k}\partial^\gamma_v f(y,v) \Big)\; dw
\\&=&x^{-i_1-j_1} \int_{0}^{\infty}
\frac{\varphi_{i_1,i_2,j_1}(\frac{w-y}x)}x\Big(\partial_y^{k}\partial^\gamma_vf(w,v)
-\partial_y^{k}\partial^\gamma_v f(y,v) \Big) dw\\&=&x^{-i_1-j_1}
\int_{\infty}^{\infty}
{\varphi_{i_1,i_2,j_1}(z)}\underbrace{\Big(\partial_y^{k}\partial^\gamma_vf(y+zx,v)
-\partial_y^{k}\partial^\gamma_v f(y,v) \Big)}_{\le C y^{\localmu
-\lambda}(1+|\ln  y|)^N x^\lambda} dz
\\ &\le & C'y^{\localmu -\lambda}(1+|\ln  y|)^N x^{k+\lambda-i-j}
 \;,
\eean
as desired.
 \qed

We continue with

\begin{Lemma}
\label{LExt2} Let $\localmu  \ge 0$ and for $0\le i \le m$ let
$f_i$ satisfy \eq{LExt1.1} with $k$ there replaced by $m-i$. There
exists $h\in y^\localmu  \fTi m$ such that \bel{derq}0\le i \le m
\qquad \partial_x^i h|_{x=0}= f_i\;.\ee If the $f_i$'s satisfy
\eq{LExt1.2} with $k=m-i$ then $h\in y^{\localmu -\lambda} \fTi
{m+\lambda}$.
\end{Lemma}

\proof We start by considering the following set of functions,
defined for $0\le i\le m$,
\bel{Emdefn} g_i(x,y,v)= \frac
{x^i}{i!} E[f_i](x,y,v)\;.\ee It follows from
Proposition~\ref{LExt1} that $g_i\in y^\localmu  x^i \fTi
{m-i}\subset y^\localmu  \fTi m$, or $g_i\in y^{\localmu -\lambda}
x^i \fTi {m+\lambda-i}\subset y^{\localmu -\lambda} \fTi
{m+\lambda}$ if \eq{ext0.3} holds. It also follows from \eq{ext1}
and from Lebesgue's differentiation theorem that for $0\le j \le
m$ the functions $\partial^j_xg_i$ extend by continuity to
continuous functions on $\{x=0\}$, and that $g_i|_{x=0}$ satisfies
\eq{LExt1.1} with $k$ there equal to $m-i$, and with the modulus
of H\"older continuity satisfying \eq{ext0.3} if this condition
was satisfied by the $f_i$'s.

Those considerations imply that the following inductive scheme is
well defined: we set $h_0=E[f_0]$ and
$$h_{i+1}(x,y,v)=h_i(x,y,v)+ \frac {x^i}{i!}
E[f_i-(\partial_x^ih_i)|_{x=0}]\;.$$ Then the function $h:=h_m$
satisfies \eq{derq}, and belongs to the spaces claimed.

\qed

\subsection{Integral operators on $\mcA$-- and $\mcC$--spaces}

For $0\le x \le y \le y_0<\infty$ set
\beal{Ione}I_1(f)(x,v^A,y)&=&\int_x^y f(s,v^A,y)ds\;, \\
I_2(f)(x,v^A,y)&=&\int_x^y f(x,v^A,s)ds\;.\eeal{Itwo} In our
arguments we will need to understand the action of $I_1$ and $I_2$
on various spaces defined above. We start with polyhomogeneous
functions:

\begin{Proposition}\label{calkA}
\begin{enumerate}
\item  Let $g\in x^\beta y^\gamma\Axyd$. Then
$$ I_1(g)\in\; y^{\beta+\gamma+1}\Ayd +
x^{\beta+1}y^\gamma\Axyd \;,$$
$$ I_2(g)\in\; x^{\beta+\gamma+1}\Axd +
x^{\beta}y^{\gamma+1}\Axyd \;.$$ It follows in particular that
$\Axy$ is stable under both integrations above.
\item Let $g\in x^\beta y^\gamma\Axd$. Then
$$I_1(g)\in\; y^{\beta+\gamma+1}\Ayd +
x^{\beta+1}y^\gamma\Axd \;,$$
$$ I_2(g)\in\; x^{\beta+\gamma+1}\Axd +
x^{\beta}y^{\gamma+1}\Axd \;.$$
\end{enumerate}
\end{Proposition}

\proof Let $f\in C_\infty(\overline{\Omega})$, $p\in\R, j\in\N$.
We start by showing that for every $m\in\N$ there exist an integer
$N$, sequences of numbers $k_{i}\in \N, \ell_{i}\in\N$, a sequence
of smooth functions $f_i$ and a function $r_m\in
C_m(\overline{\Omega})$ such that \bel{phgint}\int_x^{y}
f(s,v^A,y) s^p \ln^j s \,ds = \sum_{i=1}^N f_i \left( y^{p+k_i+1}
\ln^{\ell_i} y - x^{p+k_i+1} \ln^{\ell_i}x\right) + r_m\;.\ee
Several integrations by parts in the integral $\int v' u $ with
$v'=s^p$ and $u=\ln^j s$ yield this formula when $\partial_x f=0$.
The result for general $f$ is also obtained by integration by
parts by taking $u=f$ and  $v' = s^p \ln^j s$. Using the result
already proved with $f=1$ one obtains a $v$ with a power of $s$
higher by one. Repeating the integration by parts a finite number
of times one obtains a remainder term in $C_m(\overline{\Omega})$,
and one concludes by Proposition~\ref{Pequiv}. \qed

We continue with a study of the action of $I_1$ and $I_2$ on the
$\cDask$ spaces. Note that the action on the $\cDak$ spaces is
obtained as a special case from
$$\cDak=\cDazk\;.$$

\begin{Lemma}\label{Lintegr} Let
$\alpha,\sigma\in \R$, $k\in
\N\cup\{\infty\}$,
\begin{enumerate}
\item If  $f\in \cDask$, $\alpha < -1$, then $I_1(f)\in \cDapsk$.
\item If  $f\in \cDask$, $\alpha>-1$, then $I_1(f)\in
\cDapspkmy+\cDapskmxy$.
\end{enumerate}
\end{Lemma}

\proof 1. The case $\alpha<-1$ is obtained by straightforward
estimations.

2. For $\alpha>-1$ we write
\bel{Idec}
I_1(f)(x,v^A,y)=\underbrace{\int_0^y
f(s,v^A,y)ds}_{g_1}-\underbrace{\int_0^x f(s,v^A,y)ds}_{g_2}\;.
\ee We have $\partial^\gamma_v g_1=O(y^{\alpha+\sigma+1})$,
$\partial_x g_1=0$, $y^\ell\partial_y^\ell
\partial^\gamma_vg_1(x,v^A,y)= y^\ell\partial_y^{\ell-1}(
\partial^\gamma_vf(y,v^A,y))$, and all the estimates readily
follow.
%
%
\qed

\begin{Lemma}\label{Lintegr2} Let
$\alpha,\sigma\in \R$, $k\in
\N\cup\{\infty\}$,
\begin{enumerate}
\item If  $f\in \cDask$, $\sigma>-1$ then $I_2(f)\in
\cDaspkmxy$.
 \item If  $f\in \fCy{\sigma}{\infty}$,
then $I_2(f)\in \fT{\sigma+1}{\infty}$.
\end{enumerate}
\end{Lemma}

\proof 1.  One is tempted to argue  from point 2. of
Lemma~\ref{Lintegr2} by symmetry that $I_2(f)\in
\cDapspkmx+\cDaspkmxy$, but this is not clear, because one is not
allowed to integrate all the way to zero in $y$, as done in the
proof of 2. So we calculate directly:
\bea
\partial^\ell_x \partial_v^\beta \int_x^y f(x,v^A,s)ds & = &
\int_x^y \partial^\ell_x \partial_v^\beta f(x,v^A,s)\,ds
-\sum_{i=0}^{\ell-1} C_{\ell,i}
\partial^{\ell-1-i}_x\partial_y^i \partial_v^\beta f(x,v^A,s)\Big|_{s=x}\;,
\nonumber\\&&\label{phiderx2}\\
\partial^\ell_x \partial_y^i \partial_v^\beta
\int_x^y f(x,v^A,s)ds & = &
\partial^\ell_x \partial_y^{i-1} \partial_v^\beta
f(x,v^A,y)\;, \qquad i\ge 1\;.\qquad  \label{phiderx3}
 \eea
From
\eq{phiderx3} one immediately finds $|\partial^\ell_x \partial_y^i
\partial_v^\beta I_2(f)|\le C x^{\alpha-\ell}y^{\sigma+1-i}$ for
$i\ge 1$. Since $\sigma > -1$, the first term in \eq{phiderx2} is
estimated by $$Cx^{\alpha-\ell}|y^{\sigma+1}-x^{\sigma+1}|\le
Cx^{\alpha-\ell}y^{\sigma+1}\;,$$ as desired. Similarly, each term
in the sum is estimated by
$$Cx^{\alpha-\ell+i+1}y^{\sigma-i}\Big|_{y=x}=
Cx^{\alpha-\ell+1+\sigma}\le Cx^{\alpha-\ell}y^{\sigma+1}\;,$$ and
the result follows.

2.
 The estimates in the $\ell=0$ case are proved similarly as in point 1, see \eq{phiderx2}.
 For $1\le \ell \le k$ and $1\le i+\ell+|\gamma|\le k$ we have
$$|\partial_x^i\partial_y^\ell \partial_v^\gamma I_2(f)| = |\partial_x^i\partial_y^{\ell-1}
 \partial_v^\gamma f|\;,$$ and the desired estimate is
 straightforward.
 \qed

\subsection{Integral operators on $\mcT$-- and $\mcF$--spaces}\label{SIoT}
\begin{Proposition}
\label{PIoneT}
Let $\alpha>-1$, $\beta \ge k$. For any $\epsilon>0$ we have
$$I_1(\fffT{\alpha}\beta k)\subset  y^\epsilon\fzTi {\alpha+1-\epsilon+\beta} + \fffT
{\alpha+1-\epsilon} \beta k\;.$$
\end{Proposition}

\begin{remark}
\label{RIoT} We expect the result to remain valid with
$\epsilon=0$, but the proof below fails for this value of
$\epsilon$. In any case the current result is sufficient for our
purposes.
\end{remark}

\medskip

\proof  Write $\alpha+\beta=n-\sigma$, $n\in\N$, $\sigma\in
[0,1)$. Let $f\in\fffT{\alpha}\beta k$, set
$$g(y,v)=\int_0^y f(s,y,v)ds\;.$$ We want to show that the
function $ g$ satisfies the hypotheses needed to construct the
extensions of Section~\ref{Sext}. For $j$ such that
$\alpha+\beta-j>-1$
we have
\bel{dereq}\partial^j_y\partial^\gamma_v g (y,v)= \sum_{\ell+m=j-1}
C(\ell,m) \partial_x^\ell \partial_y^m \partial^\gamma_v f (y,y,v)
+ \int_0^y
\partial^j_y \partial^\gamma_v  f(x,y,v)dx\;. \ee From the
definition of $\fffT{\alpha}\beta k$ one thus finds, for $0\le j
\le n$,
$$|\partial^j_y\partial^\gamma_v g|\le Cy^{\alpha+\beta+1-j}|\ln y|^N
\;. $$ 
Next, we need to control the modulus of H\"older continuity of the
top order derivatives. In order to do that, consider any of the
terms appearing under the sum symbol in \eq{dereq} with $ j=n$.
For $y'\in B(y,y/2)$ we write
\bean\lefteqn{\!\!\!\!\!\!\!\!\!\!\!\!\!\!\!\!\!\!\!\!
\Big|\partial_x^\ell \partial_y^m \partial^\gamma_v\Big( f
(y,y,v)-f(y',y',v)\Big)\Big|  = \Big| \int_y^{y'}
\underbrace{\Big(\partial_x^{\ell+1}
\partial_y^m
\partial^\gamma_v f (s,s,v)+\partial_x^{\ell} \partial_y^{m+1}
\partial^\gamma_v f (s,s,v)\Big)}_{\le C y^{-\sigma}|\ln y|^N}
ds\Big|}&&
\\ &&  \phantom{xxxxx}\le C y^{-\sigma}|y'-y| |\ln y|^{N+1}\le C'
y^{1-\sigma-\lambda}|y'-y|^\lambda |\ln y|^{N+1}\eeal{dereq1} for
any $\lambda\in [0,1]$. Next, \bean \lefteqn{\int_0^y
\partial^{n}_y \partial^\gamma_v  f(x,y,v)dx-\int_0^{y'}
\partial^{n}_{y} \partial^\gamma_v  f(x,y',v)dx} &&\\ && =  \int_{y'}^y
\partial^{n}_y \partial^\gamma_v  f(x,y,v)dx+\int_0^{y'}\Big(
\partial^{n}_y \partial^\gamma_v  f(x,y,v)-\partial^{n}_y \partial^\gamma_v  f(x,y',v)\Big)dx \;.
\nonumber \\ && \eeal{dereq1a} The first term is estimated
similarly to \eq{dereq1},
 \bel{dereq2}
\Big|\int_{y'}^y
\partial^{n}_y \partial^\gamma_v  f(x,y,v)dx\Big|\le C
y^{1-\sigma-\lambda}|y'-y|^\lambda |\ln y|^{N+1}\; .\ee To control
the second we will need a weighted interpolation inequality,
obtained as follows: for $t\in [1/2,3/2]$ set
\bel{hdef} \chi(t)=
\partial^{n}_y \partial^\gamma_v  f(x,ty,v)\;.\ee Recall the
inequality~\cite[Theorem~A.5]{Hormander}
\bel{inter} \|\chi\|_{\lambda} \le C_\lambda \|\chi\|_{1}^\lambda
\|\chi\|_{0}^{1-\lambda}\;,\quad 0\le \lambda \le 1\;,\ee where
$\|\cdot\|_{\mu}$ denotes the usual H\"older $C^\mu([1/2,3/2])$
norm. Applying \eq{inter} to $\chi$ defined in \eq{hdef} gives
\bean \lefteqn{|\partial^{n}_y \partial^\gamma_v  f(x,y,v)-\partial^{n}_y \partial^\gamma_v
f(x,y',v)|}&&\\ &&\nonumber \le C\Big(\underbrace{\frac{y} x
x^{-\sigma}|\ln^N x|}_{\leftrightarrow \|\chi\|_1}\Big)^\lambda
\Big(\underbrace{x^{-\sigma}|\ln^N x|}_{\leftrightarrow
\|\chi\|_0}\Big)^{1-\lambda}y^{-\lambda} |y-y'|^\lambda \\ && = C
x^{-\lambda-\sigma} |\ln^N x| |y-y'|^\lambda\;.\eeal{fhin} It
follows that the last integral in \eq{dereq1a} converges for
$\lambda<1-\sigma$, leading to\bel{gholc} 0\le \lambda+\sigma < 1
\qquad \Big|\partial^n_y\partial^\gamma_v g (y,v) -
\partial^n_y\partial^\gamma_v g (y',v)\Big|\le C y^{1-\sigma-\lambda}|y'-y|^\lambda |\ln y|^N\; .\ee

Setting
$$g_0=
g\;,\quad g_i=0 \ \mbox{ for } \ 1\le i \le n \;,$$ by
Lemma~\ref{LExt2} there exists a function $h\in
y^{1-\sigma-\lambda}\fTi {n+\lambda}$ for any $0\le
\lambda+\sigma<1$  such that $\partial_x^i h|_{x=0}= g_i$. Define
$$\hat h = I_1(f)-h\;.$$  For $0\le i+j\le n$ we write
\beaa |\partial_x^i\partial_y^j
\partial_v^\gamma \hat h |&=& \Big|\partial_x^i\partial_y^j
\partial_v^\gamma \Big(\int_x^y f(s, \cdots) ds - h\Big) \Big|
\\
&=& \Big|\partial_x^i\partial_y^j
\partial_v^\gamma \Big(\int_0^y f(s, \cdots) ds - h- \int_0^x f(s, \cdots) ds\Big) \Big|
\\
&\le & \underbrace{\Big|\partial_x^i\partial_y^j
 \Big(\int_0^y \partial_v^\gamma f(s, \cdots) ds - \partial_v^\gamma h\Big)\Big|}_{I} +
 \underbrace{\Big|\partial_x^i
 \int_0^x \partial_y^j\partial_v^\gamma f(s, \cdots) ds \Big|}_{II}\;.
\\ && \eeaa
The estimate of $II$ is straightforward: \bel{IIest} |II|\le C
\Big(x^{\alpha+1+\beta-i-j} + x^{\alpha+1-i} y^{\beta-j} +
x^{\alpha+1 +\beta-i-k}y^{k-j}\Big)|\ln^Nx|\;.\ee
 To estimate $I$ we use Taylor's
formula,
 \beaa
\partial_x^i h(x,y,v) &=&  \sum_{\ell=0}^{n-i-j}
\frac{x^{\ell}}{\ell!}
\partial_x^{\ell+i}h(0,y,v) + \int_0^x\frac{(x-t)^{n-i-j}}{(n-i-j)!}
\partial_x^{n+1-j} h(t,y,v) dt\\ & = &
\int_0^x\frac{(x-t)^{n-i-j}}{(n-i-j)!}
\partial_x^{n+1-j} h(t,y,v) dt +
\left\{%
\begin{array}{ll}
    g(y,v), & \hbox{$i=0$;} \\
    0, & \hbox{otherwise.} \\
\end{array}%
\right.   \eeaa This gives
\bean |I| &=& \int_0^x\frac{(x-t)^{n-i-j}}{(n-i-j)!}
\underbrace{|\partial_y^j
\partial_x^{n+1-j}
\partial_v^\gamma h(t,y,v)|}_{\le C y^{1-\sigma-\lambda} |\ln^N t|
(y^{n+\lambda-n-1} + t ^{n+\lambda-n-1})\le 2Cy^{1-\sigma-\lambda}
t^{\lambda-1}|\ln^N t|}dt
\\ &\le & C' x^{n+\lambda-i-j}y^{1-\sigma-\lambda}|\ln^{N+1} x|
\;. \eeal{Iest} This calculation also proves that
$$h\in
y^{1-\sigma-\lambda}\fzTi {n+\lambda}\;.$$

For $ i+j> n$ we write
\beaa |\partial_x^i\partial_y^j
\partial_v^\gamma \hat h |&=& \Big|\partial_x^i\partial_y^j
\partial_v^\gamma \Big(\int_x^y f(s, \cdots) ds - h\Big) \Big|
\\
&\le & \underbrace{\Big|\partial_x^i\partial_y^j
 \Big(\int_x^y \partial_v^\gamma f(s, \cdots) ds\Big)\Big|}_{A} +
 \underbrace{\Big|\partial_x^i \partial_y^j\partial_v^\gamma \hat h \Big|}_{B}\;.
\\ && \eeaa
The second term is immediately estimated by a constant times
$x^{n+\lambda-i-j}y^{1-\sigma-\lambda}|\ln^N x|$. The estimation
of $A$ is simplest when $i>0$ as then we have \bel{Aest}|A| =
\Big|\partial_x^{i-1}\partial_y^j
  \partial_v^\gamma f(x, y,v) \Big|\le C
\Big(x^{\alpha+1+\beta-i-j} + x^{\alpha+1-i} y^{\beta-j} +
x^{\alpha+1 +\beta-i-k}y^{k-j}\Big)|\ln^Nx|\;.\ee When $i=0$ we
use instead \beaa|A| &= & \Big| \underbrace{\sum_{\ell+m=j-1}
C(\ell,m)
\partial_x^\ell \partial_y^m
\partial^\gamma_v f (y,y,v) }_{\le C y^{\alpha+1+\beta-j}|\ln^Ny|\le C x^{\alpha+1+\beta-j}|\ln^Nx|}+ \int_x^y
\partial^j_y \partial^\gamma_v  f(x,y,v)\Big|\;, \eeaa
which is again estimated as in \eq{Aest}. Summarising, we have
obtained
$$ |\partial_x^i\partial_y^j
\partial_v^\gamma \hat h | \le C
\Big(x^{\alpha+1+\beta-i-j} + x^{\alpha+1-i} y^{\beta-j} +
x^{\alpha+1 +\beta-i-k}y^{k-j}+
x^{\alpha+1+\beta+\lambda+\sigma-i-j}y^{1-\sigma-\lambda}\Big)|\ln^Nx|\;,
$$
for any $0\le \sigma+\lambda<1$. Setting $\lambda=
1-\sigma-\epsilon$ and estimating  $y^{\epsilon}\le y_0^\epsilon$,
$x^\alpha\le y_0^\epsilon x^{\alpha-\epsilon}$, the result
follows.
 \qed

\begin{Proposition}  \label{PIoneF}Let
$\alpha+\localkdelta>-1$. For any $\epsilon>0$ we have
\begin{enumerate}\item$I_1(x^{\localkdelta }\ln^\ell x \fzTi{\alpha})\subset
y^\epsilon\fzTi {\alpha+\localkdelta +1-\epsilon}+\fAx{x\fzTi
{\alpha}}\;.$
\item $I_1(\fAx{x^{\localkdelta}\fzTi {\alpha}})\subset
y^\epsilon\fzTi
{\alpha+\localkdelta+1-\epsilon}+\fAx{x^{\localkdelta+1}\fzTi
{\alpha}}\;.$
\end{enumerate}
\end{Proposition}

\proof The proof of point 1. is essentially the same as that of
Proposition~\ref{PIoneT}, once the integration by parts of
\eq{intpa} has been done; we outline the main steps for
completeness. We write again $\alpha+\localkdelta =n-\sigma$,
$n\in\N$, $\sigma\in [0,1)$.

Let $f\in\fzTi{\alpha}$. First we integrate by parts, obtaining
\bel{intpa} I_1(x^{\localkdelta }\ln^\ell x f) = \sum_{r=0}^N \left. C(r)
s^{\localkdelta +1} \ln^r s f(s,y,v)\right|^{s=y}_{s=x} -
\sum_{r=0}^N \int_x^y C(r) s^{\localkdelta +1} \ln^r s \partial_x
f(s,y,v) ds\;.\ee (The heuristics behind the integration by parts
is that $\partial_x f$ behaves differently from $f$ at the
boundary, in particular $\partial_xf$ vanishes at the boundary for
$\alpha>1$, while $f$ itself does not.) We set
$$g(y,v)=-\sum_{r=0}^N \int_0^y C(r)
s^{\localkdelta +1} \ln^r s \partial_x f(s,y,v) ds\;.$$ A
repetition of the argument in the proof of \ref{PIoneT} shows that
the function $ g$ satisfies the hypotheses needed to construct the
extensions of Section~\ref{Sext}, with $m$ in Lemma~\ref{LExt2}
equal to $n$, $\localmu=1-\sigma$ and any $0\le \lambda
<1-\sigma$. We set
$$g_0= g\;,\quad g_i=0 \ \mbox{ for } \
1\le i \le n \;.$$ By Lemma~\ref{LExt2} there exists  a function
$h$, belonging to the space $ y^{1-\sigma-\lambda}\fTi
{n+\lambda}$ for any $0\le \lambda+\sigma<1$, such that
$\partial_x^ih|_{x=0} = g_i$. The Taylor expansion as in the
previous proof proves further that $$h\in
y^{1-\sigma-\lambda}\fzTi {n+\lambda}\;.$$

Let us now define
$$\hat h =  I_1(x^{\localkdelta }\ln^\ell x f)-h - \sum_{r=0}^N \left. C(r) s^{\localkdelta +1}
\ln^r s f(s,y,v)\right|^{s=y}_{s=x}\;.$$ An argument as in the
proof of Proposition~\ref{PIoneT} gives $\hat h \in
y^\epsilon\fzTi {\alpha+\localkdelta +1-\epsilon}$.

What remains are the terms $$\sum_{r=0}^N \left. C(r)
s^{\localkdelta +1} \ln^r s f(s,y,v)\right|^{s=y}_{s=x}\;.$$ The
part with $s=y$ is in $\fzTi {\alpha+\localkdelta +1}\subset
y^\epsilon \fzTi {\alpha+\localkdelta +1-\epsilon}$ and the part
with $s=x$ belongs to $\fAx{x\fzTi {\alpha}}$.

Point 2 follows immediately from point 1 using the
expansion~\eq{fAF}.
\qed

\begin{Proposition}
\label{PItwo}
\begin{enumerate}\item\label{PITwoF1}$I_2(\fzTi{\alpha})\subset \fzTi{\alpha}\;.$
\item \label{PITwoF2}$I_2(\fAx{x^n\fzTi {\alpha}})\subset \fAx{x^n\fzTi {\alpha}}\;.$
\item \label{PITwoT}$I_2(\fffT \alpha \beta k) \subset \fffT \alpha {\beta+1}{k+1}\;.$
\end{enumerate}
\end{Proposition}

\proof \ref{PITwoF1}. For $f\in \fzTi{\alpha}$ define $F := I_2
(f)$. For $j\ge1$ we have $$\partial_x^i
\partial_y^j \partial_v^\gamma F(x,y,v) = \partial_x^i
\partial_y^{j-1} \partial_v^\gamma f(x,y,v)\;,$$ and the estimates
are obvious. For $j=0$ we use
$$|\partial_x^i \partial_v^\gamma F| \le \underbrace{\sum_{\ell+m=i-1}
C(\ell,m)\left|
\partial_x^\ell \partial_y^m
\partial^\gamma_v f (x,x,v)\right|}_{\le C x^{\alpha+1-i}|\ln^N x|}
+\underbrace{\int_x^y |\partial_x^i
\partial_v^\gamma f(x,s,v)|ds}_{\le C^\prime x^{\alpha-i}y|\ln^N
x|}\;.$$ In the first term we estimate $x^{\alpha+1-i}\le x
^{\alpha-i} y_0$, in the second $y\le y_0$, and the proof is done.

 \ref{PITwoF2}. The result follows from
\ref{PITwoF1} and the definition of $\fAx{F}$ (see \Eq{fAF}).

\ref{PITwoT}. For $f\in \fffT\alpha\beta k$ let us define $F :=
I_2 (f)$. We have to estimate $|\partial_x^i
\partial_y^j \partial_v^\gamma F |$. This is simplest for $j\ge
1$, we have then $$|\partial_x^i
\partial_y^j \partial_v^\gamma F | = |\partial_x^i
\partial_y^{j-1} \partial_v^\gamma f |\;,$$
 and the appropriate
estimates with $\beta\to \beta+1$, $k\to k+1$ follow from the
definition of $\fffT\alpha\beta k$. Again for $j=0$ we use
$$|\partial_x^i \partial_v^\gamma F| \le \underbrace{\sum_{\ell+m=i-1}
C(\ell,m)\left|
\partial_x^\ell \partial_y^m
\partial^\gamma_v f (x,x,v)\right|}_{\le C x^{\alpha+\beta+1-i }|\ln^N x|}
+ 
\int_x^y |\partial_x^i
\partial_v^\gamma f(x,s,v)|ds
\;.$$ The second term is estimated by
$$  C
\Big(x^{\alpha+\beta-i} y+ x^{\alpha-i} y^{\beta+1} + x^{\alpha
+\beta+1-i-k}y^{k+1}+ x^{\alpha+\beta+1-i}\Big)|\ln^{N+1}x|\;.
$$
 The calculation
$$x^{\alpha+\beta-i}y|\ln^{N+1} x| = x^{\alpha+(\beta+1)-(k+1)-i}x^{k}y|\ln^{N+1}
x|\le x^{\alpha+(\beta+1)-(k+1)-i}y^{k+1}|\ln^{N+1} x|$$ ends the
proof. \qed

\section{Polyhomogeneity of solutions for a class of linear symmetric hyperbolic
systems with smooth coefficients} In this appendix we give a
simple proof  of a counterpart of \chcite{Theorem~3.4} without
corner conditions.
The $\mcA^\delta_k$, $\mcC^\alpha_{l|0}$ and $\mcC^\alpha$ spaces
are defined as in \chcite{Sections~3.3,~3.4 and (A.2)}. We work in
the $(x, v^A, \tau)$ coordinates of \chl, see Figure~\ref{F1}. Let
us define the following space: \label{SAB}
\begin{eqnarray*}
\lefteqn{\mcTaB_k = \Bigg\{f\;:\;\forall\; 0\leq i+j+|\gamma| \leq
k } &&\\ && |\partial^i_x\partial^j_\tau\partial^\gamma_v f| \leq
\cases{C(x+2\tau)^{\alpha-i-j} & \mbox{if } $\alpha-i-j\geq0$ \cr
Cx^{\alpha-i-j} & \mbox{if } $\alpha-i-j<0$}\,\Bigg\}\;.
\end{eqnarray*}
This space is very similar, \emph{but not identical}, to the space
defined in \eq{Tadef}, because of the log factor occurring in
\eq{Tadef}; that factor is not needed in the current section.
Next, the space $\Cpdk$ is defined as the space of functions which
can be written in the form
\[\sum_{i=0}^{k}\sum_{j=0}^{N_i} f_{ij}(x+2\tau)^{i\delta} \ln^j(x+2\tau)
+ f_{k\delta+\epsilon}\;,\] for some $\epsilon > 0$, with
functions $f_{ij} \in C_\infty(\overline{\Omega})$ and
$f_{k\delta+\epsilon} \in \mcTB^{k\delta+\epsilon}_\infty$. We
will always assume $1/\delta \in \N$.

We note the following properties:
\begin{Proposition}
\label{Tprop}
\begin{enumerate} \item\label{Tprop1} $\mcTB^\alpha_k \subset \mcTB^\beta_k$ for
$\alpha\geq\beta$, \item\label{Tprop2}$\mcC^\alpha_{l|0} \subset
\mcTaB_l$,
\item\label{Tprop3} set $\hat f (x,v^A, \tau) := f(x+2\tau, v^A)$,
then for $\alpha\notin\calk$ we have $\hat f \in \mcTinfB^\alpha$
for $f \in \mcC^\alpha_\infty(\{\tau=0\})$ and $\hat f \in \Cpdk$
for $f \in \mcA^\delta_k(\{\tau=0\})$, \item\label{Tprop4} $f\in
\mcTinfB^\alpha \Rightarrow xf\in \mcTinfB^{\alpha+1}$,
\item\label{Tprop5} $f\in \mcTinfB^\alpha \Rightarrow
\partial_x f,
\partial_\tau f \in \mcTinfB^{\alpha-1}, \;
\partial_v f \in \mcTinfB^\alpha$.

\end{enumerate}
\end{Proposition}
\proof \ref{Tprop1}: This follows from the definition of
$\mcTinfB^\alpha$ and from the fact that
$(x+2\tau)^{\alpha-\beta}\leq C$.

\ref{Tprop2}: $f\in\mcC^\alpha_{l|0} $ iff for all $m, n, \gamma$
satisfying $0\leq m+n+|\gamma|\leq l$ we have $|\partial_x^m
\partial_\tau^n \partial_v^\gamma f| \leq C x^{\alpha-m-n}$.
For $\alpha-m-n<0$ this is exactly the estimate we need and for
$\alpha-m-n\geq0$ we have
$x^{\alpha-m-n}\leq(x+2\tau)^{\alpha-m-n}$.

\ref{Tprop3}: Since $f\in\mcC^\alpha({\tau=0})$ we have
$|\partial_x^m \partial_v^\gamma f(x,v^A)| \leq C x^{\alpha-m}$.
The following holds: $|\partial_x^m
\partial_\tau^n\partial_v^\gamma \hat{f}(x,v^A,\tau)| =
2^n|\partial_x^{m+n}\partial_v^\gamma f(x+2\tau,v^A)| \leq C
(x+2\tau)^{\alpha-m-n}$, which together with $(x+2\tau)^\beta\leq
x^\beta$ for $\beta<0$ yields the result. The statement concerning
$\mcA$  follows from the application of the result to the error term
$f_{k\delta+\eps}\in\mcC^{k\delta+\epsilon}_\infty$.

\ref{Tprop4}: This follows from $\partial_x^m \partial_\tau^n
\partial_v^\gamma (xf) = m \partial_x^{m-1} \partial_\tau^n
\partial_v^\gamma f + x \partial_x^m \partial_\tau^n
\partial_v^\gamma f$.

\ref{Tprop5}: Obvious. \qed

\remark Let $\gamma\in\R$ and $f\in\mcTinfB^\alpha$. It may happen
that $x^\gamma f \notin \mcTinfB^{\alpha+\gamma}$ or
$(x+2\tau)^\gamma f \notin \mcTinfB^{\alpha+\gamma}$.

The behavior of $\mcTaB$ spaces under integration is summarised in
the following proposition:
\begin{Proposition}\label{calkt}
Let $f\in\mcTinfB^\alpha$, $\alpha\in\R\setminus\calk$. Then
\[I_1(f) := \int_x^{x+2\tau}f(s,v^A,\frac12 x+\tau-\frac12s) \,ds \;\in\; \mcTinfB^{\alpha+1}\;,\]
\[I_2(f) := \int_0^{\tau}f(x,v^A,s) \,ds \;\in\; \mcTinfB^{\alpha}\;.\]
\end{Proposition}
\proof First we note that
$$ \partial_x I_1(f) = I_1(\partial_x f)\;,\qquad
\partial_\tau I_1(f) = f(x+2\tau,v^A,0) +
I_1(\partial_\tau f)\;,$$ and for the latter integral
$$ \partial_x I_2(f) = I_2(\partial_x f)\;,\qquad
\partial_\tau I_2(f) = f(x,v^A,\tau).$$ Then the result is easily
obtained, with the following formula being the key estimate:
$$I_1(x^\alpha) =
\frac{(x+2\tau)^{\alpha+1}}{\alpha+1} -
\frac{x^{\alpha+1}}{\alpha+1} \leq \left\{\begin{array}{ll}
  C(x+2\tau)^{\alpha+1} & \alpha+1>0 \\
  Cx^{\alpha+1}& \alpha+1<0\\
\end{array}
\right. \;.$$\qed

 It is obvious that the
standard space of polyhomogeneous functions $\mcA^\delta_k$ is
closed with respect to $x\partial_x$, the same is true for
$x^\beta\mcA^\delta_k$. The $\Cpdk$ space does not have this
property (terms like $\frac{x}{x+2\tau}$ may appear). In order to
fix this we introduce another space, denoted by $\eCp$, consisting
of functions which can be written as a finite sum of terms of the
form
\[f_{ij\ell}x^i(x+2\tau)^{j\delta} \ln^\ell(x+2\tau)\;,\]
with $i,\ell\in\N$, $j\in\calk$, $0\leq i+j\delta \leq k\delta$,
plus an error term
$f_{k\delta+\epsilon}\in\mcTinfB^{k\delta+\epsilon}$. The space
$\eCp$ as well as $(x+2\tau)^\beta\eCp$ is closed with respect to
$x\partial_x, x\partial_\tau$ and $\partial_v$.

We use the symbols $\zmcAk^\delta$, $\zCpdk$ and $\zeCp$ for spaces
of function which can be written as sums of a finite number of terms
as above \emph{without} the error term.

If we restrict ourselves to one value of $\tau \geq 0$ then there
is a correspondence between the spaces $\mcA^\delta_k$ and $\eCp$.
Obviously any function belonging to $\mcA^\delta_k$ belongs also
to $\zmcAk^\delta + \Cpdk$ and $\zmcAk^\delta + \eCp$. The
relation in the opposite direction is given by the following
Lemma:
\begin{Lemma}
Let $\mathring\tau\geq0$ and let $f\in(\mcA^\delta_k +
\eCp)|_{\tau=\mathring\tau}$. Then $f\in \mcA^\delta_{k}$.
\end{Lemma}
\proof The ``error term' $r\in \mcTinfB^{k\delta+\eps}$ in $f$
equals
$$ f-\sum_{i,j,\ell} f_{ij\ell}x^i(x+2\mathring\tau)^{j\delta}
\ln^\ell(x+2\mathring\tau) -\sum_{i,j} f_{ij}x^{i\delta}\ln^j
x\;,$$ with both sums being finite. From the definition of our
spaces we have
\be|\partial_x^n \partial_v^\alpha r| \leq C \cdot \cases{
(x+2\mathring\tau)^{k\delta+\eps-n} & \mbox{if }\
$k\delta+\eps-n\geq 0$,\cr x^{k\delta+\eps-n} & \mbox{if }\
$k\delta+\eps-n <0$.}\ee We want to show that there exist $\hat
f_{ij}$ such that
\be|\partial_x^n \partial_v^\alpha
(f-\sum_{ij}\hat f_{ij}x^{i\delta}\ln^j x)| \leq C \cdot
x^{k\delta+\eps-n} \;.\ee The estimate for $k\delta+\eps-n<0$ is
trivial so we focus on the first case. Let us notice that
$k\delta+\eps = m +\alpha$ with $m=\lfloor k\delta +\eps\rfloor$
and $\alpha\in [0,1)$. Now, $r$ has $m$ $x$-derivatives
$\partial_x^n \partial_v^\alpha r $ bounded by a constant and the
$(m+1)$-th $x$-derivative $\partial_x^{m+1} \partial_v^\alpha r $
is bounded by $C\cdot x^{\alpha-1}$. Therefore we can
Taylor-expand $r$ up to the order $m$. The rest is of order
$m+\alpha$: \be |r(x,v) - \sum_{i=0}^m \frac{\partial_x^i
r(x=0,v)}{i!}x^i| \leq C\cdot x^{m+\alpha}\;\ee  and the proof is
easily completed.\qed

The following proposition gives the behaviour of polyhomogeneous
spaces under integration:
\begin{Proposition}\label{calkAo}
\begin{enumerate}
\item  Let $f\in(x+2\tau)^\beta\zeCp$. There exists
$f_{\beta+k\delta+\epsilon}\in\mcTinfB^{\beta+k\delta+\epsilon}$,
for some $\epsilon>0$,  such that
$$ I_2(f)  \;\in\;
(x+2\tau)^\beta\zeCp + x^\beta\zmcAk^\delta +
f_{\beta+k\delta+\epsilon};,$$
$$ I_1(f) \;\in\;
(x+2\tau)^\beta\zeCp +f_{\beta+k\delta+\epsilon}\;.$$ \item Let
$f\in x^\beta\zmcAk^\delta$. There exists
$f_{\beta+k\delta+\epsilon}\in\mcTinfB^{\beta+k\delta+\epsilon}$
such that
$$ I_1(f) \;\in\;
(x+2\tau)^\beta\zeCp + x^\beta\zmcAk^\delta
+f_{\beta+k\delta+\epsilon}\;.$$
\end{enumerate}
\end{Proposition}
\proof The proposition readily follows from the next two
lemmata.\qed

\begin{Lemma}
Let $f$ be a smooth function, $p\in\R, k\in\N$. For all
$\alpha\in\R$ there exists $N\in\N$, sequences of numbers
$A_{i}\in p+1+\N, B_{i}\in\N$ ,a sequence of smooth functions
$f_i$ and a function $f_\alpha\in\mcTinfB^{\alpha}$ such that
$$I_1(f x^p \ln^k x) =
\sum_{i=1}^N f_i \left( (x+2\tau)^{A_i} \ln^{B_i} (x+2\tau) -
x^{A_i} \ln^{B_i}x\right) + f_\alpha\;.$$
\end{Lemma}
\proof Integration by parts yields the result for $f=1$. The
result for general $f$ is also obtained by integration by parts:
$$\int f v^\prime = fv - \int f^\prime v\;,$$
where $v^\prime = s^p \ln^k s$. Using the result for $f=1$ one
gets that $v$ has a power of $s$ one higher than $v^\prime$.
Repeating this integration a finite number of times yields a
result with an error term in $\mcTinfB^n$ and $n$ high enough.
\qed

\begin{Lemma}
Let $f$ be smooth, $i\geq -1$ and $\alpha>0$. Then there exists
$n\in \N$, a sequence of smooth functions $f_j$ and $f_\alpha \in
\mcTinfB^{\alpha}$ such that
$$ I_1(f x^i) = \sum_{j=0}^n f_j \left(
(x+2\tau)^{i+j+1} - x^{i+j+1} \right) + f_\alpha  \;.$$
\end{Lemma}
\proof Since $f$ is smooth it can be expanded in powers of $x$ to
any order and the result follows. \qed

\begin{Proposition}\label{calkbphi}
Let $\varphi$ be a solution of $$\partial_\tau\varphi + b\varphi =
c\;,$$ with $$b\in C_\infty(\overline{\Omega})\;.$$ Let us assume
that $$\varphi(x,v,0) \in x^\beta\mcA^\delta_k\;$$ and
$$c\in(x+2\tau)^\beta\zeCp + x^\beta\zmcAk^\delta + \mcTinfB^\alpha +
C_\infty(\overline{\Omega})\;.$$ Then
$$\varphi\in(x+2\tau)^\beta\zeCp + x^\beta\zmcAk^\delta
+\mcTinfB^{\min(\alpha,
\beta+k\delta+1)}+ C_\infty(\overline{\Omega})\;.$$
\end{Proposition}
\proof The solution of the equation at hand may be expressed as
$$\varphi(\cdot,\tau) = R(\cdot,\tau)\varphi(\cdot,0) +
\int_0^\tau R(\cdot,s)c(\cdot,s) \,ds\;,$$ where $R(x,v,\tau)$ is
the resolvent of the equation $\partial_\tau\varphi = -b\varphi$.
It is a standard result that for $b\in
C_\infty(\overline{\Omega})$ we have also $R\in
C_\infty(\overline{\Omega})$. Then the result follows from the
analysis of the above formula and Propositions \ref{calkt} and
\ref{calkAo}.
 \qed

Now we are ready to pass to the proof of polyhomogeneity of
solutions:

\begin{Theorem}
Let $\alpha, \beta \in \R$, $k \in \N$, and let $(\varphi, \psi) \in
\mcC^{\alpha}_{\infty}$ be a solution of
\begin{eqnarray*}
  \partial_\tau\varphi + B_{11}\varphi + B_{12}\psi&=& L_{11}\varphi +
  L_{12}\psi + a
  \;,\\ e_+\psi+
  B_{21}\varphi + B_{22}\psi &=&   L_{21}\varphi +
  L_{22}\psi + b \label{s1bcphg}\;,
\end{eqnarray*}
where $e_+ = \partial_\tau - 2\partial_x$. Suppose that
 \be L_{ij}=
L_{ij}^{A}\partial_A + xL_{ij}^\tau
\partial_\tau + xL_{ij}^x\partial_x\,,
\ee with \be
  L_{11}^{\mu}\in x C_\infty(\overline{\Omega})\;, \quad L_{21}^{\mu}\;,
  L_{12}^{\mu}\;, L_{22}^{\mu}\in C_\infty(\overline{\Omega}) \;,
\ee and that
\begin{deqarr}&
    B_{ab}\in C_\infty(\overline{\Omega})\;,&
\\ &a,b \in C_\infty(\overline{\Omega})\;,
\qquad\varphi(0), \psi(0) \in x^\beta\mcA_k^\delta(M_{x_0})\;.&
  \end{deqarr}
Then
$$\varphi, \psi \in
x^\beta\zmcAk^\delta+(x+2\tau)^\beta\zeCp
+C_\infty(\overline{\Omega}) + \mcTinfB^{\beta+k\delta+\epsilon}$$
for some $\epsilon>0$.
\end{Theorem}
\proof The proof is very similar to the proof of Theorem
\chcite{Theorem~3.4}. First we notice that $$(\varphi, \psi) \in
\mcTinfB^{\alpha}\;,$$ which is due to Proposition \ref{Tprop}
point \ref{Tprop2}. For the purpose of the proof it will be
convenient to use the following notation:
$$ \soln {\alpha} := \mcTinfB^{\alpha} + C_\infty + x^\beta\zmcAk^\delta +
(x+2\tau)^\beta \zeCp\;.$$ To prove the theorem we need to show
\be\label{solnest}(\varphi,\psi) \in
\soln{\beta+k\delta+\epsilon}\;.\ee We rewrite the equations at
hand as
\begin{eqnarray*}
\partial_\tau \varphi  + B_{11}\varphi& = & c_1\;,
\\
e_+ \psi  & = & c_2\;,
\end{eqnarray*}
where
\begin{eqnarray*}
c_1&:=& L_{11}\varphi + L_{12}\psi + a -B_{12}\psi \;,
\\c_2 &:=&
L_{21}\varphi + L_{22}\psi + b-B_{21}\varphi -B_{22}\psi\;.
\end{eqnarray*}
Let us start with the second of the equations. The integration
yields
$$
 \psi(x,v^A,\tau) = \psi(x+2\tau,v^A,0) +
 \frac 12
I_1(c_2)\;.
$$
 We have
$$c_2\in \mcTinfB^{\alpha} + C_\infty \subset
\soln{\alpha}\;.$$ Propositions \ref{calkt} and \ref{calkAo}
together with \ref{Tprop} point \ref{Tprop1} yield
$$ I_1(c_2) \in \soln{\min(\alpha+1,\beta+k\delta+\epsilon)}\;.$$
{}From Proposition \ref{Tprop} point \ref{Tprop3} we have
$$ \psi(x+2\tau,v^A,0) \in \soln{\beta+k\delta+\epsilon}\;.$$
Therefore $$\psi \in
\soln{\min(\alpha+1,\beta+k\delta+\epsilon)}\;.$$ Now we estimate
$c_1$ to be in $\soln{\min(\alpha+1,\beta+k\delta+\epsilon)} $ and
use Proposition \ref{calkbphi} to get
$$\varphi \in \soln{\min(\alpha+1,\beta+k\delta+\epsilon)}\;.$$
Then we repeat this procedure, with $c_2$ now in
$\soln{\min(\alpha+1,\beta+k\delta+\epsilon)}$ and get first
$\psi$ then $\varphi$ in the appropriate $\soln{}$ space with the
error-term index increased by one. After a finite number of steps
we get \eq{solnest}. \qed

\remark If we make an additional assumption that $L^\tau_{ij} =
L^x_{ij}=0$, we obtaint the following variant of the result:
$$\varphi, \psi \in x^\beta\zmcAk^\delta+(x+2\tau)^\beta\zCpdk
+C_\infty(\overline{\Omega}) +
\mcTinfB^{\beta+k\delta+\epsilon}\;.$$ The $\zeCp$ space has been
replaced by $\zCpdk$.

\def\cprime{$'$} \def\cprime{$'$} \def\cprime{$'$}
\providecommand{\bysame}{\leavevmode\hbox
to3em{\hrulefill}\thinspace}
\providecommand{\MR}{\relax\ifhmode\unskip\space\fi MR }
\providecommand{\MRhref}[2]{%
  \href{http://www.ams.org/mathscinet-getitem?mr=#1}{#2}
} \providecommand{\href}[2]{#2}

%



\end{document}